\documentclass[a4paper,11pt]{amsart}
\usepackage{mathrsfs}
\usepackage{cases}
\usepackage{epic}
\usepackage{amsfonts}
\usepackage{graphicx}
\usepackage{amsmath}
\usepackage{amssymb, upgreek}
\usepackage{bm}
\usepackage{latexsym,todonotes}
\usepackage{pdflscape}
\usepackage[all]{xypic}
\usepackage[all]{xy}
\usepackage{color}
\usepackage{colordvi}
\usepackage{multicol}
\usepackage[linktocpage=true]{hyperref}
\hypersetup{colorlinks,linkcolor=blue,urlcolor=cyan,citecolor=blue}
\newcommand{\bvs}{\mathbf{\varsigma}}
\newcommand{\vs}{\varsigma}
\allowdisplaybreaks

\def \cI {\mathcal I}

\def \cJ{\mathcal J}

\begin{document}
	\input xy
	\xyoption{all}
	
	\def \vth{{\theta}}
	\newtheorem{innercustomthm}{{\bf Main~Theorem}}
	\newenvironment{customthm}[1]
	{\renewcommand\theinnercustomthm{#1}\innercustomthm}
	{\endinnercustomthm}
	
	\newtheorem{innercustomcor}{{\bf Corollary}}
	\newenvironment{customcor}[1]
	{\renewcommand\theinnercustomcor{#1}\innercustomcor}
	{\endinnercustomthm}
	
	\newtheorem{innercustomprop}{{\bf Proposition}}
	\newenvironment{customprop}[1]
	{\renewcommand\theinnercustomprop{#1}\innercustomprop}
	{\endinnercustomthm}
	
	\newcommand{\LaK}{\Lambda_{\texttt{Kr}}}
	\def \BH{\mathbb{H}}
	\newcommand{\oH}{\ov{\BH}}
	
	\newcommand{\iadd}{\operatorname{iadd}\nolimits}
	\newcommand{\Gr}{\operatorname{Gr}\nolimits}
	\newcommand{\FGS}{\operatorname{FGS}\nolimits}
	\newcommand{\tUiD}{\operatorname{{}^{{Dr}}\tUi}\nolimits}
	\renewcommand{\mod}{\operatorname{mod}\nolimits}
	\newcommand{\proj}{\operatorname{proj}\nolimits}
	\newcommand{\inj}{\operatorname{inj.}\nolimits}
	\newcommand{\rad}{\operatorname{rad}\nolimits}
	\newcommand{\Span}{\operatorname{Span}\nolimits}
	\newcommand{\soc}{\operatorname{soc}\nolimits}
	\newcommand{\ind}{\operatorname{inj.dim}\nolimits}
	\newcommand{\Ginj}{\operatorname{Ginj}\nolimits}
	\newcommand{\res}{\operatorname{res}\nolimits}
	\newcommand{\np}{\operatorname{np}\nolimits}
	\newcommand{\Fac}{\operatorname{Fac}\nolimits}
	\newcommand{\Aut}{\operatorname{Aut}\nolimits}
	\newcommand{\DTr}{\operatorname{DTr}\nolimits}
	\newcommand{\TrD}{\operatorname{TrD}\nolimits}
	
	\newcommand{\Mod}{\operatorname{Mod}\nolimits}
	\newcommand{\R}{\operatorname{R}\nolimits}
	\newcommand{\End}{\operatorname{End}\nolimits}
	\newcommand{\lf}{\operatorname{l.f.}\nolimits}
	\newcommand{\Iso}{\operatorname{Iso}\nolimits}
	\newcommand{\aut}{\operatorname{Aut}\nolimits}
	\newcommand{\Ui}{{\mathbf U}^\imath}
	\newcommand{\UU}{{\mathbf U}\otimes {\mathbf U}}
	\newcommand{\UUi}{(\UU)^\imath}
	\newcommand{\tUU}{{\tU}\otimes {\tU}}
	\newcommand{\tUUi}{(\tUU)^\imath}
	\newcommand{\tUi}{\widetilde{{\mathbf U}}^\imath}
	\newcommand{\sqq}{{\bf v}}
	\newcommand{\sqvs}{\sqrt{\vs}}
	\newcommand{\dbl}{\operatorname{dbl}\nolimits}
	\newcommand{\swa}{\operatorname{swap}\nolimits}
	\newcommand{\Gp}{\operatorname{Gp}\nolimits}
	
	\newcommand{\U}{{\mathbf U}}
	\newcommand{\tU}{\widetilde{\mathbf U}}
	\newcommand{\fgm}{{\rm mod}^{{\rm fg}}}
	\newcommand{\fgmz}{\mathrm{mod}^{{\rm fg},\Z}}
	\newcommand{\fdmz}{\mathrm{mod}^{\Z}}
	\newcommand{\tMHL}{{}^\imath\widetilde{\ch}(\bfk \QK)}
	\newcommand{\ov}{\overline}
	\newcommand{\und}{\underline}
	\newcommand{\tk}{\widetilde{k}}
	\newcommand{\tH}{\widetilde{H}}
	\newcommand{\tK}{\widetilde{K}}
	\newcommand{\tTT}{\operatorname{\widetilde{\texttt{\rm T}}}\nolimits}
	
	\def \iQ{Q^\imath}

	\newcommand{\utM}{\operatorname{\cm\ch}\nolimits}
	\newcommand{\tM}{\operatorname{\cs\cd\widetilde{\ch}}\nolimits}
	\newcommand{\rM}{\operatorname{\cm\ch_{\rm{red}}}\nolimits}
	\newcommand{\utMH}{\cs\cd\ch(\Lambda^\imath)}
	\newcommand{\tMH}{\cs\cd\widetilde{\ch}(\Lambda^\imath)}
	\newcommand{\tCMH}{{\cc\widetilde{\ch}(\bfk Q,\btau)}}
	\newcommand{\QK}{Q_{{\texttt{Kr}}}}
	\newcommand{\rMH}{\operatorname{\cm\ch_{\rm{red}}(\Lambda^\imath)}\nolimits}
	\newcommand{\utMHg}{\operatorname{\ch(Q,\btau)}\nolimits}
	\newcommand{\tMHg}{\operatorname{\widetilde{\ch}(Q,\btau)}\nolimits}
	\newcommand{\tMHk}{{\widetilde{\ch}(\bfk Q,\btau)}}
	\newcommand{\rMHg}{\operatorname{\ch_{\rm{red}}(Q,\btau)}\nolimits}
	
	\newcommand{\rMHd}{\operatorname{\cm\ch_{\rm{red}}(\Lambda^\imath)_{\bvsd}}\nolimits}
	\newcommand{\tMHd}{\operatorname{\cs\cd\widetilde{\ch}(\Lambda^\imath)_{\bvsd}}\nolimits}
	
	\newcommand{\tMHl}{\cs\cd\widetilde{\ch}({\bs}_\ell\Lambda^\imath)}
	\newcommand{\rMHl}{\cm\ch_{\rm{red}}({\bs}_\ell\Lambda^\imath)_{\bvsd}}
	\newcommand{\tMHi}{\cs\cd\widetilde{\ch}({\bs}_i\Lambda^\imath)}
	\newcommand{\rMHi}{\cm\ch_{\rm{red}}({\bs}_i\Lambda^\imath)_{\bvsd}}
	\newcommand{\tMHgi}{\widetilde{\ch}({\bs}_i Q,\btau)}

	\newcommand{\utGpg}{\operatorname{\ch^{\rm Gp}(Q,\btau)}\nolimits}
	\newcommand{\tGpg}{\operatorname{\widetilde{\ch}^{\rm Gp}(Q,\btau)}\nolimits}
	\newcommand{\rGpg}{\operatorname{\ch_{red}^{\rm Gp}(Q,\btau)}\nolimits}
	
	\def \blx{x}
	\newcommand{\colim}{\operatorname{colim}\nolimits}
	\newcommand{\gldim}{\operatorname{gl.dim}\nolimits}
	\newcommand{\cone}{\operatorname{cone}\nolimits}
	\newcommand{\rep}{\operatorname{rep}\nolimits}
	\newcommand{\Ext}{\operatorname{Ext}\nolimits}
	\newcommand{\Tor}{\operatorname{Tor}\nolimits}
	\newcommand{\Hom}{\operatorname{Hom}\nolimits}
	\newcommand{\Top}{\operatorname{top}\nolimits}
	\newcommand{\Coker}{\operatorname{Coker}\nolimits}
	\newcommand{\thick}{\operatorname{thick}\nolimits}
	\newcommand{\rank}{\operatorname{rank}\nolimits}
	\newcommand{\Gproj}{\operatorname{Gproj}\nolimits}
	\newcommand{\Len}{\operatorname{Length}\nolimits}
	\newcommand{\RHom}{\operatorname{RHom}\nolimits}
	\renewcommand{\deg}{\operatorname{deg}\nolimits}
	\renewcommand{\Im}{\operatorname{Im}\nolimits}
	\newcommand{\Ker}{\operatorname{Ker}\nolimits}
	\newcommand{\Coh}{\operatorname{Coh}\nolimits}
	\newcommand{\Id}{\operatorname{Id}\nolimits}
	\newcommand{\Qcoh}{\operatorname{Qch}\nolimits}
	\newcommand{\CM}{\operatorname{CM}\nolimits}
	\newcommand{\sgn}{\operatorname{sgn}\nolimits}
	\newcommand{\Gdim}{\operatorname{G.dim}\nolimits}
	\newcommand{\fpr}{\operatorname{\mathcal{P}^{\leq1}}\nolimits}
	\newcommand{\QJ}{Q_{\texttt J}}
	\newcommand{\tMHLJ}{{}^\imath\widetilde{\ch}(\bfk \QJ)}	
	\newcommand{\For}{\operatorname{{\bf F}or}\nolimits}
	\newcommand{\coker}{\operatorname{Coker}\nolimits}
	\renewcommand{\dim}{\operatorname{dim}\nolimits}
	\newcommand{\rankv}{\operatorname{\underline{rank}}\nolimits}
	\newcommand{\dimv}{{\operatorname{\underline{dim}}\nolimits}}
	\newcommand{\diag}{{\operatorname{diag}\nolimits}}
	\newcommand{\qbinom}[2]{\begin{bmatrix} #1\\#2 \end{bmatrix} }
	
	\renewcommand{\Vec}{{\operatorname{Vec}\nolimits}}
	\newcommand{\pd}{\operatorname{proj.dim}\nolimits}
	\newcommand{\gr}{\operatorname{gr}\nolimits}
	\newcommand{\id}{\operatorname{Id}\nolimits}
	\newcommand{\Res}{\operatorname{Res}\nolimits}
	\def \tT{\widetilde{\mathcal T}}
	\def \tTL{\tT(\Lambda^\imath)}
	\newcommand{\iH}{{}^\imath\widetilde{\ch}}
	\newcommand{\mbf}{\mathbf}
	\newcommand{\mbb}{\mathbb}
	\newcommand{\mrm}{\mathrm}
	\newcommand{\cbinom}[2]{\left\{ \begin{matrix} #1\\#2 \end{matrix} \right\}}
	\newcommand{\dvev}[1]{{B_1|}_{\ev}^{{(#1)}}}
	\newcommand{\dv}[1]{{B_1|}_{\odd}^{{(#1)}}}
	\newcommand{\dvd}[1]{t_{\odd}^{{(#1)}}}
	\newcommand{\dvp}[1]{t_{\ev}^{{(#1)}}}
	\newcommand{\ev}{\bar{0}}
	\newcommand{\odd}{\bar{1}}
	\newcommand{\Iblack}{\I_{\bullet}}
	\newcommand{\wb}{w_\bullet}
	\newcommand{\Uidot}{\dot{\bold{U}}^{\imath}}
	\def \scrM{\mathscr M}
	\def \scrf{\mathscr F}
	\def \scrt{\mathscr T}
	\newcommand{\kk}{h}
	\newcommand{\la}{\lambda}
	\newcommand{\LR}[2]{\left\llbracket \begin{matrix} #1\\#2 \end{matrix} \right\rrbracket}
	\newcommand{\ff}{B}
	\newcommand{\pdim}{\operatorname{proj.dim}\nolimits}
	\newcommand{\idim}{\operatorname{inj.dim}\nolimits}
	\newcommand{\Gd}{\operatorname{G.dim}\nolimits}
	\newcommand{\Ind}{\operatorname{Ind}\nolimits}
	\newcommand{\add}{\operatorname{add}\nolimits}
	\newcommand{\pr}{\operatorname{pr}\nolimits}
	\newcommand{\oR}{\operatorname{R}\nolimits}
	\newcommand{\oL}{\operatorname{L}\nolimits}
	\newcommand{\ext}{{ \mathfrak{Ext}}}
	\newcommand{\Perf}{{\mathfrak Perf}}
	\def\scrP{\mathscr{P}}
	\def \bfk {k}
	\newcommand{\bk}{{\mathbb K}}
	\newcommand{\cc}{{\mathcal C}}
	\newcommand{\gc}{{\mathcal GC}}
	\newcommand{\dg}{{\rm dg}}
	\newcommand{\ce}{{\mathcal E}}
	\newcommand{\cs}{{\mathcal S}}
	\newcommand{\cl}{{\mathcal L}}
	\newcommand{\cf}{{\mathcal F}}
	\newcommand{\cx}{{\mathcal X}}
	\newcommand{\cy}{{\mathcal Y}}
	\newcommand{\ct}{{\mathcal T}}
	\newcommand{\cu}{{\mathcal U}}
	\newcommand{\cv}{{\mathcal V}}
	\newcommand{\cn}{{\mathcal N}}
	\newcommand{\mcr}{{\mathcal R}}
	\newcommand{\ch}{{\mathcal H}}
	\newcommand{\ca}{{\mathcal A}}
	\newcommand{\cb}{{\mathcal B}}
	\newcommand{\ci}{{\I}_{\btau}}
	\newcommand{\cj}{{\mathcal J}}
	\newcommand{\cm}{{\mathcal M}}
	\newcommand{\cp}{{\mathcal P}}
	\newcommand{\cg}{{\mathcal G}}
	\newcommand{\cw}{{\mathcal W}}
	\newcommand{\co}{{\mathcal O}}
	\newcommand{\cq}{{Q^{\rm dbl}}}
	\newcommand{\cd}{{\mathcal D}}
	\newcommand{\ck}{\widetilde{\mathcal K}}
	\newcommand{\calr}{{\mathcal R}}
	\newcommand{\iLa}{\Lambda^{\imath}}
	\newcommand{\La}{\Lambda}
	\newcommand{\ol}{\overline}
	\newcommand{\ul}{\underline}
	\newcommand{\st}{[1]}
	\newcommand{\ow}{\widetilde}
	\renewcommand{\P}{\mathbf{P}}
	\newcommand{\pic}{\operatorname{Pic}\nolimits}
	\newcommand{\Spec}{\operatorname{Spec}\nolimits}
	\newcommand{\wt}{\text{wt}\,}
	
	\newtheorem{theorem}{Theorem}[section]
	\newtheorem{acknowledgement}[theorem]{Acknowledgement}
	\newtheorem{algorithm}[theorem]{Algorithm}
	\newtheorem{assumption}[theorem]{Assumption}
	\newtheorem{axiom}[theorem]{Axiom}
	\newtheorem{case}[theorem]{Case}
	\newtheorem{claim}[theorem]{Claim}
	\newtheorem{conclusion}[theorem]{Conclusion}
	\newtheorem{condition}[theorem]{Condition}
	\newtheorem{conjecture}[theorem]{Conjecture}
	\newtheorem{construction}[theorem]{Construction}
	\newtheorem{corollary}[theorem]{Corollary}
	\newtheorem{criterion}[theorem]{Criterion}
	\newtheorem{definition}[theorem]{Definition}
	\newtheorem{example}[theorem]{Example}
	\newtheorem{exercise}[theorem]{Exercise}
	\newtheorem{lemma}[theorem]{Lemma}
	\newtheorem{notation}[theorem]{Notation}
	\newtheorem{problem}[theorem]{Problem}
	\newtheorem{proposition}[theorem]{Proposition}
	\newtheorem{solution}[theorem]{Solution}
	\newtheorem{summary}[theorem]{Summary}
	\numberwithin{equation}{section}
	
	\theoremstyle{remark}
	\newtheorem{remark}[theorem]{Remark}
	\newcommand{\Pd}{\pi_*}
	\def \bvs{{\boldsymbol{\varsigma}}}
	\def \bvsd{{\boldsymbol{\varsigma}_{\diamond}}}
	\def \btau{\tau}
	\newcommand{\LL}{\texttt{L}}
	\newcommand{\RR}{\texttt{R}}
	
	\def \hA{\widehat{\mathfrak{sl}}_2}
	\def \y{B}
	\def \haB{\widehat{B}}
	\def \bp{{\mathbf p}}
	\def \bq{{\bm q}}
	\def \bv{{v}}
	\def \bs{{\bm s}}
	\def \de {\delta}
	\def \bt {\mathbf{t}}
	\newcommand{\PL}{\bbP^1_{\bfk}}
	\def\bbP{{\mathbb P}}
	\def \bfK{{\mathbf K}}
	\def\coh{{\rm coh}}
	\def \bbZ{\mathbb Z}
	\def \BG{{\mathbb G}}
	\def \BS{{\mathbb S}}
	\def \BF{{\digamma}}
	\newcommand{\tCMHg}{\cc\widetilde{\ch}(Q,\btau)}
	\newcommand{\bfv}{\mathbf{v}}
	\def \bA{{\mathbf A}}
	\def \ba{{\mathbf a}}
	\def \bL{{\mathbf L}}
	\def \bF{{\mathbf F}}
	\def \bS{{\mathbf S}}
	\def \bC{{\mathbf C}}
	\def \bU{{\mathbf U}}
	\def \bc{{\mathbf c}}
	\def \fpi{\mathfrak{P}^\imath}
	\def \Ni{N^\imath}
	\def \fp{\mathfrak{P}}
	\def \fg{\mathfrak{g}}
	\def \fk{\fg^\theta}  
	
	\def \II{\I_0}
	\def \ul{\underline}
	\def \bla{\boldsymbol{\lambda}}
	\def \X{\mathbb{X}}
	
	\def \fn{\mathfrak{n}}
	\def \fh{\mathfrak{h}}
	\def \fu{\mathfrak{u}}
	\def \fv{\mathfrak{v}}
	\def \fa{\mathfrak{a}}
	\def \fq{\mathfrak{q}}
	\def \Z{{\Bbb Z}}
	\def \F{{\Bbb F}}
	\def \SS{\Bbb S}
	\def \D{{\Bbb D}}
	\def \R{{\Bbb R}}
	\def \C{{\Bbb C}}
	\def \N{{\Bbb N}}
	\def \Q{{\Bbb Q}}
	\def \G{{\Bbb G}}
	\def \P{{\Bbb P}}
	\def \K{{\Bbb K}}
	\def \bK{{\Bbb K}}
	\def \mfp{\mathfrak p}
	\newcommand{\haH}{\widehat{H}}
	\newcommand{\haT}{\widehat{\Theta}}
	\def \E{{\Bbb K}}
	\def \A{{\Bbb A}}
	\def \L{{\Bbb L}}
	\def \I{{\Bbb I}}
	\def \BH{{\Bbb H}}
	\def \T{{\Bbb T}}
	\def \bfk{\mathbf k}
	\def \TT{\bold{T}}
	\newcommand {\lu}[1]{\textcolor{red}{$\clubsuit$: #1}}
	\newcommand{\tMHX}{\operatorname{{}^\imath\widetilde{\ch}(\PL)}\nolimits}
	\newcommand{\tMHW}{\operatorname{{}^\imath\widetilde{\ch}(\X_\bfk)}\nolimits}
	\newcommand{\nc}{\newcommand}
	\newcommand{\browntext}[1]{\textcolor{brown}{#1}}
	\newcommand{\greentext}[1]{\textcolor{green}{#1}}
	\newcommand{\redtext}[1]{\textcolor{red}{#1}}
	\newcommand{\bluetext}[1]{\textcolor{blue}{#1}}
	\newcommand{\brown}[1]{\browntext{ #1}}
	\newcommand{\green}[1]{\greentext{ #1}}
	\newcommand{\red}[1]{\redtext{ #1}}
	\newcommand{\blue}[1]{\bluetext{ #1}}
	\newcommand{\arxiv}[1]{\href{http://arxiv.org/abs/#1}{\tt
			arXiv:\nolinkurl{#1}}}
	
	\newcommand{\wtodo}{\todo[inline,color=orange!20, caption={}]}
	\newcommand{\lutodo}{\todo[inline,color=green!20, caption={}]}

	\title[$\imath$Hall algebras and $\imath$quantum groups]{$\imath$Hall algebras and $\imath$quantum groups}
	
	\author[Ming Lu]{Ming Lu}
	\address{Department of Mathematics, Sichuan University, Chengdu 610064, P.R.China}
	\email{luming@scu.edu.cn}

	\author[Weiqiang Wang]{Weiqiang Wang}
	\address{Department of Mathematics, University of Virginia, Charlottesville, VA 22904, USA}
	\email{ww9c@virginia.edu}


	\maketitle
	
	\vspace{-1em}
	
	\begin{quote}\begin{center}
			{\em In memory of Brian Parshall}
		\end{center}
	\end{quote}   
	
	\begin{abstract}
		We survey some recent development on the theory of $\imath$Hall algebras. Starting from $\imath$quivers (aka quivers with involutions), we construct a class of 1-Gorenstein algebras called $\imath$quiver algebras, whose semi-derived Hall algebras give us $\imath$Hall algebras. We then use these $\imath$Hall algebras to realize quasi-split $\imath$quantum groups arising from quantum symmetric pairs. Relative braid group symmetries on $\imath$quantum groups are realized via reflection functors. In case of Jordan $\imath$quiver, the $\imath$Hall algebra is commutative and connections to $\imath$Hall-Littlewood symmetric functions are developed. In case of $\imath$quivers of diagonal type, our construction amounts to a reformulation of Bridgeland-Hall algebra realization of the Drinfeld double quantum groups (which in turn generalizes Ringel-Hall algebra realization of halves of quantum groups). Many rank 1 and rank 2 computations are supplied to illustrate the general constructions. We also briefly review $\imath$Hall algebras of weighted projective lines, and use them to realize Drinfeld type presentations of $\imath$quantum loop algebras. 
	\end{abstract}

	\setcounter{tocdepth}{2}
	\tableofcontents
	
	\section{Introduction}
	\subsection{Hall algebras}

	The theory of Hall algebras has a rich history and diverse applications. We give a quick overview of 3 (interconnected) classes of Hall algebras. 
	
	This first Hall algebra (which was popularized by Macdonald's book \cite{Mac95} but goes back to Steinitz decades earlier) is an associative algebra with a (Hall) basis parameterized by the isoclasses of finite modules over a discrete valuation ring (such as the ring of  $p$-adic integers or the power series ring $\bfk [[x]]$ over a finite field $\bfk$ of $q$ elements), or equivalently, parameterized by partitions. The structure constants arise from counting the extensions between these modules. In this case, the Hall algebra is commutative and isomorphic to the ring $\Lambda$ of symmetric functions. The Hall basis corresponds to a distinguished basis of $\Lambda$ known as Hall-Littlewood functions. From a modern point of view, the Steinitz-Hall algebra is associated to the Jordan quiver (consisting of a vertex and a loop). 
	
	Ringel-Hall algebra associated to a Dynkin quiver was introduced by \cite{Rin90}, whose structure constants count the extensions of quiver representations over $\bfk$. Ringel showed that this version of Hall algebra is isomorphic to half a quantum group associated to the underlying Dynkin diagram of the quiver. Lusztig \cite{Lus90, Lus93} developed a geometric version of Ringel-Hall algebra, which allows him to construct a canonical basis of half the quantum group with remarkable positivity property. Ringel's construction was extended by J.A.~Green \cite{Gr95} to acyclic quivers in connection to quantum groups of Kac-Moody type, who also provides a categorical interpretation of the comultiplication. Reflection functors on Hall algebras provide a realization of braid group symmetries on quantum groups, and in addition, they help to construct PBW bases through Hall algebra approach; see \cite{Lus93, Rin96, CX99, Z00, XY01}. 
	
	In spite of various attempts, the search for a Hall algebra realization of the whole quantum group $\U$ (or the Drinfeld double $\tU$) was largely unsuccessful until the work of Bridgeland \cite{Br} 20 years later. We mention some of these attempts. One first (non-categorical) approach was to use the Drinfeld double of the twisted extended Hall algebras from  hereditary categories to realize the whole quantum groups (see \cite{X97,SV99}). Also relevant is the (not entirely categorical either) construction of the  Hall Lie algebras from root categories (see \cite{PX00, Hub06}) to realize the Kac-Moody Lie algebras. Yet another construction is the derived Hall algebras from derived categories or more general triangulated categories with certain finiteness conditions (see \cite{T06,XX08}, for a different version see \cite{Ka98}).  Bridgeland's construction starts with the Hall algebra of the exact category of $\Z_2$-graded complexes of projective modules of quivers and then applies localizations to the acyclic complexes.  
	
	Before describing a third class of Hall algebras, let us digress briefly to make some preparation. Affine Lie algebras can be viewed as a special class of Kac-Moody algebras and also as central extensions of loop algebras. Amazingly, the affine quantum groups admit a current presentation (known as Drinfeld presentation) besides the Serre presentation given by Drinfeld and Jimbo; see \cite{Dr87, Be94, Da15}. 
	
	Ringel's construction can be formally extended to  hereditary abelian categories, which contain as distinguished examples the categories of quiver representations and the categories of coherent sheaves of (weighted) projective lines over $\bfk$. The Hall algebra of the projective line was shown in a pioneering work of Kaparanov \cite{Ka97} to be isomorphic to half quantum affine $\hA$ in its Drinfeld (current) presentation. Kaparanov's work has been extended in \cite{Sch04} (also \cite{BS13}) who developed the connection between Hall algebras of weighted projective lines and halves of quantum loop algebras. 
	
	The versions of Hall algebras following Steinitz, Ringel, and Bridgeland are reviewed in Section~\ref{sec:Hall}; we refer to \cite{DDPW} and \cite{Sch06} for in-depth expositions and extensive references on Hall algebras. Other important classes of Hall algebras, such as derived Hall algebras and cohomological Hall algebras, will not be considered in this survey paper.

	\subsection{Quantum groups}
	
	Drinfeld-Jimbo quantum groups $\U$ or Drinfeld double $\tU$ are by now well-known subjects \cite{Dr87, Lus90}. Some fundamental constructions on quantum groups which are more relevant to this survey include PBW bases, canonical bases, Hall algebra realizations, braid group actions, and Drinfeld presentation for affine quantum groups.
	
	Quantum groups admit a natural generalization known as $\imath$quantum groups.

	\subsection{$\imath$Quantum groups}
	
	Symmetric pairs can be classified in terms of Satake diagrams denoted by $(\I =\Iblack \cup \I_\circ, \tau)$, where $\I =\Iblack \cup \I_\circ$ is a suitable bicolored partition of the Dynkin diagram $\I$ together with a diagram automorphism $\tau$ such that $\tau^2 =\Id$; cf. \cite{OV90}. 
	Starting with Satake diagrams, Letzter \cite{Let99, Let02} introduced quantum symmetric pairs $(\U, \Ui)$ of finite type, where $\Ui =\Ui_\bvs$ (call an $\imath$quantum group) is a coideal subalgebra of a Drinfeld-Jimbo quantum group $\U$ depending on parameters $\bvs$; Kolb \cite{Ko14} improved and extended the constructions of quantum symmetric pairs to Kac-Moody type. In this survey, we shall restrict ourselves to the quasi-split $\imath$quantum groups $\Ui$, associated to $\I=\I_\circ$ and $\Iblack =\emptyset$. The authors introduced in \cite{LW22} the universal (quasi-split) $\imath$quantum group  $\tUi$, which is a coideal subalgebra of the Drinfeld double $\tU$. The $\imath$quantum groups $\Ui_\bvs$ for various parameters are obtained from $\tUi$ by central reductions. 
	
	Drinfeld-Jimbo quantum groups can be viewed as $\imath$quantum groups associated to diagonal Satake diagrams (recall a diagonal Satake diagram is a union of 2 identical copies of a Dynkin diagram with an involution swapping corresponding vertices). This viewpoint has been fruitful in the $\imath$program (cf. \cite{BW18a}) which seeks to generalize all fundamental (algebraic, geometric, and categorical) constructions in quantum groups to  $\imath$quantum groups; see \cite{W22} for a survey.
	
	In contrast to quantum groups (for which we have $\U =\U^-\U^0\U^+$), the $\imath$quantum groups in general do not admit a natural triangular decomposition. For quantum groups, there is a unique rank 1 case (i.e., the quantum $\mathfrak{sl}_2$) and there is a uniform Serre relation associated to all rank 2 Dynkin diagrams. For quasi-split $\imath$quantum groups, there are 3 (real) rank 1 cases, and accordingly there are several Serre type relations which have to be formulated separately. 
	
	$\imath$Divided powers arise as $\imath$canonical basis for the split rank 1 $\imath$quantum group (with a distinguished parameter) \cite{BW18a, BeW18}, and they have since been generalized to $\Ui$ for arbitrary parameters \cite{CLW21a} and also to the universal $\imath$quantum groups $\tUi$ \cite{LW20}. A Serre presentation for quasi-split $\imath$quantum groups $\tUi$ of arbitrary Kac-Moody type was obtained in \cite{CLW21a} (generalizing the earlier presentation of Letzter \cite{Let02} in finite type and some relation in \cite{BK15}), where a new key ingredient is the so-called $\imath$Serre relation in terms of $\imath$divided powers. 
	
	The basics on quantum groups and $\imath$quantum groups are reviewed in Section~\ref{sec:QG}. 

	\subsection{$\imath$Hall algebras}
	
	This survey is mainly concerned about the constructions of $\imath$Hall algebras and their use in realizing the quasi-split $\imath$quantum groups. The reflection functors on $\imath$Hall algebras are then used to realize the relative braid group action on $\imath$quantum groups. We also develop the connection between the $\imath$Hall algebra of the Jordan quiver and ring of symmetric functions.  

	\subsubsection{Semi-derived Hall algebras}
	
	With the help of the Riedtmann-Peng formula, the definition of Hall algebras has been extended to exact categories by Hubery \cite{Hub06} (also cf. \cite{PX00}). Inspired by Bridgeland's construction \cite{Br}, M.~Gorsky \cite{Gor2} constructed semi-derived Hall algebras for Frobenius categories; see also \cite{Gor1} for the semi-derived Hall algebras constructed for $\Z_2$-graded complexes of an exact category admitting enough projective objects. 
	
	More recently, motivated by the works of Bridgeland and Gorsky, the first author and Peng \cite{LP21} formulated the {\em semi-derived Ringel-Hall algebras} for hereditary abelian categories, which may not have enough projectives.
	In \cite{Lu22}, the construction in \cite{LP21} is further extended to weakly 1-Gorenstein exact categories. Note that the module categories of $\imath$quiver algebras below are 1-Gorenstein but not hereditary. 
	
	
	We give an overview of semi-derived Hall algebras in Section~\ref{sec:Semi-derivedHall}. 

	\subsubsection{$\imath$Quiver algebras}
	
	Starting with $\imath$quivers $(Q, \tau)$ (aka quivers $Q$ with involutions $\tau$ preserving the arrows of $Q$, where $\tau=\id$ is allowed), we constructed in \cite{LW22} a class of 1-Gorenstein algebras $\Lambda^\imath$; in case when $\tau =\Id$,  the $\imath$quiver algebras can be identified with $\bfk Q \otimes \bfk[\varepsilon]/(\varepsilon^2)$, and these algebras also appeared in \cite{RZ} and \cite{GLS} from very different perspectives. The $\imath$quiver algebras can be infinite-dimensional if the quiver $Q$ is not acyclic; but in all cases, they are 1-Gorenstein and admit favorable homological properties, cf. \cite{LW22, LW20}. 
	The $\imath$Hall algebra $\tMHk$ is by definition the twisted semi-derived Hall algebra of (the module category of) the $\imath$quiver algebra $\Lambda^\imath$.
	
	The $\imath$quiver algebra associated to the diagonal $\imath$quiver $(Q\cup Q, \text{swap})$ is isomorphic to $\bfk Q \otimes R_2$, where $R_2$ is the exterior $\bfk$-algebra in 2 variables $\varepsilon, \varepsilon'$;
	see \S\ref{subsec:i-QA}. We show that the $\imath$Hall algebra in this case is isomorphic to Drinfeld double $\tU$, and this provides a reformulation of the main theorem of Bridgeland \cite{Br}. 
	
	The $\imath$quiver algebras and their basic homological properties were presented in Section~\ref{sec:iQA}. 

	\subsubsection{$\imath$Hall algebras for Jordan quiver}
	
	A general theory of quiver representations is developed often for quivers without loops, and so representations of the Jordan quiver requires some separate discussions. The framework of $\imath$quivers does allow to formulate Jordan $\imath$quiver naturally; see Example~\ref{example 2}(5) and Figure~\ref{fig:Jordan}. 
	
	It is shown in \cite{LRW21} that the $\imath$Hall algebra of Jordan quiver is commutative and it is isomorphic to the ring of symmetric functions (now in 2 variables $t$ and $\theta$), extending the classic Steinitz-Hall isomorphism. The $\imath$Hall basis now corresponds to $\imath$Hall-Littlewood functions $Q_\la^\imath$ (or its modified version $H_\la^\imath$), which admit a formulation via raising and lowering operators. One may view the (modified) HL functions as the leading term in the variable $\theta$ of the (modified) $\imath$HL functions, as 
	$Q_\la^\imath$ and $H_\la^\imath$ at $\theta=0$ reduce to $Q_\la$ and $H_\la$. When setting $\theta=1$, the $\imath$HL functions are identified with type C deformed universal characters introduced in \cite{SZ06}. Several $\imath$Pieri rules for $\imath$HL functions were obtained through mixed combinatorial and $\imath$Hall approaches; the combinatorial approach follows and generalizes \cite{Tam11}. 
	
	We present the $\imath$Hall algebra of Jordan quiver and its connection to ring of symmetric functions in Section~\ref{sec:iJordan}. 

	\subsubsection{$\imath$Hall algebras for $\imath$quivers}
	
	The $\imath$Hall algebra $\tMHk$ for a Dynkin $\imath$quiver is shown by the authors \cite{LW22} to be isomorphic to the universal $\imath$quantum group $\tUi$. The relations (including the Serre relations) of $\tUi$ can be verified directly in the framework of $\tMHk$, as there are only a few local rank 2 cases to consider. 
	The connection between $\imath$Hall algebras for Dynkin $\imath$quivers and $\imath$quantum groups of finite type is developed in Section~\ref{sec:iHallDynkin}. 
	
	It requires some major work to extend the above construction for Dynkin $\imath$quivers to a class of virtually acyclic (including all acyclic) $\imath$quivers; in this generality, we obtain an injective homomorphism from the $\imath$quantum group $\tUi$ to the corresponding $\imath$Hall algebra \cite{LW20}. The verification of $\imath$Serre relations (among others) in an $\imath$Hall algebra is highly nontrivial. To that end, we first need to understand the  $\imath$divided powers in terms of an $\imath$Hall basis. Then we express the $\imath$Serre formula in terms of an $\imath$Hall basis, through nontrivial homological computations. Finally, we reduce the proof of the vanishing of various coefficients in the $\imath$Hall basis for the $\imath$Serre formula into some new $q$-binomial identities. 
	
	The connection between $\imath$Hall algebras for virtually acyclic $\imath$quivers and $\imath$quantum groups of Kac-Moody type is developed in Section~\ref{sec:iHallKM}. 

	\subsubsection{Reflection functors}
	
	Quiver representations have been generalized to representations of modulated graphs by Dlab-Ringel, and the latter can be adapted to $\imath$quivers. BGP type reflection functors can also be formulated in the setting of representations of modulated graphs, and they can be realized in terms of tilting modules; see \cite{Li12, GLS, LW21a, LW22b}.
	
	For a sink $\ell$ of $Q_0$, let $\bs_\ell Q$ be the quiver obtained from $Q$ by reversing the arrows ending at $\ell$ and $\btau \ell$. Denote by $\iLa$ (respectively, $\bs_\ell \iLa$)  the $\imath$quiver algebra for the $\imath$quiver $(Q,\btau)$ (respectively, $(\bs_\ell Q,\btau)$). A reflection functor from $\mod(\Lambda^{\imath})$ to $\mod(\bs_\ell \iLa )$ can then be defined \cite{LW21a, LW22b}. This induces an isomorphism of $\imath$Hall algebras, $\Gamma_{\ell}:\widetilde{\ch}(\bfk Q,\btau)  \stackrel{\sim}{\longrightarrow} \widetilde{\ch}(\bfk \bs_\ell Q,\btau)$, with explicit formulas on generators; and when composed with a Fourier transformation, this further induces an automorphism of $\tMHk$.
	
	Under the $\imath$Hall algebra realization of the $\imath$quantum group $\tUi$, we translate the (reflection functor) automorphisms of $\tMHk$ to (relative braid group) symmetries on $\tUi$, which are given by explicit formulas on generators; this confirms substantial cases of \cite[Conjecture 6.5]{CLW21a} and \cite[Conjecture~ 3.7]{CLW21b}. These conjectures have been most recently completely settled by Weinan Zhang \cite{Z22} who developed further the approach initiated in \cite{WZ22}. Existence of relative braid group symmetries on $\imath$quantum groups was originally conjectured in \cite{KP11}. 
	
	Reflection functors on $\imath$Hall algebras and relative braid group actions on $\tUi$ are presented in Section~\ref{sec:braid}. 

	\subsubsection{$\imath$Hall algebras for (weighted) projective lines}
	
	Drinfeld (current) presentation of affine quantum groups were obtained in \cite{Dr87, Be94, Da15}.
	Recently, the Drinfeld type presentations of split affine $\imath$quantum groups were obtained in \cite{LW21c} for ADE type and then in \cite{Z21} for BCFG type. The affine $\imath$quantum group of split type $A_1^{(1)}$ is also known as $q$-Onsager algebra. 
	
	For a hereditary abelian category $\ca$, let $\cc_{\Z_1}(\ca)$ be the category of $1$-periodic complexes of $\ca$. Denote its $\imath$Hall algebra by ${}^\imath\widetilde{\ch}(\ca)$, which is the twisted semi-derived Ringel-Hall algebra of  $\cc_{\Z_1}(\ca)$. 
	
	When $\ca$ is the category of coherent sheaves on the projective line $\bbP^1_\bfk$, We show in \cite{LRW20} that the $\imath$Hall algebra  ${}^\imath\widetilde{\ch}(\ca)$ is isomorphic to the (universal) $q$-Onsager algebra in its Drinfeld type presentation. The isomorphism between two presentations of $q$-Onsager algebra is conceptually explained by the derived equivalence between category of modules over Kronecker $\imath$quiver algebra and category $\cc_{\Z_1}(\ca)$.
	
	The above picture admits a vast generalization \cite{LR21} when $\ca$ is the category of coherent sheaves on a weighted projective line. In this case, one needs to replace $q$-Onsager algebra by $\imath$quantum loop algebras, which are by definition a generalization of Drinfeld presentation of split affine $\imath$quantum groups.
	
	We review the Drinfeld type presentation of affine $\imath$quantum groups and $\imath$Hall algebras of weighted projective lines in Section~\ref{sec:iWPL}.

	
	%
	%
	\subsection{Notations}
	
	We list the notations which are often used throughout the paper.
	\smallskip
	
	$\triangleright$ $\N,\Z,\Q$, $\C$ -- sets of nonnegative integers, integers, rational  and complex numbers,
	
	$\triangleright$ $\bfk$ -- a finite field of $q$ elements,
	
	$\triangleright$ $|\cs|$ -- cardinality of a finite set $\cs$.

	\medskip
	
	For a quiver algebra $A=\bfk Q/I$ (not necessarily finite-dimensional), we always identify left $A$-modules with representations of $Q$ satisfying relations in $I$. A representation $V=(V_i,V(\alpha))_{i\in Q_0,\alpha\in Q_1}$ of $A$ is called {\em nilpotent} if for each oriented cycle $\alpha_m\cdots\alpha_1$ at a vertex $i$, the $\bfk$-linear map $V(\alpha_m)\cdots V(\alpha_1):V_i\rightarrow V_i$ is nilpotent. We denote
	
	
	$\triangleright$ $\proj(A)$ -- category of finitely generated projective $A$-modules,
	
	$\triangleright$ $\mod(A)$ -- category of finite-dimensional nilpotent $A$-modules,
	
	
	
	$\triangleright$ $\cd^b(\mod(A))$ -- bounded derived category for $\mod(A)$,
	
	$\triangleright$ ${\rm proj.dim}_AM$ -- projective dimension of an $A$-module $M$,
	
	$\triangleright$ ${\rm inj.dim}_AM$ -- injective dimension of $M$,
	
	$\triangleright$ $\cp^{\leq m}(A)$ -- subcategory of $\mod(A)$ with projective dimension $\leq m$, for $m\geq0$,
	
	$\triangleright$ $\cp^{<\infty}(A)$ -- subcategory of $\mod(A)$ with finite projective dimension,
	
	$\triangleright$ $\Gproj(A)$ -- category of Gorenstein-projective $A$-modules.
	
	\medskip
	
	Let $Q=(Q_0,Q_1)$ be a quive, and $A=\bfk Q/I$ be a quiver algebra.
	For $i \in Q_0$, we denote
	
	$\triangleright$ $e_i$ -- the primitive idempotent of $A$,
	
	$\triangleright$ $S_i$ -- the $1$-dimensional simple nilpotent $A$-module supported at $i$,
	
	$\triangleright$ $P_i$ -- the projective cover of $S_i$,
	
	$\triangleright$ $I_i$ --  the injective hull of $S_i$.

	\medskip
	
	For an additive category $\ca$ and $M\in \ca$, we denote
	
	$\triangleright$ $\add M$ -- subcategory of $\ca$ whose objects are the direct summands of finite direct sums of copies of $M$,
	
	$\triangleright$ $\Ind (\ca)$ --  set of the isoclasses of indecomposable objects in $\ca$,
	
	$\triangleright$ $\Iso(\ca)$ -- set of the isoclasses of objects in $\ca$,
	
	$\triangleright$ $\aut(M)$ -- automorphism group of $M$.

	\medskip
	
	For an exact category $\ca$, we denote
	
	$\triangleright$ $K_0(\ca)$ --  Grothendieck group of $\ca$,
	
	$\triangleright$ $\widehat{A}$ -- the class in $K_0(\ca)$ of $A \in \ca$,
	
	$\triangleright$ $\cc_{\Z_m}(\ca)$ -- the category of $\Z_m$-graded complexes,
	
	$\triangleright$ $\cd_{sg}(\ca)$ -- the singularity category of $\ca$.
	
	\medskip
	
	For various Hall algebras, we denote
	
	$\triangleright$ $\ch (\ca)$ -- Hall algebra of the exact category $\ca$,
	
	$\triangleright$ $\ch (\bfk Q)$ -- Hall algebra of the path algebra $\bfk Q$,
	
	$\triangleright$ $\widetilde{\ch} (\bfk Q)$ -- twisted Hall algebra of the path algebra $\bfk Q$,
	
	$\triangleright$ $D\widetilde{\ch} (\bfk Q)$ -- Bridgeland's Hall algebra of the path algebra $\bfk Q$ (with $Q$ acyclic).

	\smallskip
	Associated to the $\imath$quiver $(Q,\btau)$ (aka quiver with involution), we denote
	
	$\triangleright$ $\La^\imath$ -- $\imath$quiver algebra, 
	
	$\triangleright$ $\utMH$ -- semi-derived Hall algebra for $\La^\imath$ (i.e., for $\mod (\La^\imath)$),
	
	$\triangleright$ $\ct (\La^\imath)$ -- quantum torus, 
	
	$\triangleright$ $\tMH=\iH(\bfk Q,\btau)$ -- $\imath$Hall algebra 
	(=twisited semi-derived Hall algebra for $\La^\imath$),
	
	$\triangleright$ $\tTL$ -- twisted quantum torus. 
	
	\medskip
	
	$\triangleright$ $\QJ$--Jordan quiver,
	
	$\triangleright$ $\PL$ -- the projective line over $\bfk$,
	
	$\triangleright$	$\coh(\PL)$ -- category of coherent sheaves of $\PL$,
	
	$\triangleright$	$\X$ -- weighted projective line over $\bfk$,
	
	$\triangleright$	$\coh(\X)$ -- catogory of coherent sheaves of $\X$.

	\medskip
	For quantum algebras, we denote
	
	$\triangleright$ $\U$ -- quantum group,
	
	$\triangleright$ $\tU$ -- Drinfeld double (a variant of $\U$ with doubled Cartan subalgebra),
	
	$\triangleright$ $\Ui=\Ui_{\bvs}$ -- a right coideal subalgebra of $\U$, depending on parameter $\bvs \in (\Q(v)^\times)^\I$,
	
	$\triangleright$ $(\U, \Ui)$ -- quantum symmetric pair,
	
	$\triangleright$ $(\tU, \tUi)$ -- a universal quantum symmetric pair, such that $\Ui$ is obtained from $\tUi$ by a central reduction,
	
	$\triangleright$ $\tUiD$ -- $\imath$quantum loop algebra.
	
	%
	\vspace{2mm}
	\noindent {\bf Acknowledgement.}
	ML thanks Liangang Peng for guidance and continuing encouragement, and thanks University of Virginia, Shanghai Key Laboratory of Pure Mathematics and Mathematical Practice, East China Normal University for hospitality and support.
	ML is partially supported by the Science and Technology Commission of Shanghai Municipality (grant No. 18dz2271000) and the National Natural Science Foundation of China (grant No. 12171333). WW is partially supported by the NSF grant DMS-2001351.

	\section{Quantum groups and $\imath$quantum groups}
	\label{sec:QG}
	
	In this section, we recall Drinfeld double quantum groups, and review the (quasi-split) universal $\imath${}quantum groups from \cite{LW22} whose central reductions are the $\imath$quantum groups \`a la Letzter and Kolb \cite{Ko14}. We then give a Serre presentation of the (universal) $\imath$quantum groups via the $\imath$divided powers, following \cite{CLW21a}.

	\subsection{Quantum groups and Drinfeld doubles}
	\label{subsection Quantum groups}
	
	Let $Q$ be a quiver (without loops) with vertex set $Q_0= \I$. The underlying graph of $Q$ with orientations of edges forgotten is called a Dynkin diagram, denoted by $D_Q$. 
	Let $n_{ij}$ be the number of edges connecting vertex $i$ and $j$. Let $C=(c_{ij})_{i,j \in \I}$ be the symmetric generalized Cartan matrix of $D_Q$, defined by $c_{ij}=2\delta_{ij}-n_{ij}.$ Let $\fg$ be the Kac-Moody Lie algebra associated to $C$, with simple roots denoted by $\alpha_i$ ($i\in\I $); for this survey of Hall algebras, we can simply take $\fg$ to be the derived subalgebra of the Kac-Moody algebra. The root lattice is denoted by $\Z^\I:=\bigoplus_{i\in\I}\Z\alpha_i$. The {\em simple reflection} $s_i:\Z^{\I}\rightarrow\Z^{\I}$ is defined to be $s_i(\alpha_j)=\alpha_j-c_{ij}\alpha_i$, for $i,j\in \I$.
	Denote the Weyl group by $W =\langle s_i\mid i\in \I\rangle$.
	
	Let $\bv$ be an indeterminant. Write $[A, B]=AB-BA$. Denote, for $r \in \N$ and $m\in \Z$,
	\[
	[r]=\frac{\bv^r-\bv^{-r}}{\bv-\bv^{-1}},
	\quad
	[r]!=\prod_{i=1}^r [i], \quad \qbinom{m}{r} =\frac{[m][m-1]\ldots [m-r+1]}{[r]!}.
	\]
	Then $\tU = \tU_\bv(\fg)$ is the $\Q(\bv)$-algebra generated by $E_i,F_i, \tK_i,\tK_i'$, $i\in \I$, where $\tK_i, \tK_i'$ are invertible, subject to the following relations:
	\begin{align}
		[E_i,F_j]= \delta_{ij} \frac{\tK_i-\tK_i'}{\bv-\bv^{-1}},  &\qquad [\tK_i,\tK_j]=[\tK_i,\tK_j']  =[\tK_i',\tK_j']=0,
		\label{eq:KK}
		\\
		\tK_i E_j=\bv^{c_{ij}} E_j \tK_i, & \qquad \tK_i F_j=\bv^{-c_{ij}} F_j \tK_i,
		\label{eq:EK}
		\\
		\tK_i' E_j=\bv^{-c_{ij}} E_j \tK_i', & \qquad \tK_i' F_j=\bv^{c_{ij}} F_j \tK_i',
		\label{eq:K2}
	\end{align}
	and the quantum Serre relations, for $i\neq j \in \I$,
	\begin{align}
		& \sum_{r=0}^{1-c_{ij}} (-1)^r  E_i^{(r)} E_j  E_i^{(1-c_{ij}-r)}=0,
		\label{eq:serre1} \\
		& \sum_{r=0}^{1-c_{ij}} (-1)^r   F_i^{(r)} F_j  F_i^{(1-c_{ij}-r)}=0.
		\label{eq:serre2}
	\end{align}
	Here we have used the divided powers
	\[
	F_i^{(n)} =F_i^n/[n]!, \quad E_i^{(n)} =E_i^n/[n]!, \quad \text{ for } n\ge 1, i\in \I.
	\]
	Note that $\tK_i \tK_i'$ are central in $\tU$ for all $i$.
	We endow $\tU$ a Hopf algebra structure with the comultiplication $\Delta: \widetilde{\U} \longrightarrow \widetilde{\U} \otimes \widetilde{\U}$ given by
	\begin{align}  \label{eq:Delta}
		\begin{split}
			\Delta(E_i)  = E_i \otimes 1 + \tK_i \otimes E_i, & \quad \Delta(F_i) = 1 \otimes F_i + F_i \otimes \tK_{i}', \\
			\Delta(\tK_{i}) = \tK_{i} \otimes \tK_{i}, & \quad \Delta(\tK_{i}') = \tK_{i}' \otimes \tK_{i}'.
		\end{split}
	\end{align}
	The Chevalley involution $\omega$ on $\tU$ is given by
	\begin{align}  \label{eq:omega}
		\omega(E_i)  = F_i,\quad  \omega(F_i) = E_i,\quad \omega(\tK_{i}) = \tK_{i}' , \quad \omega(\tK_{i}') =\tK_{i}, \quad \forall i\in \I.
	\end{align}
	
	Analogously as for $\tU$, the quantum group $\bU$ is defined to be the $\Q(v)$-algebra generated by $E_i,F_i, K_i, K_i^{-1}$, $i\in \I$, subject to the  relations modified from \eqref{eq:KK}--\eqref{eq:serre2} with $\tK_i$ and $\tK_i'$ replaced by $K_i$ and $K_i^{-1}$, respectively. The comultiplication $\Delta$ and Chevalley involution $\omega$ on $\U$ are obtained by modifying \eqref{eq:Delta}--\eqref{eq:omega} with $\tK_i$ and $\tK_i'$ replaced by $K_i$ and $K_i^{-1}$, respectively (cf. \cite{Lus93}; beware that our $K_i$ has a different meaning from $K_i \in \U$ therein.)
	
	The algebra $\U$ can be obtained from $\tU$ by a central reduction; that is, $\U$ is isomorphic to a quotient algebra of $\tU$ by the ideal $( \tK_i \tK_i'- 1 \mid \forall i\in \I )$.
	
	Let $\widetilde{\bU}^+$ be the subalgebra of $\widetilde{\bU}$ generated by $E_i$ $(i\in \I)$, $\widetilde{\bU}^0$ be the subalgebra of $\widetilde{\bU}$ generated by $\tK_i, \tK_i'$ $(i\in \I)$, and $\widetilde{\bU}^-$ be the subalgebra of $\widetilde{\bU}$ generated by $F_i$ $(i\in \I)$, respectively.
	The subalgebras $\bU^+$, $\bU^0$ and $\bU^-$ of $\bU$ are defined similarly. Then both $\widetilde{\bU}$ and $\bU$ have triangular decompositions:
	$
	\widetilde{\bU} =\widetilde{\bU}^+\otimes \widetilde{\bU}^0\otimes\widetilde{\bU}^-,
	\, 
	\bU =\bU^+\otimes \bU^0\otimes\bU^-.
	$ 
	Clearly, ${\bU}^+\cong\widetilde{\bU}^+$, ${\bU}^-\cong \widetilde{\bU}^-$, and ${\bU}^0 \cong \widetilde{\bU}^0/(\tK_i \tK_i' -1 \mid   i\in \I)$.
	
	\begin{remark}
		\begin{enumerate}
			\item 
			$\tU$ (respectively, $\U$) is a $\Z^\I$-graded algebra with $\wt E_i=\alpha_i$, $\wt F_i=-\alpha_i$, and $\wt \tK_i=0=\wt \tK_i'$ (respectively, $\wt K_i=0 =\wt K_i^{-1}$). 
			\item
			There is a second $\Z^\I$-grading on $\tU$ (which do not descend to $\U$) by setting $\deg E_i =\alpha_i=\det F_i$, and $\deg \tK_i=\alpha_i+\alpha_{\tau i}=\deg \tK_i'$.
		\end{enumerate}
	\end{remark}

	\begin{example}[Rank 1]
		There is a unique quantum group of rank $1$ corresponding to $\I$ which consists of a single vertex. In this case, $\U=\U_v(\mathfrak{sl}_2)$ is generated by $E,F,K^{\pm1}$ such that $KK^{-1}=K^{-1}K=1$, $KE=v^2EK$, $KF=v^{-2}FK$, and $[E,F]=\frac{K-K^{-1}}{v-v^{-1}}$.
	\end{example}
	
	\begin{example}[Rank 2]
		The quantum groups $\U$ of rank $2$  are associated to the Dynkin diagram $\I$ consisting of two vertices $1$ and $2$ with $a$ edges between them; in this case  $c_{12} =-a$, and we draw the diagram as 
		\[\xymatrix{1\ar@{-}[r]|-{a}& 2}\]
		In particular, if $a=1$, then $\U=\U_v(\mathfrak{sl}_3)$ is of type $A_2$; and if $a=2$, then $\U=\U_v(\hA)$ is of affine type $A_1^{(1)}$.
	\end{example}

	\subsection{(Universal) $\imath$quantum groups $\Ui$ and $\tUi$}
	\label{subsec:iQG}
	
	For a generalized Cartan matrix $C=(c_{ij})$, let $\Aut(C)$ be the group of all permutations $\btau$ of the set $\I$ such that $c_{ij}=c_{\btau i,\btau j}$. An element in $\Aut(C)$ is the same as a diagram automorphism of the Dynkin diagram associated to $C$. An element $\btau\in\Aut(C)$ is called an \emph{involution} if $\btau^2=\Id$; we allow $\tau =\Id$. Such a $\tau$ induces an automorphism of $\fg$, denoted again by $\tau$. 
	
	Let $\btau$ be an involution in $\Aut(C)$. We define $\widetilde{\bU}^\imath$ 
	to be the $\Q(v)$-subalgebra of $\tU$ generated by
	\[
	B_i= F_i +  E_{\btau i} \tK_i',
	\qquad \tk_i = \tK_i \tK_{\btau i}', \quad \forall i \in \I.
	\]
	Indeed, $\tUi$ is a right coideal subalgebra of $\tU$, i.e., $\Delta (\tUi) \subset \tUi\otimes \tU$; this can be checked on generators of $\tUi$. Let $\tU^{\imath 0}$ be the commutative $\Q(v)$-subalgebra (called the Cartan subalgebra) of $\tUi$ generated by $\tk_i$, for $i\in \I$.
	The elements
	\begin{align}
		\label{eq:central}
		\tk_i\; (i= \btau i)
		\qquad
		\tk_i \tk_{\btau i}\;  (i\neq \btau i)
	\end{align}
	are central in $\tUi$.
	
	
	
	Let 
	\[
	\bvs=(\vs_i \mid i\in \I) \in  (\Q(\bv)^\times)^{\I}
	\]
	be such that $\vs_i=\vs_{\btau i}$ whenever $i\in \I$ satisfies $c_{i, \btau i}=0$.
	Let $\Ui_{\bvs}$ be the $\Q(v)$-subalgebra of $\bU$ generated by
	\[
	B_i= F_i+\vs_i E_{\btau i}K_i^{-1},
	\quad
	k_j= K_jK_{\btau j}^{-1},
	\qquad  \forall i \in \I, j \in \ci.
	\]
	We often write $\Ui$ for $\Ui_{\bvs}$ by dropping the parameter $\bvs$ when there is no confusion. 
	It is known \cite{Let99, Ko14} that $\bU^\imath$ is a right coideal subalgebra of $\bU$, i.e., $\Delta (\Ui) \subset \Ui\otimes \U$; and $(\bU,\Ui)$  
	is called a (quasi-split) \emph{quantum symmetric pair}, and it specializes at $v=1$ to (the enveloping algebra version of) the symmetric pair $(\fg, \fg^\theta)$, where $\theta=\omega \circ \btau$. (The Chevalley involution $\omega$ and the involution $\tau$ on $\fg$ commute.) 
	
	We call $\Ui$ a (quasi-split) $\imath$quantum group and $\tUi$ a (quasi-split) universal $\imath$quantum group. A quasi-split $\imath$quantum group is called {\em split} if in addition we have $\tau=\Id$.
	
	\begin{remark}
		More general $\imath$quantum groups can be defined associated to Satake diagrams $(\I =\I_{\bullet} \cup \I_\circ, \tau)$ \cite{Let99, Ko14}, and so are the universal $\imath$quantum groups (cf. \cite{WZ22}). For the purpose of $\imath$Hall algebra realization in this paper, we only consider quasi-split universal $\imath$quantum groups as defined above, which are associated to Satake diagrams with $\I_\bullet =\emptyset$.  
	\end{remark}
	
	The algebras $\Ui_{\bvs}$, for various parameters $\bvs$, are obtained from $\tUi$ by central reductions.
	
	\begin{proposition}[\text{\cite[Propositon 6.2]{LW22}}]
		\label{prop:QSP12}
		The $\Q(v)$-algebra $\Ui$ is isomorphic to the quotient of $\tUi$ by the ideal $I^\imath$ generated by
		$
		\tk_i - \vs_i \; (\text{for } i =\btau i)$ and
		$\tk_i \tk_{\btau i} - \vs_i \vs_{\btau i}  \;(\text{for } i \neq \btau i).
		$ 
		The isomorphism $\Ui \rightarrow \tUi/ I^\imath$ is given by sending 
		\[
		B_i \mapsto B_i, \quad 
		k_j \mapsto \vs_{\btau j}^{-1} \tk_j, \quad
		k_j^{-1} \mapsto \vs_{ j}^{-1} \tk_{\btau j}, \quad 
		\forall i\in \I, \, j\in \I\setminus\ci.
		\]
	\end{proposition}
	
	The {\em rank} of $\tUi$ is by definition the number of $\btau$-orbits in $\I$.

	\begin{example}[Split type of rank $1, 2$]
		\begin{enumerate}
			\item 
			The split (universal) $\imath$quantum groups $\tUi$ and $\Ui$ of rank $1$ are associated to $\I$ which consists of a single vertex $1$. In this case, $\tUi$ is a commutative algebra generated by $B, \tk^{\pm1}$, while $\Ui$ is isomorphic to the polynomial algebra $\Q(v)[B]$. 
			
			\item
			The split $\imath$quantum groups $\tUi$ and $\Ui$ of rank $2$ are associated to the following (Satake) diagram
			
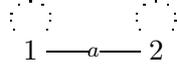
\vspace{1cm}
			\[\xymatrix@C=0.5cm{1\ar@(ru,lu)@{.}[]\ar@{-}[rr]|-{a}&&2\ar@(ru,lu)@{.}[]}\]
			For $a=2$, the $\imath$quantum groups $\Ui$ and $\tUi$ are known as $q$-Onsager algebra.
		\end{enumerate}
	\end{example}
	
	\begin{example}[Quasi-split type of rank $1, 2$]
		\begin{enumerate}
			\item 
			The quasi-split $\imath$quantum groups $\tUi$ and $\Ui$ of  rank $1$ are associated to the following (Satake) diagram 
			\[\xymatrix{1 \ar@/^0.6pc/@{.}[rr] \ar@<-0.5ex>@{-}[rr]|-{a}& &2  }\]
			\item
			The quasi-split $\imath$quantum groups $\tUi$ and $\Ui$ of rank $2$ are associated to the following (Satake) diagrams: 
			
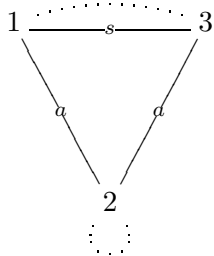
\begin{equation*}
				\xymatrix{&1\ar@<-0.5ex>@{-}[rr]|-{s}\ar@/^0.6pc/@{.}[rr] \ar@{-}[ddr]|-{a}  &&3 \ar@{-}[ddl]|-{a} 
					\\
					&&&
					\\
					&&2\ar@(rd,ld)@{.}[]&}  \qquad\qquad\qquad 
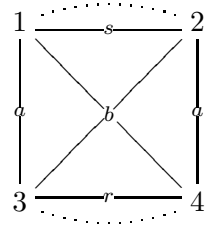
\xymatrix{&1 \ar@/^0.6pc/@{.}[rr]\ar@{-}[dd]|-{a} \ar@<-0.5ex>@{-}[rr]|-{s} \ar@{-}[dr] && 2 \ar@{-}[dd]|-{a} \ar@{-}[ddll]|-{b} \\
					&&\ar@{-}[dr]&\\
					&3 \ar@/_0.6pc/@{.}[rr] \ar@<0.5ex>@{-}[rr]|-{r}& &4 }
				\vspace{1cm}
			\end{equation*}
			In the second diagram above, 
			the number of edges between 1 and 4 (and respectively, between 2 and 3) is $b$. 
		\end{enumerate}
	\end{example}

	\subsection{A Serre presentation of $\tUi$}\label{subsection:iserre presentation}
	
	For  $\ell\in \I$ with $\btau \ell\neq \ell$,  we define the  $\imath${}divided power of $B_\ell$ as
	\begin{align}
		\label{eq:iDP1}
		B_\ell^{(m)}:=B_\ell^{m}/[m]!, \quad \forall m\ge 0, \qquad (\text{if } \ell \neq \btau \ell).
	\end{align}
	
	Let $i\in \I$ with $\btau i= i$. The {\em $\imath${}divided power} of $B_i$ in $\Ui_\bvs$ for a distinguished parameter $\vs_i =v^{-1}$ was originally introduced in \cite{BW18a, BeW18}, and then generalized for a general parameter $\vs_i$ in \cite{CLW21a}. Replacing $\vs_i$ by $\tk_i$ (see Proposition~\ref{prop:QSP12}) and following \cite[(2.20)-(2.21)]{CLW21b}, we define the {\em $\imath${}divided powers} of $B_i$ in $\tU$ to be 
	\begin{eqnarray}
		&&\ff_{i,\odd}^{(m)}=\frac{1}{[m]!}\left\{ \begin{array}{ccccc} B_i\prod_{s=1}^k (B_i^2-[2s-1]^2 v \tk_i) & \text{if }m=2k+1,\\
			\prod_{s=1}^k (B_i^2-[2s-1]^2 v \tk_i) &\text{if }m=2k; \end{array}\right.
		\label{eq:iDPodd}\\
		&&\ff_{i,\ev}^{(m)}= \frac{1}{[m]!}\left\{ \begin{array}{ccccc} B_i\prod_{s=1}^k (B_i^2-[2s]^2 v \tk_i) & \text{if }m=2k+1,\\
			\prod_{s=1}^{k} (B_i^2-[2s-2]^2 v \tk_i) &\text{if }m=2k. \end{array}\right.
		\label{eq:iDPev}
	\end{eqnarray}

	Denote
	\[
	(a;x)_0=1, \qquad (a;x)_n =(1-a)(1-ax)  \cdots (1-ax^{n-1}), \quad   n\ge 1.
	\]
	
	The following theorem is an upgrade of \cite[Theorem~3.1]{CLW21a} for $\Ui$ to the setting of a universal $\imath$quantum group $\tUi$. In the setting of $\Ui$, the relations \eqref{relation3} and \eqref{relation5} below were obtained earlier in \cite{Let99}-\cite{Ko14} and \cite{BK15}, respectively, while the $\imath$Serre relation \eqref{relation6} involving $\imath$divided powers \eqref{eq:iDPodd}--\eqref{eq:iDPev} was established in \cite{CLW21a}.
	
	\begin{theorem} [\text{Serre presentation for $\tUi$, \cite[Theorem 4.2]{LW20}}]
		\label{thm:Serre}
		Fix $\ov{p}_i\in \Z_2$ for each $i\in \I$. The $\Q(v)$-algebra $\tUi$ has a presentation with generators $B_i$, $\tk_i$ $(i\in \I)$ and the relations \eqref{relation1}--\eqref{relation5} below: for $\ell \in \I$, and $i\neq j \in \I$,
		\begin{align}
			\tk_i \tk_\ell =\tk_\ell \tk_i,
			\quad
			\tk_i B_\ell & = v^{c_{\btau i,\ell} -c_{i \ell}} B_\ell \tk_i,
			\label{relation1}
			\\
			B_iB_{j}-B_jB_i &=0, \quad \text{ if }c_{ij} =0 \text{ and }\btau i\neq j,\label{relation2}
			\\
			\sum_{n=0}^{1-c_{ij}} (-1)^n  B_{i, \overline{p_i}}^{(n)}B_j B_{i,\overline{c_{ij}}+\overline{p}_i}^{(1-c_{ij}-n)} &=0,\quad   \text{ if }\btau i=i \neq j,
			\label{relation6}
			\\
			\sum_{n=0}^{1-c_{ij}} (-1)^nB_i^{(n)}B_jB_i^{(1-c_{ij}-n)} &=0, \quad \text{ if } i\neq j \neq \btau i\neq i, \label{relation3}
			\\
			\sum_{n=0}^{1-c_{i,\btau i}} (-1)^{n+c_{i,\btau i}}B_i^{(n)}B_{\btau i}B_i^{(1-c_{i,\btau i}-n)}& =\frac{1}{v-v^{-1}} \times
			\label{relation5}     \\
			\left(v^{c_{i,\btau i}} (v^{-2};v^{-2})_{-c_{i,\btau i}} \right.  &
			\left. B_i^{(-c_{i,\btau i})} \tk_i
			-(v^{2};v^{2})_{-c_{i,\btau i}}B_i^{(-c_{i,\tau i})} \tk_{\btau i}  \right),
			\text{ if } \btau i \neq i.
			\notag
		\end{align}
		(Note for $\tUi$ of split type, only the relations \eqref{relation1}--\eqref{relation6} are needed.)
	\end{theorem}

	By the Serre presentation of $\tUi$ in Theorem~\ref{thm:Serre}, letting $\deg B_i=\alpha_i$ and $\deg \tk_i=0$, for $i\in\I$, endows $\tUi$ a $\N\I$-filtered algebra structure. Let $\tU^{\imath,\gr}$ be the associated graded algebra (with generators $\bar B_i$ and so on).
	
	\begin{remark}
		There exists a natural algebra monomorphism
		$\phi: \U^-\longrightarrow \tU^{\imath,\gr}$ which maps
		$F_i\mapsto \bar B_i$, for $i\in\I$. Moreover, the multiplication map  $\Im\phi \times \widetilde{\bU}^{\imath 0} \rightarrow \tU^{\imath,\gr}$ is an isomorphism of vector spaces. These are $\tUi$-variants of some basic  results on the structure of $\Ui$, which were established originally by Letzter \cite{Let02} and Kolb \cite{Ko14}.
	\end{remark}
	
	\section{Hall algebras}
	\label{sec:Hall}
	
	In this section, we introduce (Steinitz-)Hall algebra of the Jordan quiver and its connection to symmetric functions. We formulate Ringel-Hall algebras of quivers to realize halves of quantum groups, and then formulate Bridgeland-Hall algebras to realize the Drinfeld double quantum groups. 

	\subsection{Steinitz-Hall and Hall-Littlewood functions}
	\label{subsec:SHconstruction}
	
	Let $\co$ be a discrete valued ring with the maximal principal ideal $\mfp$ and a finite residue field $\bfk$ of $q$ elements.
	
	Let $\ca$ be the category of all finite $\co$-modules (i.e., modules with finitely many elements).
	Then $\ca$ has only one simple module, $\bfk$ (which is identified with $\co/\mfp$). The indecomposable objects of $\ca$ are $\co/\mfp^m$ ($m\geq1$) up to isomorphisms. Given a partition $\lambda=(\lambda_1,\dots,\lambda_n)$, we denote
	\begin{align*}
		M(\lambda):=\co/\mfp^{\lambda_1}\oplus \cdots\oplus \co/\mfp^{\lambda_n}.    
	\end{align*}
	Then every object $M\in\ca$ is isomorphic to $M(\lambda)$ for some unique partition $\lambda$.
	
	For any partitions $\lambda,\mu,\nu$, let $F_{\mu,\nu}^\lambda$ denote the number of sumodules $N$ of $M(\lambda)$ such that $N\cong M(\nu)$, and $M(\lambda)/N\cong M(\mu)$. 
	Define $\ch=\bigoplus_\lambda \Q \fu_\lambda$, a $\Q$-linear space with a basis indexed by partitions. We endow $\ch$ with an algebra structure given by
	\begin{align*}
		\fu_\mu*\fu_\nu:=\sum_{\lambda} F_{\mu,\nu}^\lambda \fu_\lambda, \quad \forall \mu, \nu.
	\end{align*}
	We call $\ch$ a (Steinitz-) Hall algebra, following \cite{Mac95}; it will be understood as Hall algebra of the Jordan quiver, see \S\ref{sec:iJordan}.
	
	Let $\Lambda$ be the $\Q$-algebra of symmetric functions in infinitely many variables. It is well known that  $\Lambda$ is a polynomial ring $\Q[e_r \mid r\ge 1]$, where $e_r$ denotes the $r$th elementary symmetric function. Let $t$ be an indeterminant. Denote by $\Lambda[t]:=\Lambda\otimes_\Q \Q[t]$.

	We recall the Hall-Littlewood (HL) functions in variables $x=(x_1, x_2, \ldots)$, following \cite[Pages 208--211]{Mac95}. Denote by $h_r$ the $r$th complete symmetric function, for $r\ge 0$.
	For an indeterminate $u$, denote
	\begin{align}
		H(u)&= \sum_{r=0}^\infty h_r u^r = \prod_{i\ge 1} \frac1{1-ux_i},
		\\
		Q(u)&= \sum_{r=0}^\infty Q_{r} u^r=H(u)/H(tu).
	\end{align}
	
	We set $q_r =Q_{r}$, and $q_\la =q_{\la_1} q_{\la_2} \cdots$, for any composition $\la=(\la_1,\la_2,\dots)$.
	Let $u_1,u_2,\cdots$ be independent indeterminates.  Let
	\begin{align}
		\label{eq:F}
		F(u)=\frac{1-u}{1-tu}=1+ (1-t^{-1}) \sum_{r\geq1} t^r u^r.
	\end{align}
	Then the HL symmetric function $Q_\la$ is the coefficient of $u^\la=u_1^{\la_1} u_2^{\la}\cdots$ in
	\begin{align}
		\label{eq:HL}
		Q(u_1,u_2,\dots)= \prod_{i\geq1}Q(u_i) \prod_{i<j} F(u_iu_j^{-1}).
	\end{align}
	
	Let $1\le i<j$. The raising operator $R_{ij}$ acts on integer sequences by
	\[
	R_{ij} \alpha =(\ldots, \alpha_i+1, \ldots, \alpha_j-1, \ldots).
	\]
	The actions of raising operators on $q_\alpha$ are given by letting
	\[
	R_{ij} q_\alpha =q_{R_{ij}\alpha}.
	\]
	
	According to \cite[p.212]{Mac95}, the generating function \eqref{eq:HL} for HL functions $Q_\la$ can be restated that
	\begin{align}  \label{eq:Rq}
		Q_\la = \prod_{i<j} \frac{1 -R_{ij}}{1 -tR_{ij}} q_\la
		=\Big(1 + (t-1) \sum_{r\ge 1} t^{r-1}  R_{ij}^r  \Big) q_\la.
	\end{align}
	Then $\{Q_\la\mid \la \mbox{ is a partition}\}$ forms a $\Q[t]$-basis of $\Lambda[t]$.

	For a composition $\alpha$, we denote
	\begin{align}
		\label{eq:number}
		\begin{split}
			|\alpha| &=\alpha_1 +\alpha_2+\cdots,
			\qquad
			n(\lambda)=\sum_i(i-1)\lambda_i,
		\end{split}
	\end{align}	
	and denote by $m_i(\lambda)$ the number of times $i$ occurs as a part of $\lambda$.

	For any parition $\la$, define $P_\la\in\Lambda[t]$ such that $Q_\la=b_\lambda(t) P_\la$, where 
	\begin{align}
		\label{eq:varphi:r}
		\varphi_r(t) &= (1-t)(1-t^2)\cdots (1-t^r);
		\\
		\label{def:bla}
		b_\la(t) &= \prod_{i\geq1} \varphi_{m_i(\la)}(t).
	\end{align}
	In fact, $P_\la\in\Lambda[t]$. 
	
	We view all functions in $\Lambda[t]$ to be in $\Lambda$ by setting $t=q^{-1}$. 
	
	\begin{theorem} [\text{Steinitz, Hall, \cite[II, III]{Mac95}}]
		\quad
		\begin{enumerate}
			\item 
			$\ch$ is a commutative associative algebra with the identity $\fu_\emptyset$.
			\item
			There is an isomorphism of algebras $\phi_{(q)}:\ch\rightarrow \Lambda$ which maps 
			$\fu_{(1^r)}\mapsto q^{-{r\choose2}}e_{r}$, for $r\ge 1$.
			\item 
			For any partition $\la$, we have
			$\phi_{(q)}(\fu_\la)=q^{-n(\la)}P_\la$.
		\end{enumerate}
	\end{theorem}
	
	In fact, $|\aut(M(\lambda))|=q^{|\lambda|+2n(\la)}b_\lambda(q^{-1})$; see e.g. \cite[Lemma 2.8]{Sch06}. Then we have
	\begin{align}
		\phi_{(q)}(|\aut(M(\lambda))|\cdot\fu_\la)= q^{|\lambda|+n(\lambda)} Q_\lambda,
	\end{align}
	for any partition $\lambda$.

	\subsection{Ringel-Hall and halves of quantum groups}
	
	Let $\ce$ be an essentially small exact category in the sense of Quillen, linear over  $\bfk=\F_q$. For the basics on exact categories, we refer to \cite{Buh} and references therein. 
	In particular, abelian categories are exact categories.	
	
	Assume that $\ce$ has finite morphism and extension spaces, i.e.,
	\[
	|\Hom(M,N)|<\infty,\quad |\Ext^1(M,N)|<\infty,\,\,\forall M,N\in\ce.
	\]
	
	Given objects $M,N,L\in\ce$, define $\Ext^1(M,N)_L\subseteq \Ext^1(M,N)$ as the subset parameterizing extensions whose middle term is isomorphic to $L$. We define the {\em Ringel-Hall algebra} $\ch(\ce)$ (or {\em Hall algebra} for short) to be the $\Q$-vector space whose basis is formed by the isoclasses $[M]$ of objects $M$ in $\ce$, with the multiplication defined by (see \cite{Br})
	\begin{align}
		\label{eq:mult}
		[M]\diamond [N]=\sum_{[L]\in \Iso(\ce)}\frac{|\Ext^1(M,N)_L|}{|\Hom(M,N)|}[L].
	\end{align}
	
	For any three objects $L,M,N$, let
	\begin{align}
		\label{eq:Fxyz}
		F_{MN}^L:= \big |\{X\subseteq L \mid X \cong N,  L/X\cong M\} \big |.
	\end{align}
	The Riedtmann-Peng formula states that
	\[
	F_{MN}^L= \frac{|\Ext^1(M,N)_L|}{|\Hom(M,N)|} \cdot \frac{|\Aut(L)|}{|\Aut(M)| |\Aut(N)|},
	\]
	where $\Aut(M)$ denotes the automorphism group of $M$. For any object $M$, 
	let
	\begin{align}
		\label{eq:doublebrackets}
		[\![M]\!]:=\frac{[M]}{|\Aut(M)|}.
	\end{align}
	Then the Hall multiplication \eqref{eq:mult} can be reformulated to be
	\begin{align}
		[\![M]\!]\diamond [\![N]\!]=\sum_{[\![L]\!]}F_{M,N}^L[\![L]\!],
	\end{align}
	which is the version of Hall multiplication used in \cite{Rin90}.
	
	Let $Q=(\I,Q_1)$ be a quiver (without loops), and $\bfk=\F_q$. 
	Let $\ca=\mod(\bfk Q)$ be the category of finite-dimensional nilpotent $\bfk Q$-modules. 
	Then $\ca$ is a hereditary abelian category. Let $\langle-,-\rangle$ be the Euler form of $\ca$, i.e., 
	\begin{align*}
		\langle M,N\rangle=\dim_\bfk\Hom(M,N)-\dim_\bfk\Ext^1(M,N), \,\,\forall M,N\in\ca.
	\end{align*}
	Throughout the paper we shall denote
	\begin{align}
		\label{v2q}
		\sqq=\sqrt{q}.
	\end{align}
	The twisted Hall algebra $\widetilde{\ch}(\ca)$ is defined over $\ch(\ca)\otimes_\Q \Q(\sqq)$ with multiplication twisted by Euler form:
	\begin{align*}
		[L]*[M]=\sqq^{\langle L,M \rangle} [L]\diamond[M].
	\end{align*}
	Denote by $S_i$ the $1$-dimensional simple nilpotent module supported at $i \in \I$. 
	
	Recall the $\Q(v)$-algebra $\U^-$, and consider the base change $\U^-_{\sqq} :=\Q(\sqq) \otimes_{\Q(v)} \U^-$.
	
	\begin{theorem}[\text{Ringel \cite{Rin90}, Green \cite{Gr95}}]
		\label{thm:Ringel-Green}
		Let $Q=(\I,Q_1)$ be a quiver, and   $\ca=\mod(\bfk Q)$. Then there exists an injective homomorphism of $\Q(\sqq)$-algebras $\phi: \U^-_{\sqq}\longrightarrow \widetilde{\ch}(\ca)$, which maps $F_i$ to $[\![S_i]\!]$, for all $i\in\I$. Moreover, $\phi$ is an isomorphism if and only if $Q$ is of Dynkin type. 
	\end{theorem}
	
	\begin{proof}
		The verification that $\phi(F_i)$ satisfy the quantum Serre relations is given below in Example~\ref{rank2ex} (also cf. Examples~\ref{rank1ex}--\ref{rankA2ex}). Hence we have a morphism $\phi: \U^-_{\sqq}\rightarrow \widetilde{\ch}(\ca)$ which maps $F_i$ to $[\![S_i]\!]$, for all $i\in\I$.
		
		We now sketch the proof of injectivity of $\phi$.
		
		\underline{Step 1}. We identify $\U^-$ with $\bf f$ in \cite[\S1.2]{Lus93}. Then
		$\U^-=\oplus_{\mu\in\N\I} \U^-_\mu$ is a graded bialgebra with $\deg F_i=\alpha_i$ and a comultiplication $r: \U^- \rightarrow \U^- \otimes \U^-$ such that $r(F_i) =F_i\otimes 1 + 1\otimes F_i$, for each $i\in \I$. By \cite[\S1.2]{Lus93}, there exists a non-degenerate symmetric bilinear form $(\cdot, \cdot)$ on $\U^-_{\sqq}$ (which is slightly renormalized here) such that 
		\begin{enumerate}
			\item  $(1,1)=1$, and $(F_i, F_j)= \delta_{i,j}  \frac{1}{\sqq^{2}-1}$, for $i,j\in \I$.
			\item  $(x, yz)=(r(x), y\otimes z)$ for $x,y,z\in\U^-_{\sqq}$; it is understood that 
			$(a\otimes a', b\otimes b')=(a, b)(a', b')$.
		\end{enumerate}
		
		\underline{Step 2}. $\widetilde{\ch}(\ca)$ is a bialgebra with comultiplication $r'$ defined by
		$$
		[\![A]\!] \mapsto \sum_{[B],[C]}\sqq^{\langle B,C\rangle}\frac{|\Ext^1_\ca(B,C)_A|}{|\Hom_\ca(B,C)|}[\![B]\!]\otimes[\![C]\!].
		$$
		There is a non-degenerate bilinear pairing $(\cdot, \cdot)'$ on $\widetilde{\ch}(\ca)$ defined by
		\begin{align*}
			([\![M]\!],[\![N]\!])' :=\delta_{[M],[N]}\frac{1}{|\Aut(M)|}.
		\end{align*}
		Note that $([\![S_i]\!], [\![S_j]\!])'=\delta_{ij}\frac1{q-1} =\delta_{ij}\frac1{\sqq^2 -1}$. 
		
		\underline{Step 3}. We have $(\phi(x),\phi(y))' =(x,y)$ for any $x,y\in\U^-_{\sqq}$. 
		
		\underline{Step 4}. For any $x\in\ker \phi$, we have $(x,y)=0$ for any $y\in\U^-_{\sqq}$, and then $x=0$ since $(\cdot, \cdot)$ is non-degenerate. \end{proof}
	
	\begin{example} \label{rank1ex}
		Let $Q$ be the quiver consisting of a single vertex with no arrow. Then $\bfk Q=\bfk$. Denote by $S$ the (1-dimensional) simple $\bfk Q$-module. For any $l\geq1$, we shall write
		\[
		[l S] =[\underbrace{S\oplus \cdots \oplus S}_{l}],
		\qquad
		[S]^{*l} = \underbrace{[S]* \cdots * [S]}_{l}, \qquad 
		[\![S]\!]^{*l} = \underbrace{[\![S]\!]* \cdots * [\![S]\!]}_{l}.
		\]
		One shows by induction on $l$ that 
		\begin{align}
			[S]^{*l}= \sqq^{-\frac{l(l-1)}{2}}[l S].  \label{eq:DPS}
		\end{align}
		Note that $|\aut(S^{\oplus l})|=(q-1)^l\sqq^{\frac{3l(l-1)}{2}}[l]_\sqq^!$ for any $l\geq1$. 
		So
		\begin{align*}
			[\![S]\!]^{*l}&= \sqq^{-\frac{l(l-1)}{2}}\frac{|\Aut(S^{\oplus l})|}{(q-1)^l}[\![l S]\!]=\sqq^{l(l-1)}][l]_\sqq^! [\![l S]\!].
		\end{align*}
		and then
		\begin{align}
			\label{HallDP}
			[S]^{(l)}:=\frac{[S]^{*l}}{[l]^!_\sqq}=\sqq^{-\frac{l(l-1)}{2}}\frac{[l S]}{[l]^!_\sqq},\qquad [\![S]\!]^{(l)}:=\frac{[\![S]\!]^{*l}}{[l]^!_\sqq}=\sqq^{l(l-1)}[\![l S]\!].
		\end{align}
	\end{example}  
	
	Before tackling the general rank 2 quivers, let us first warm up with the simplest cases.
	\begin{example} \label{rankA2ex}
		\begin{enumerate}
			\item 
			Assume $Q$ is the $A_1\times A_1$ quiver with vertices labelled by $1$ and $2$ and no arrows. Then we have 
			\begin{align*}
				[S_1]*[S_2]=[S_1\oplus S_2]=[S_2]*[S_1],
			\end{align*}
			i.e., $[S_1]$ and $[S_2]$ commute with each other, since 
			$$\Hom_{\bfk Q}(S_1,S_2)=0=\Hom_{\bfk Q}(S_2,S_1),\qquad \Ext^1_{\bfk Q}(S_1,S_2)=0=\Ext^1_{\bfk Q}(S_2,S_1).$$
			
			\item
			Assume $Q$ is the $A_2$ quiver with vertices labelled by $1$ and $2$ and one arrow from $1$ to $2$.
			Then $\dim_\bfk\Ext^1_{\bfk Q}(S_1,S_2)=1$, and
			$\dim_\bfk\Ext^1_{\bfk Q}(S_2,S_1)=0$. 
			Let $I_2$ be the indecomposable injective module which fits into the following short exact sequence of $\bfk Q$-modules:
			\begin{align*}
				0 \longrightarrow S_2 \longrightarrow I_2 \longrightarrow S_1 \longrightarrow 0.
			\end{align*}
			Note $I_2$ is also projective (i.e., it is the projective cover of $S_1$). Then we have
			\begin{align*}
				[S_1]*[S_2] &= \sqq^{-1}[S_1\oplus S_2] +\sqq^{-1}(q-1)[I_2],
				\\
				[S_2]*[S_1]&= [S_1\oplus S_2],
			\end{align*}
			and thus
			\begin{align*}
				\sqq [S_1]*[S_2]-[S_2]*[S_1]=(q-1)[I_2].
			\end{align*}
			We have 
			\begin{align*}
				[I_2]*[S_1] &=\sqq^{-1} [S_1\oplus I_2], \quad [S_1]*[I_2]=[S_1\oplus I_2];
				\\
				[I_2]*[S_2] &=[S_2\oplus I_2], \quad
				[S_2]*[I_2]=\sqq^{-1}[S_2\oplus I_2], 
			\end{align*}
			by noting that
			\begin{align*}
				\dim_\bfk \Hom_{\bfk Q}(I_2,S_1)& =1=\dim_\bfk\Hom_{\bfk Q}(S_2,I_2), 
				\\
				\dim_\bfk \Hom_{\bfk Q}(I_2,S_2)& =0=\dim_\bfk \Hom_{\bfk Q}(S_1, I_2).
			\end{align*}

			Therefore we have
			\begin{align*}
				&[\![S_1]\!]*[\![S_1]\!]*[\![S_2]\!]-(\sqq+\sqq^{-1})[\![S_1]\!]*[\![S_2]\!]*[\![S_1]\!] +[\![S_2]\!]*[\![S_1]\!]*[\![S_1]\!]
				\\
				&= \frac{1}{(q-1)^3}\Big([S_1]*\big( [S_1]*[S_2]-\sqq^{-1}[S_2]*[S_1]\big)- \big(\sqq [S_1]*[S_2]-[S_2]*[S_1]\big)*[S_1]  \Big)  
				\\
				&= \frac{1}{(q-1)^2}\big(\sqq^{-1}[S_1]*[I_2]-[I_2]*[S_1]\big)=0.
			\end{align*}
			
			An entirely similar computation shows that 
			\begin{align*}
				[\![S_2]\!]*[\![S_2]\!]*[\![S_1]\!]-(\sqq+\sqq^{-1})[\![S_2]\!]*[\![S_1]\!]*[\![S_2]\!] +[\![S_1]\!]*[\![S_2]\!]*[\![S_2]\!]
				=0.
			\end{align*}
			That is, $[\![S_1]\!]$ and $[\![S_2]\!]$ satisfy the quantum Serre relation of type $A_2$. 
		\end{enumerate}
	\end{example}  
	
	Now we treat the general rank 2 quivers.
	\begin{example} [General rank 2 quivers]
		\label{rank2ex}
		Consider a rank 2 quiver $Q$ with $a$ arrows from vertex $1$ to vertex $2$, and $b$ arrows from $2$ to $1$:
		\begin{align*}
			Q=(\xymatrix{1\ar@<1ex>[r]|-{a} & 2\ar@<1ex>[l]|-{b}}),
			\qquad \text{where } a+b=-c_{12}.
		\end{align*}
		
		For any $\bfk Q$-module $M=(M_i, M(\alpha_j), M(\beta_k), )_{i=1,2;1\leq j\leq a,1\leq k\leq b}$ such that $\dim_\bfk M_2=1$, we define
		\begin{align}
			\label{eq:UW}
			U_M:= \bigcap_{1\leq i\leq a} \Ker M(\alpha_i),\qquad W_M:= \sum_{j=1}^b \Im M(\beta_j),
		\end{align}
		and let
		\begin{align}  \label{eq:uw}
			u_M:=\dim_\bfk U_M,\qquad w_M:=\dim_\bfk W_M.
		\end{align}
		For any $r\geq0$, define
		\begin{align}
			\cI_{r}:&= \{[M]\in\Iso(\mod(\bfk Q))\mid \widehat{M}=r\widehat{S_1}+\widehat{S_2}, W_M\subseteq U_M\},
			\\
			\cJ_r:&= \{[M]\in\Iso(\mod(\bfk Q))\mid \widehat{M}=r\widehat{S_1}+\widehat{S_2}, S_2\in\add \Top(M)\}.
		\end{align}
		
		We shall write $[s S_1] = [S_1^{\oplus{s}}].$
		For any $s,t\geq0$, we have
		\begin{align*}
			[s S_1]*[S_2]*[t S_1]
			&= \sqq^{-tb} \sum\limits_{[U]\in\cJ_t} |\Ext^1(S_2,S_1^{\oplus{t}})_U| [S_1^{\oplus{s}}]*[U]
			\\
			&= \sqq^{-tb} \sum\limits_{[U]\in\cJ_t} |\Ext^1(S_2,S_1^{\oplus{t}})_U|  \sum\limits_{[M]\in\cI_{s+t}}\sqq^{sa-st}q^{s(t-a)}
			\frac{|\Ext^1(S_1^{\oplus{s}},U)_M|}{|\Hom(S_1^{\oplus{s}},U)|} \cdot [M].
		\end{align*}
		It is clear that $F_{S_2,S_1^{\oplus t}}^U\neq0$ if and only if $[U]\in\cJ_t$; in this case, $F_{S_2,S_1^{\oplus t}}^U=1$. By the Riedtmann-Peng formula we have
		\begin{align*}
			|\Ext^1(S_2,S_1^{\oplus{t}})_U | 
			&= \frac{|\Aut(S_2)|\cdot|\Aut(S_1^{\oplus t})|}{|\Aut(U)|} F_{S_2,S_1^{\oplus t}}^U= \frac{|\Aut(S_2)|\cdot|\Aut(S_1^{\oplus t})|}{|\Aut(U)|},
			\\
			\lvert \Ext^1(S_1^{\oplus{s}},U)_M \rvert 
			&= F_{S_1^{\oplus{s}},U}^{M}\cdot \lvert \Hom(S_1^{\oplus{s}},U)\rvert
			\frac{|\Aut{(S_1^{\oplus{s}})}|\cdot|\Aut(U)|}{|\Aut(M)|}.
		\end{align*}
		Therefore we have 
		\begin{align*}
			[s S_1]*[S_2]*[t S_1]
			&= \sqq^{st-sa-tb}(q-1) \sum_{[M]\in\cI_{s+t}} |\Aut(S_1^{\oplus s})||\Aut(S_1^{\oplus t})| \big( \sum_{[U]\in\cJ_t} F_{S_1^{\oplus{s}},U}^{M} \big) [\![M]\!].
		\end{align*}
		A simple computation shows that
		\begin{align*}
			\sum_{[U]\in\cJ_t} F_{S_1^{\oplus{s}},U}^{M}
			&= |\{f: S_1^{\oplus{s}}\hookrightarrow M\}|\cdot \frac{1}{|\Aut(S_1^{\oplus{s}})|}
			\\
			&= |\Gr(t-w_M,u_M-w_M)|\\
			&= \sqq^{(u_M-t)(t-w_M)} \qbinom{u_M-w_M}{t-w_M}_\sqq.
		\end{align*}
		Hence, we can rewrite the previous formula as 
		\begin{align}
			\label{HallDP2}
			[\![s S_1]\!]*[\![S_2]\!]*[\![t S_1]\!]
			&= \sqq^{st-sa-tb} \sum_{[M]\in\cI_{s+t}}   \sqq^{(u_M-t)(t-w_M)} \qbinom{u_M-w_M}{t-w_M}_\sqq [\![M]\!].
		\end{align}
		Then by \eqref{HallDP} and \eqref{HallDP2}, we have
		\begin{align*}
			&\sum_{t=0}^{a+b+1} (-1)^t	[\![ S_1]\!]^{(a+b+1-t)}*[\![S_2]\!]*[\![ S_1]\!]^{(t)}
			\\
			&= \sum_{[M]\in\cI_{a+b+1}} \sqq^{(a+b+1)b-u_Mw_M} \sum_{t=0}^{a+b+1}(-1)^t  \sqq^{-(2b+1-u_M-w_M)t}\qbinom{u_M-w_M}{t-w_M}_\sqq [\![M]\!].
		\end{align*}	
		
		Since $u_M\geq b+1>w_M$	for any $[M]\in \cI_{a+b+1}$, we have $1-u_M-w_M\leq 2b+1-u_M-w_M\leq u_M+w_M-1$, and by a standard $v$-binomial identity (cf. \cite[1.3.1(c)]{Lus93})
		\begin{align*}
			\sum_{t=0}^{a+b+1}(-1)^t  \sqq^{-(2b+1-u_M-w_M)t}\qbinom{u_M-w_M}{t-w_M}_\sqq =0.
		\end{align*}
		Therefore, we have concluded that 
		\begin{align*}
			&\sum_{t=0}^{a+b+1} (-1)^t	[\![ S_1]\!]^{(a+b+1-t)}*[\![S_2]\!]*[\![ S_1]\!]^{(t)}=0.
		\end{align*}
		The analogous identity with $[\![ S_1]\!], [\![ S_2]\!]$ switched holds since the above identity is symmetric with respect to $a,b$. 
	\end{example}

	\subsection{Bridgeland-Hall algebras and quantum groups}
	\label{Bridgeland}
	
	Let $\ca$ be a $\bfk$-linear exact category with finite-dimensional Hom-spaces and $\Ext^1$-spaces. 
	Let $\cc_{\Z_2}(\ca)$ be the exact category of $\Z_2$-graded complexes over $\ca$. Namely, an object $M$ of this category is a diagram with objects and morphisms in $\ca$:
	$$\xymatrix{ M^0 \ar@<0.5ex>[r]^{d^0}& M^1 \ar@<0.5ex>[l]^{d^1}  },\quad d^1d^0=d^0d^1=0.$$
	All indices of components of $\Z_2$-graded objects will be understood modulo $2$.
	A {morphism} $s=(s^0,s^1):M\rightarrow N$ is a diagram
	\[\xymatrix{  M^0 \ar@<0.5ex>[r]^{d^0} \ar[d]^{s^0}& M^1 \ar@<0.5ex>[l]^{d^1} \ar[d]^{s^1} \\
		N^0 \ar@<0.5ex>[r]^{e^0}& N^1 \ar@<0.5ex>[l]^{e^1} }  \]
	with $s^{i+1}d^i=e^is^i$.
	
	The shift functor on complexes is an involution $*$ on $\cc_{\Z_2}(\ca)$
	which shifts the grading and changes the sign of the differential as follows:
	$$\xymatrix{ \big(M^0 \ar@<0.5ex>[r]^{d^0}& M^1 \big) \ar@<0.5ex>[l]^{d^1} \ar[r]^{*} & \big(M^1 \ar@<0.5ex>[r]^{-d^1} & M^0\big) \ar@<0.5ex>[l]^{-d^0}  }.$$
	
	Similar to ordinary complexes, we can define the $i$th homology group for $M$, denoted by $H^i(M)$, for $i\in\Z_2$. A complex is called acyclic if its homology group is zero. The subcategory of all acyclic $\Z_2$-graded complexes over $\ca$ is denoted by $\cc_{\Z_2,ac}(\ca)$.
	
	For any object $X\in\ca$, we define
	\begin{align}
		\label{stalks}
		\begin{split}
			K_X& := \big(\xymatrix{ X \ar@<0.5ex>[r]^{1}& X \ar@<0.5ex>[l]^{0}  }\big),\qquad \,\, K_X^*:= \big(\xymatrix{ X \ar@<0.5ex>[r]^{0}& X \ar@<0.5ex>[l]^{1}  } \big),
			\\
			C_X& := \big(\xymatrix{ 0 \ar@<0.5ex>[r]& X \ar@<0.5ex>[l]  }\big),\qquad \quad C_X^*:=\big(\xymatrix{ X\ar@<0.5ex>[r]& 0 \ar@<0.5ex>[l]  }\big)
		\end{split}
	\end{align}
	in $\cc_{\Z_2}(\ca)$.  Note that $K_X,K_X^*$ are acyclic complexes.	Let $\res: \cc_{\Z_2}(\ca)\rightarrow \ca\times \ca$ be the forgetful functor, which maps $M=\big(\xymatrix{ M^0 \ar@<0.5ex>[r]^{d^0}& M^1 \ar@<0.5ex>[l]^{d^1}  }\big)$ to $(M^0, M^{1})$.
	
	From now on, we shall impose a stronger assumption that $\ca$ is a hereditary abelian category.
	
	\begin{lemma} [\text{\cite[Proposition 2.3]{LP21}}]
		\label{proposition extension 2 zero}
		Let $\ca$ be a hereditary abelian category.
		For any $K\in\cc_{\Z_2}(\ca)$ with $K$ acyclic, we have $
		\pdim_{\cc_{\Z_2}(\ca)} (K)\leq 1, \text{ and }{ \rm inj.dim}_{\cc_{\Z_2}(\ca)} (K)\leq 1.$
	\end{lemma}
	
	By Lemma \ref{proposition extension 2 zero}, for any $K,M\in \cc_{\Z_2}(\ca) $ with $K$ acyclic, the following pairings are well defined:
	\begin{align*}
		\langle K,M\rangle=\dim\Hom_{\cc_{\Z_2}({\ca})}(K,M)-\dim\Ext^1_{\cc_{\Z_2}({\ca})}(K,M)
	\end{align*}
	and
	\begin{align*}
		\langle M,K\rangle =\dim\Hom_{\cc_{\Z_2}({\ca})}(M,K)-\dim\Ext^1_{\cc_{\Z_2}({\ca})}(M,K).
	\end{align*}
	They induce bilinear pairings (call {\em Euler forms}) between the Grothendieck groups $K_0(\cc_{\Z_2,ac}(\ca))$ and $K_0\big(\cc_{\Z_2}(\ca) \big)$:
	$$\langle\cdot,\cdot\rangle:K_0(\cc_{\Z_2,ac}({\ca}))\times K_0\big(\cc_{\Z_2}(\ca) \big)\longrightarrow \Z,$$
	and
	$$\langle\cdot,\cdot\rangle:K_0\big(\cc_{\Z_2}(\ca) \big)\times K_0(\cc_{\Z_2,ac}({\ca}))\longrightarrow \Z.$$
	We can use the same notation above, since the two forms coincide on $K_0(\cc_{\Z_2,ac}({\ca}))\times K_0(\cc_{\Z_2,ac}({\ca}))$.

	Let $K_0(\ca)$ be the Grothendieck group of $\ca$.
	For any $A\in\ca$, we denote by $\widehat{A}$ the corresponding element in the Grothendieck group $K_0(\ca)$.
	We also use $\langle \cdot,\cdot\rangle$ to denote the Euler form of $\ca$, i.e.,
	\begin{align*}
		\langle \widehat{A}, \widehat{B}\rangle=\dim\Hom_\ca(A,B)-\dim\Ext^1_\ca(A,B), \text{ for any }A,B\in\ca.
	\end{align*}
	Let $(\cdot,\cdot)$ be the symmetrized Euler form of $\ca$, i.e., $(\widehat{A},\widehat{B})= \langle \widehat{A}, \widehat{B}\rangle+\langle \widehat{B}, \widehat{A}\rangle$.
	
	
	\begin{proposition}[\text{\cite[Proposition 2.4, Corollary 2.5]{LP21}}]
		\label{lema euler form}
		For any $A ,B \in\ca$, we have 
		\begin{align}
			&\langle C_A,K_B\rangle=\langle C^*_A,K_B^*\rangle=\langle \widehat{A}, \widehat{B}\rangle,\quad
			\langle K_B,C^*_A\rangle=\langle K_B^*,C_A\rangle=\langle \widehat{B},\widehat{A}\rangle;
			\\
			&\langle K_B,C_A\rangle=\langle C_A^*,K_B\rangle=\langle C_A,K_B^*\rangle=\langle K_B^*,C_A^*\rangle=0;
			\\
			&\langle K_{A}, K_{B}\rangle= \langle K_{A}^*, K_{B}^*\rangle=\langle K_{A}, K_{B}^*\rangle=\langle K_{A}^*, K_{B}\rangle=\langle \widehat{A},\widehat{B}\rangle.
		\end{align}
	\end{proposition}
	
	In the remaining of this subsection, we shall assume that $\ca$ is a hereditary abelian category with enough projective objects.	Let $\cp$ be the full subcategory of $\ca$ formed by projective objects. 
	Let $\cc_{\Z_2}(\cp)$ be the category of $\Z_2$-graded complexes over $\cp$.  Then $K_P$, $K_P^*$ are projective objects in $\cc_{\Z_2}(\ca)$, and $\cc_{\Z_2}(\cp)$ is a Frobenius category with $K_P$, $K_P^*$ ($P\in\cp$) as projective-injective objects; see \S\ref{subsec:1-Gorenstein} for the definition of Frobenius categories. 
	For any $X\in\ca$, let $0\rightarrow P_X\stackrel{f}{\rightarrow} Q_X\rightarrow X\rightarrow0$ be the minimal projective resolution of $X$.
	Define 
	\begin{align}
		\label{GX}
		G_X &:=\big(\xymatrix{ P_X \ar@<0.5ex>[r]^{f}& Q_X \ar@<0.5ex>[l]^{0}  }\big),\qquad \,\, G_X^*:=\big(\xymatrix{ Q_X \ar@<0.5ex>[r]^{0}& P_X \ar@<0.5ex>[l]^{f}  }\big).
	\end{align}

	We consider $\widetilde{\ch}(\cc_{\Z_2}(\cp))$, the twisted Ringel-Hall algebra of $\cc_{\Z_2}(\cp)$, i.e., 
	\begin{align}
		\widetilde{\ch}(\cc_{\Z_2}(\cp))=\bigoplus_{[M]\in\Iso(\cc_{\Z_2}(\cp))}\Q(\sqq)[M],
	\end{align}
	with multiplication given by
	\begin{align*}
		[L]*[M]=\sqq^{\langle \res L,\res M\rangle} \sum_{[N]} \frac{|\Ext^1(L,M)_N|}{|\Hom(L,M)|}[N].
	\end{align*}
	
	Since $K_P,K_P^*\in\cc_{\Z_2}(\cp)$ are projective-injective objects, we have
	\begin{align*}
		[K_P]*[M]&= \sqq^{\langle P,M^0\oplus M^1\rangle}q^{-\langle K_P,M\rangle}[K_P\oplus M]\\
		&= \sqq^{\langle P,M^0\rangle-\langle P,M^1\rangle}[K_P\oplus M],
	\end{align*}
	and similarly
	\begin{align*}
		[M]*[K_P]= \sqq^{\langle M^1,P\rangle-\langle M^0,P\rangle}[K_P\oplus M].
	\end{align*}
	Therefore, we have
	\begin{align}
		\label{eq:KM}
		[K_P]*[M]=\sqq^{(P,M^0)-(P,M^1)}[M]*[K_P].
	\end{align}
	So one can show that $\cs_\cp:=\{a[K_P]\mid a\in\Q(\sqq)^\times, P\in\cp\}$ is a multiplicatively closed subset with the identity $[0]$, which is a right Ore, right reversible subset. 
	
	By definition, the {\em Bridgeland-Hall algebra}, denoted by
	$D\widetilde{\ch}(\ca)$, is the right localization of $\widetilde{\ch}(\cc_{\Z_2}(\cp))$ with respect to $\cs_\cp$.
	
	Denote by 
	$$E_X=\sqq^{\langle P,X\rangle}[K_{P_X}]^{-1}*[G_X],\quad F_X=E_X^*,$$
	and
	$$[K_X]:=[K_{Q_X}]*[K_{P_X}]^{-1},\quad[K^*_X]:=[K^*_{Q_X}]*[K^*_{P_X}]^{-1} ,$$
	where $0\rightarrow P_X\rightarrow Q_X\rightarrow X\rightarrow0$ is a projective resolution of $X$.

	Note that $E_X$, $F_X$, $[K_X]$, $[K_X^*]$ do not depend on the projective resolutions of $X$  (not necessarily minimal).
	In fact, $[K_X]$ and $[K_X^*]$ only depend on the class $\widehat{X}\in K_0(\ca)$. So for any $\alpha\in K_0(\ca)$, one can define $[K_\alpha],[K_\alpha^*]$ in $D\widetilde{\ch}(\ca)$. 
	
	In case $\ca=\mod(\bfk Q)$ for an acyclic quiver $Q$, we shall denote its Bridgeland-Hall algebra by $D\widetilde{\ch}(\bfk Q)$. We shall also consider the base change $\tU_{\sqq} =\Q(\sqq) \otimes_{\Q(v)} \tU.$
	
	\begin{theorem}[\cite{Br}]
		\label{thm:Bridgeland}
		There is an injective homomorphism of $\Q(\sqq)$-algebras
		$\Psi:\tU_{\sqq} \longrightarrow D\widetilde{\ch}(\bfk Q)$ such that 
		\begin{align}
			\Psi(E_i)=\frac{\sqq}{q-1}E_{S_i},\quad \Psi(F_i)=\frac{-1}{q-1}F_{S_i},\quad \Psi(K_i)=[K_{S_i}],\quad \Psi(K_i')=[K_{S_i}^*],\quad \forall i\in \I.
		\end{align}
	\end{theorem}

	\begin{proof}
		We give a sketched proof here. 
		First, to show that $\Psi$ is a homomorphism, we prove that $\Psi$ 
		preserves the defining relations for $\tU$.
		
		For the Serre relations, we only prove \eqref{eq:serre2} for $F$'s, as the same proof applies for the Serre relations for $E$'s. With the help of Theorem \ref{thm:Ringel-Green}, we only need to check that \begin{align*}
			I^-:\widetilde{\ch}(\bfk Q)\longrightarrow D\widetilde{\ch}(\bfk Q),\quad [M]\mapsto F_M
		\end{align*} 
		is an algebra homomorphism. 
		
		For any $X,Z\in\mod(\bfk Q)$, we have two (minimal) projective resolution:
		\begin{align*}
			0\longrightarrow P_X\longrightarrow Q_X\longrightarrow X\longrightarrow0,\qquad  0\longrightarrow P_Z\longrightarrow Q_Z\longrightarrow Z\longrightarrow0,
		\end{align*}
		which give rise to two exact sequences:
		\begin{align}
			\label{eq:GPresol}
			0\longrightarrow K_{P_X}\rightarrow G_X\longrightarrow C_X\longrightarrow0,\qquad   0\longrightarrow K_{P_Z}\rightarrow G_Z\longrightarrow C_Z\longrightarrow0.
		\end{align}
		
		We make the following
		\vspace{2mm}
		
		{\bf Claim.} We have $\Ext^1_\ca(X,Z)\cong \Ext^1_{\cc_{\Z_2}(\ca)}(C_X,C_Z)\cong \Ext^1_{\cc_{\Z_2}(\ca)}(G_X,G_Z)$. 
		
		Let us prove the Claim. 
		Using \eqref{eq:GPresol}, we have the following exact sequences:
		\begin{align*}
			0\longrightarrow \Hom_{\cc_{\Z_2}(\ca)}& (K_{P_X},C_Z)\longrightarrow \Hom_{\cc_{\Z_2}(\ca)}(C_X,C_Z)\longrightarrow \Hom_{\cc_{\Z_2}(\ca)}(G_X,C_Z)\longrightarrow0,
			\\
			&0\longrightarrow	\Ext^1_{\cc_{\Z_2}(\ca)}(G_X,G_Z)\longrightarrow\Ext^1_{\cc_{\Z_2}(\ca)}(G_X,C_Z)\longrightarrow0,
		\end{align*}
		since $K_{P_X}$ is projective in $\cc_{\Z_2}(\ca)$, and $K_{P_Z}$ is projective-injective in $\cc_{\Z_2}(\cp)$. It follows that $ \Ext^1_{\cc_{\Z_2}(\ca)}(C_X,C_Z)\cong \Ext^1_{\cc_{\Z_2}(\ca)}(G_X,G_Z)$ by  $\Hom_{\cc_{\Z_2}(\ca)}(K_{P_X},C_Z)=0$. The Claim follows.
		
		It follows that $I^-([X]*[Z])=I^-([X])*I^-([Z])$ by definition and the Claim above. Hence the Serre relations are preserved by $\Psi$.

		Similarly, one checks that
		$$\Ext^1_{\cc_{\Z_2}(\ca)}(G_X,G_Z^*)\cong \Hom_\ca(X,Z)\cong \Ext^1_{\cc_{\Z_2}(\ca)}(G_X^*,G_Z).$$
		So we have $[E_{S_i},F_{S_j}]=0$, for any $i\neq j$. 
		
		On the other hand, we have
		$\Ext^1_{\cc_{\Z_2}(\ca)}(G_{S_i},G_{S_i}^*)\cong \bfk$, for any $i\in \I$. If an exact sequence $0\rightarrow G_{S_i}^*\rightarrow M\rightarrow G_{S_i}\rightarrow0$ is non-split, then $M$ is acyclic and $M\cong K_{P_{S_i}}^*\oplus K_{Q_{S_i}}^*$. So we have  $[E_{S_i},F_{S_i}]=(q-1)([K_{{S_i}}^*]-[K_{{S_i}}])$, as desired.
		
		It remains to prove the injectivity of $\Psi$. Any $M\in\cc_{\Z_2}(\ca)$ is isomorphic to $C_A\oplus C_B^*\oplus K_P\oplus K_Q^*$; see \cite[Lemma~ 4.2]{Br}. So $D\widetilde{\ch}(\bfk Q)$ has a basis consisting of elements
		$$[C_A\oplus C_B^*]*[K_\alpha]*[K_\beta^*], \quad \forall [A],[B]\in\Iso(\mod(\bfk Q)), \alpha,\beta\in K_0(\mod(\bfk Q)).$$
		It implies that $D\widetilde{\ch}(\bfk Q)$ has a basis consisting of elements
		$$[C_A] *[C_B^*]*[K_\alpha]*[K_\beta^*], \quad \forall [A],[B]\in\Iso(\mod(\bfk Q)), \alpha,\beta\in K_0(\mod(\bfk Q));$$
		see \cite[Lemma 4.7]{Br} or \cite[Theorem 3.25]{LP21}. Then the injectivity of $\Psi$ follows from this and Theorem~ \ref{thm:Ringel-Green}. 
	\end{proof}

	\section{Semi-derived Ringel-Hall algebras}
	\label{sec:Semi-derivedHall}
	
	In this section, we formulate weakly $1$-Gorenstein exact categories and provide several examples. We then introduce the semi-derived Ringel-Hall algebras for such categories.

	\subsection{Weakly $1$-Gorenstein exact categories}
	\label{subsec:1-Gorenstein}
	
	Let $\ce$ be an essentially small exact category in the sense of Quillen,
	linear over a finite field $\bfk=\F_q$.
	For an exact category $\ce$, we introduce the following subcategories of $\ca$:
	\begin{align*}
		\cp^{\leq i}(\ce) &=\{X\in\ce\mid\pd X\leq i\},
		\\
		{\mathcal I}^{\leq i}(\ce) &=\{X\in\ce\mid\ind X\leq i\}, \quad \forall i\in\N,
		\\
		\cp^{<\infty}(\ce) &= \{X\in\ce\mid\pd X<\infty\},
		\\
		{\mathcal I}^{<\infty}(\ce) &= \{X\in\ce\mid\ind X<\infty\}.
	\end{align*}
	The category $\ce$ is called \emph{weakly Gorenstein} if $\cp^{<\infty}(\ce)={\mathcal I}^{<\infty}(\ce)$, and $\ce$ is a \emph{weakly $d$-Gorenstein} exact category if $\ca$ is weakly Gorenstein and $\cp^{<\infty}(\ce)=\cp^{\leq d}(\ce)$, ${\mathcal I}^{<\infty}(\ce)={\mathcal I}^{\leq d}(\ce)$.
	
	\begin{lemma}[Iwanaga's Theorem]
		Let $\ce$ be a weakly Gorenstein exact category with enough projectives and injectives. Then $\cp^{<\infty}(\ce)=\cp^{\leq d}(\ce)$ if and only if ${\mathcal I}^{<\infty}(\ce)={\mathcal I}^{\leq d}(\ce)$.
	\end{lemma}
	
	Throughout this section, we always assume that $\ce$ is an exact category satisfying the following conditions:
	\begin{itemize}
		\item[(Ea)] $\ce$ is essentially small, with finite morphism spaces, and finite extension spaces,
		\item[(Eb)] $\ce$ is linear over $\bfk=\F_q$,
		\item[(Ec)] $\ce$ is weakly $1$-Gorenstein,
		\item[(Ed)] For any object $X\in\ca$, there exists an object $P_X\in\cp^{<\infty}(\ce)$ and a deflation $P_X\twoheadrightarrow X$.
	\end{itemize}
	In this case, we have  $\cp^{<\infty}(\ce)=\cp^{\leq1}(\ce)={\mathcal I}^{<\infty}(\ce)={\mathcal I}^{\leq 1}(\ce)$.

	\begin{example}
		\begin{enumerate}
			\item 
			Any $\bfk$-linear hereditary abelian category with finite-dimensional morphism spaces and extension spaces satisfies (Ea)-(Ed).
			\item
			For any finite-dimensional $1$-Gorenstein algebra $\Lambda$ over $\bfk$ (i.e., $\ind{}_\Lambda \Lambda\leq1$ and $\ind \Lambda_\Lambda\leq1$), $\ce=\mod(\Lambda)$ satisfies (Ea)--(Ed).
		\end{enumerate}
	\end{example}
	
	An exact category $\ce$ is called a {\em Frobenius category}, provided that the class of projective objects coincides with the class of injective objects and it has enough projective objects. The importance of Frobenius categories lies in that they give rise naturally to triangulated categories; see \cite{Ha88}.
	
	\begin{example}
		Any $\bfk$-linear Frobenius category   with finite-dimensional morphism spaces and extension spaces satisfies (Ea)-(Ed).
		
		Let $A$ be a finite-dimensional algebra over $\bfk$. The notion of Gorenstein projective modules is introduced in \cite{EJ}. Let $\Gproj(A)$ be the category consisting of all finitely generated Gorenstein projective modules. Then $\Gproj(A)$ is a Frobenius category satisfies (Ea)-(Ed).
	\end{example}
	
	\begin{example}\label{example 1}
		Let $\ca$ be a hereditary abelian $\bfk$-category (not necessarily with enough projective objects). Let $\cc_{\Z_n}(\ca)$ be the category of $\Z_n$-graded complexes (also called $n$-periodic complexes), for $n\geq1$. Denote by $\cc_{ac,\Z_n}(\ca)$ the subcategory of acyclic complexes in
		$\cc_{\Z_n}(\ca)$. Then  $\cc_{\Z_n}(\ca)$ satisfies (Ea)-(Ed) with $\cp^{<\infty}(\cc_{\Z_n}(\ca))=\cc_{ac,\Z_n}(\ca)=\mathcal{I}^{<\infty}(\cc_{\Z_n}(\ca))$; see  \cite[Proposition 2.3]{LP21} and \cite[Corollary 2.4]{LRW20}.
	\end{example}
	
	By definition, the {\em singularity category} $\cd_{sg}(\ce)$ is the Verdier quotient of $\cd^b(\ce)$ with respect to the thick subcategory generated by $\cp^{<\infty}(\ce)$. 
	
	\subsection{Definition of semi-derived Ringel-Hall algebra}

	Given an exact category $\ce$ satisfying (Ea)-(Ed), its Ringel-Hall algebra $\ch(\ce)$ can be too big to be interesting. We apply some suitable quotient and localization constructions below to $\ch(\ce)$ to obtain a new algebra of suitable size, which has a natural (categorical) basis. As we shall see, for some distinguished classes of exact categories, the resulting algebras provide realizations of quantum groups and $\imath$quantum groups.
	
	Let $I$ be the two-sided ideal of $\ch(\ce)$ generated by all differences 
	\begin{align}
		\label{def:I}
		[M]-[N], \text{ if }\widehat{M}=\widehat{N}\in K_0(\ce), \text{ and }M\cong N \text{ in }\cd_{sg}(\ce).
	\end{align}
	Denote
	\begin{equation}
		\label{eq:Sca}
		\cs_{\ce} := \{ a[K] \in \ch(\ce)/I \mid a\in \Q^\times, K\in \cp^{\leq1}(\ce)\}.
	\end{equation}
	
	\begin{proposition}[\cite{LP21,Lu22}]\label{prop:Def of MRH}
		Let $\ce$ be an exact category satisfying (Ea)-(Ed). Then the multiplicatively closed subset $\cs_{\ce}$ is a right Ore, right reversible subset of $\ch(\ce)/I$. Equivalently, the right localization of
		$\ch(\ce)/I$ with respect to $\cs_{\ce}$ exists, and will be denoted by $(\ch(\ce)/I)[\cs_{\ce}^{-1}]$.
	\end{proposition}

	\begin{definition}[\cite{LP21,Lu22}]
		For any exact category $\ce$ satisfying (Ea)-(Ed), $(\ch(\ce)/I)[\cs_{\ce}^{-1}]$ is called the semi-derived Ringel-Hall algebra of $\ce$, and denoted by $\cs\cd\ch(\ce)$.
		
		The quantum torus $\ct(\ce)$ is defined to be the subalgebra of $\cs\cd\ch(\ce)$ generated by $[M]$, for $M$ in $\fpr(\ce)$.
	\end{definition}

	%
	
	\begin{example}
		Let $\ca$ be a $\bfk$-linear hereditary abelian category with finite-dimensional morphism spaces and extension spaces. Then
		$\cs\cd\ch(\ca)$ is the group algebra of $K_0(\ca)$ with its multiplication twisted by the Euler form, i.e.,
		$$[M]\diamond[N]=q^{\langle M,N\rangle}[M\oplus N],$$
		for any $M,N\in\ca$. 
	\end{example}

	\begin{example}
		Let $\cf$ be a Frobenius category with $\cp(\cf)$ the subcategory of projective-injective objects. The semi-derived Ringel-Hall algebra of $\cf$ was introduced by M.~Gorsky \cite{Gor2}.  
	\end{example}

	\begin{example}
		Let $\Lambda$ be a finite-dimensional $1$-Gorenstein algebra. Then the natural embedding $\Gproj(\Lambda)\subseteq \mod(\Lambda)$ induces an isomorphism \begin{align*}
			F:\cs\cd\ch \big(\Gproj(\Lambda)\big) \stackrel{\cong}{\longrightarrow} \cs\cd\ch(\Lambda). 
		\end{align*}
		Moreover, for any $\Lambda$-module $M$, we have an exact sequence
		$$0\longrightarrow P_M\longrightarrow G_M\longrightarrow M\longrightarrow0,$$
		where $G_M\in\Gproj(\Lambda)$ and $P_M$ is projective. 
		Then the inverse of $F$ is given by $[M] \mapsto q^{-\langle M, P_M \rangle}[G_M]\diamond [P_M]^{-1}$; see \cite[Theorem A.15]{Lu22}.
	\end{example}

	\begin{example}
		For a $\bfk$-linear hereditary abelian category $\ca$, let $\cc_{\Z_n}(\ca)$ be the category of $\Z_n$-graded complexes of modules over $\ca$, for $n\geq1$. Then $\cs\cd\ch(\cc_{\Z_n}(\ca))$ was studied in \cite{LP21, LinP19}. 
	\end{example}

	\subsection{Bridgeland's Theorem reformulated} 
	

	Let $\ca$ be a hereditary abelian $\bfk$-category. Recall that the category $\cc_{\Z_2}(\ca)$ is defined in \S\ref{Bridgeland}. 
	Let $\cs\cd\widetilde{\ch}\big(\cc_{\Z_2}(\ca) \big)$ be the twisted semi-derived Ringel-Hall algebra of $\cc_{\Z_2}(\ca)$ over $\Q(\sqq)$, i.e., its multiplication is given by
	\begin{align}
		\label{eq:twistedprod}
		[M]*[N]:=\sqq^{\langle \res M,\res N\rangle}[M]\diamond[N],\quad \forall M,N\in\cc_{\Z_2}(\ca).
	\end{align}
	
	Let $\cp$ be the subcategory of $\ca$ consisting of projective objects. If $\ca$ admits enough projectives, then the natural embedding $\cc_{\Z_2}(\cp)\subseteq \cc_{\Z_2}(\ca)$ induces an algebra isomorphism (cf. \cite[Theorem 5.13]{LP21})
	\begin{align*}
		\phi:D\widetilde{\ch}(\ca) \stackrel{\cong}{\longrightarrow}  \cs\cd\widetilde{\ch}\big(\cc_{\Z_2}(\ca) \big).
	\end{align*} 
	In fact, $D\widetilde{\ch}(\ca)$ is just the twisted semi-derived Ringel-Hall algebra of $\cc_{\Z_2}(\cp)$ by noting that 
	$\cc_{\Z_2}(\cp)$ is a Frobenius category.
	
	We further specialize to the distinguished case when $\ca=\mod(\bfk Q)$ for a finite quiver $Q$ (not necessarily acyclic); recall $\mod(\bfk Q)$ may not have enough projectives. Denote by $\cs\cd\widetilde{\ch}_{\Z_2}(\bfk Q)$ the twisted semi-derived Ringel-Hall algebra of $\cc_{\Z_2}(\mod(\bfk Q))$.
	
	\begin{theorem}[Bridgeland's Theorem reformulated; see \cite{LW22}]
		There exists an injective homomorphism of $\Q(\sqq)$-algebras
		$\Psi:\tU_{\sqq} \longrightarrow \cs\cd\widetilde{\ch}_{\Z_2}(\bfk Q)$ such that 
		\begin{align}
			\Psi(E_i)=\frac{\sqq}{q-1}[C_{S_i}],\quad \Psi(F_i)=\frac{-1}{q-1}[C^*_{S_i}],\quad \Psi(K_i)=[K_{S_i}],\quad \Psi(K_i')=[K_{S_i}^*],\quad \forall i\in \I.
		\end{align}
	\end{theorem}

	\begin{proof}
		We remark that the verification that $\Psi$ is an algebra homomorphism is much easier than in the setting of Theorem~ \ref{thm:Bridgeland}.
		
		Let $\ca=\mod(\bfk Q)$.
		The natural embeddings $\iota: \ca\rightarrow \cc_{\Z_2}(\ca)$ given by $X\mapsto C_X$ makes $\ca$ a full subcategory of $\cc_{\Z_2}(\ca)$ which is closed under taking extensions. So $\iota$ induces an algebra homomorphism $\widetilde{\ch}(\ca) \rightarrow\cs\cd\widetilde{\ch}_{\Z_2}(\bfk Q)$ sending $[X]\mapsto [C_X]$. So $\Psi$ preserves the quantum Serre relation \eqref{eq:serre1} by Theorem \ref{thm:Ringel-Green}. Similarly, $\Psi$ preserves the other quantum Serre relation \eqref{eq:serre2}.
		
		Let $i\neq j$. Since  $\dim_\bfk\Ext^1_{\cc_{\Z_2}(\ca)}(C_{S_i},C_{S_j}^*)=\dim_\bfk \Hom_\ca(S_i,S_j)= 0$, we have $[C_{S_i},C_{S_j}^*]=0$. 
		
		For any non-split exact sequence
		$0\rightarrow C_{S_i}^*\rightarrow M\rightarrow C_{S_i} \rightarrow0$, we have
		$M\cong K_{S_i}^*$. So we have
		$[C_{S_i},C_{S_i}^*]=(q-1)([K_{S_i}^*]-[K_{S_i}])$.
		
		Finally, by Proposition~ \ref{lema euler form}, we see that the remaining relations \eqref{eq:KK}-\eqref{eq:K2} for $\tU$ are preserved by $\Psi$. 
	\end{proof}

	\section{$\imath$Quiver algebras}
	\label{sec:iQA}

	In this section, starting from quivers with involutions, we formulate a class of 1-Gorenstein algebras known as  $\imath$quiver algebras. In case of acyclic quivers, these algebras are finite-dimensional. We present some basic homological properties of these algebras. 

	\subsection{$\imath$Quivers}
	\label{subsec:i-quivers}
	
	Let $\bfk$ be a field.
	Let $Q=(Q_0,Q_1)$ be a quiver. 
	An {\em involution} of $Q$ is defined to be an automorphism $\btau$ of the quiver $Q$ such that $\btau^2=\Id$. In particular, we allow the {\em trivial} involution $\Id:Q\rightarrow Q$. An involution $\btau$ of $Q$ induces an involution of the path algebra $\bfk Q$, again denoted by $\btau$.
	A quiver together with a specified involution $\btau$, $(Q, \btau)$, will be called an {\em $\imath$quiver}. 	The {\em rank} of an $\imath$quiver $(Q, \btau)$ is by definition the number of $\btau$-orbits in $Q_0$.
	
	Let $Z_m$ be the quiver with $m$ vertices and $m$ arrows which forms an oriented cycle. The vertex set of $Z_m$ is $\{0,1,\dots,m-1\}$. Let $R_m$ be the radical square zero selfinjective Nakayama algebra of $Z_m$, i.e.,  $R_m:=\bfk Z_m/J$, where $J$ denotes the ideal of $\bfk Z_m$ generated by all paths of length two. In particular,
	
	$\triangleright$ $R_1$ is isomorphic to the truncated polynomial algebra $\bfk[\varepsilon]/(\varepsilon^2)$;
	
	$\triangleright$ $R_2$ is the radical square zero of the path algebra of $\xymatrix{1 \ar@<0.5ex>[r]^{\varepsilon} & 1' \ar@<0.5ex>[l]^{\varepsilon'}}$, i.e., $\varepsilon' \varepsilon =0 =\varepsilon\varepsilon '$.
	
	Recall that $\cc_{\Z_m}(\mod(\bfk Q))$ is the category of the $\Z_m$-graded complexes over $\mod(\bfk Q)$ for any $m\geq 1$, see \cite{Br, LP21}. 
	We have $\cc_{\Z_m}(\mod(\bfk Q))\cong \mod(\bfk Q\otimes_\bfk R_m)$ for any $m \geq 1$.
	
	Define a $\bfk$-algebra
	\begin{equation}
		\label{eq:La}
		\Lambda=\bfk Q\otimes_\bfk R_2.
	\end{equation}
	
	Let $Q^{\sharp}$ be the quiver such that
	\begin{itemize}
		\item the vertex set of $Q^{\sharp}$ consists of 2 copies of the vertex set $Q_0$, $\{i,i'\mid i\in Q_0\}$;
		\item the arrow set of $Q^{\sharp}$ is
		\[
		\{\alpha: i\rightarrow j,\alpha': i'\rightarrow j'\mid(\alpha:i\rightarrow j)\in Q_1\}\cup\{ \varepsilon_i: i\rightarrow i' ,\varepsilon'_i: i'\rightarrow i\mid i\in Q_0 \}.
		\]
	\end{itemize}
	We call $Q^{\sharp}$ the {\em double framed quiver} associated to the quiver $Q$.

	Let $I^{\sharp}$ be the  ideal of $\bfk Q^{\sharp}$ generated by
	\begin{itemize}
		\item
		(Nilpotent relations) $\varepsilon_i \varepsilon_i'$, $\varepsilon_i'\varepsilon_i$ for any $i\in Q_0$;
		\item
		(Commutative relations) $\varepsilon_j' \alpha' -\alpha\varepsilon_i'$, $\varepsilon_j \alpha -\alpha'\varepsilon_i$ for any $(\alpha:i\rightarrow j)\in Q_1$.
	\end{itemize}
	Then the algebra $\La$ can be realized as
	\begin{equation}  \label{eq:La=QI}
		\Lambda\cong \bfk Q^{\sharp} \big/ I^{\sharp}.
	\end{equation}
	Let $Q$ (respectively, $Q'$) be the full subquiver of $Q^{\sharp}$ formed by all vertices $i$ (respectively, $i'$) for $i\in Q_0$. Then $Q\sqcup Q'$ is a subquiver of $Q^{\sharp}$.
	
	The category $\mod(\Lambda)$ is isomorphic to $\cc_{\Z_2}(\mod(\bfk Q))$.
	
	\begin{example}
		\label{exDoubleQuiver}
		(a)  Let $Q=(\xymatrix{ 1\ar[r]^{\alpha} &2})$. Then its double framed quiver $Q^{\sharp}$ is
		\begin{center}\setlength{\unitlength}{0.7mm}
			
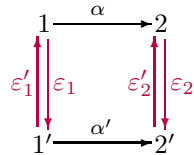
\vspace{-2cm}
			\begin{equation*}
				\begin{picture}(100,40)(0,20)
					\put(49,8){\small $1'$}
					\put(50,31){\small $1$}
					\put(72,8){\small $2'$}
					\put(72,31){\small $2$}
					
					\put(53,10){\vector(1,0){18.5}}
					\put(53,32.5){\vector(1,0){18.5}}
					
					\put(60,12.5){$_{\alpha'}$}
					\put(60,35){$_\alpha$}
					\color{purple}
					\put(50,13){\vector(0,1){17}}
					\put(52,29.5){\vector(0,-1){17}}
					\put(72,13){\vector(0,1){17}}
					\put(74,29.5){\vector(0,-1){17}}
					
					\put(45,20){\small $\varepsilon_1'$}
					\put(53,20){\small $\varepsilon_1$}
					\put(67,20){\small $\varepsilon_2'$}
					\put(75,20){\small $\varepsilon_2$}
				\end{picture}
			\end{equation*}
			\vspace{0.2cm}
		\end{center}
		and $I^{\sharp}$ is generated by all possible quadratic relations
		\begin{eqnarray*}
			&&\varepsilon_1\varepsilon_1', \,\,\,\,\varepsilon_1'\varepsilon_1, \,\,\,\,\varepsilon_2'\varepsilon_2, \,\,\,\,\varepsilon_2\varepsilon_2',\,\,\,\, \alpha' \varepsilon_1 -\varepsilon_2 \alpha,\,\,\,\, \alpha \varepsilon_1'- \varepsilon_2' \alpha'.
		\end{eqnarray*}
		
		(b) Let $Q=(\xymatrix{ 1\ar[r]^{\alpha} &2 & 3\ar[l]_{\beta} })$. Then its double framed quiver $Q^{\sharp}$ is
		\begin{center}\setlength{\unitlength}{0.7mm}
			
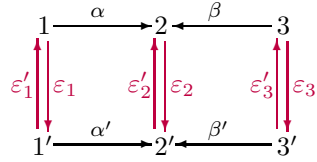
\vspace{-2cm}
			\begin{equation*}
				\begin{picture}(100,40)(0,20)
					\put(49,8){\small $1'$}
					\put(50,31){\small $1$}
					\put(72,8){\small $2'$}
					\put(72,31){\small $2$}
					\put(95,8){\small $3'$}
					\put(95,31){\small $3$}
					
					\put(53,10){\vector(1,0){18.5}}
					\put(53,32.5){\vector(1,0){18.5}}
					
					\put(94,10){\vector(-1,0){18}}
					\put(94,32.5){\vector(-1,0){18.5}}
					
					\put(60,12.5){$_{\alpha'}$}
					\put(60,35){$_\alpha$}
					\put(82,12.5){$_{\beta'}$}
					\put(82,35){$_\beta$}
					\color{purple}
					\put(50,13){\vector(0,1){17}}
					\put(52,29.5){\vector(0,-1){17}}
					\put(72,13){\vector(0,1){17}}
					\put(74,29.5){\vector(0,-1){17}}
					\put(95,13){\vector(0,1){17}}
					\put(97,29.5){\vector(0,-1){17}}
					
					\put(45,20){\small $\varepsilon_1'$}
					\put(53,20){\small $\varepsilon_1$}
					\put(67,20){\small $\varepsilon_2'$}
					\put(75,20){\small $\varepsilon_2$}
					\put(90,20){\small $\varepsilon_3'$}
					\put(98,20){\small $\varepsilon_3$}
				\end{picture}
			\end{equation*}
			\vspace{0.2cm}
		\end{center}
		and $I^{\sharp}$ is generated by all possible quadratic relations
		\begin{eqnarray*}
			&&\varepsilon_i'\varepsilon_i, \,\,\,\,\varepsilon_i\varepsilon_i', \quad \forall 1\leq i\leq 3\\
			&&\alpha' \varepsilon_1 -\varepsilon_2 \alpha,\,\,\,\, \alpha\varepsilon_1'- \varepsilon_2' \alpha',\,\,\,\,\beta \varepsilon_3' -\varepsilon_2' \beta',\,\,\,\, \beta' \varepsilon_3- \varepsilon_2 \beta.
		\end{eqnarray*}
		
		(c) Let $Q$ be the quiver such that $Q_0= \{1,2\}$, and $Q_1=\emptyset$. Then $\Lambda\cong R_2\times R_2$.
		
		(d) Let $\QJ$ be the Jordan quiver. Then its double framed quiver $\QJ^\sharp$ is 
		\begin{center}\setlength{\unitlength}{0.7mm}
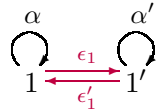
			\begin{picture}(50,13)(0,0)
				
				\put(0,9){\small $\alpha$}
				\put(20,9){\small $\alpha'$}
				\qbezier(-1,1)(-3,3)(-2,5.5)
				\qbezier(-2,5.5)(1,9)(4,5.5)
				\qbezier(4,5.5)(5,3)(3,1)
				\put(3.1,1.4){\vector(-1,-1){0.3}}
				\qbezier(19,1)(17,3)(18,5.5)
				\qbezier(18,5.5)(21,9)(24,5.5)
				\qbezier(24,5.5)(25,3)(23,1)
				\put(23.1,1.4){\vector(-1,-1){0.3}}
				
				\put(0,-4){$1$}
				\put(19,-4){$1'$}
				\color{purple}
				\put(4,0){\vector(1,0){14}}
				\put(18,-2){\vector(-1,0){14}}
				\put(10,0){$^{\epsilon_1}$}
				\put(10,-8){$^{\epsilon_1'}$}

			\end{picture}
			\vspace{0.3cm}
		\end{center}
		and $I^\sharp$ is generated by all possible quadratic relations
		$\alpha\epsilon_1'=\epsilon_1'\alpha'$, $\alpha'\epsilon_1=\epsilon_1\alpha$.

		(e) Let  $\QK: \xymatrix{0\ar@<0.5ex>[r]^\alpha \ar@<-0.5ex>[r]_\beta& 1}$ be the Kronecker quiver. Then its double framed quiver 
		$\QK^\sharp$ is 
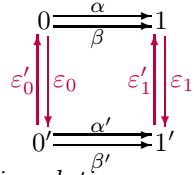
		\begin{center}\setlength{\unitlength}{0.7mm}
			\vspace{-2cm}
			\begin{equation*}
				\begin{picture}(100,40)(0,20)
					\put(49,8){\small $0'$}
					\put(50,31){\small $0$}
					\put(72,8){\small $1'$}
					\put(72,31){\small $1$}
					
					\put(53,11){\vector(1,0){18.5}}
					\put(53,9){\vector(1,0){18.5}}	\put(53,33.5){\vector(1,0){18.5}}
					\put(53,31.5){\vector(1,0){18.5}}
					\put(60,12.5){$_{\alpha'}$}
					\put(60,5.5){$_{\beta'}$}
					\put(60,35){$_\alpha$}
					\put(60,29){$_\beta$}
					\color{purple}
					\put(50,13){\vector(0,1){17}}
					\put(52,29.5){\vector(0,-1){17}}
					\put(72,13){\vector(0,1){17}}
					\put(74,29.5){\vector(0,-1){17}}
					
					\put(45,20){\small $\varepsilon_0'$}
					\put(53,20){\small $\varepsilon_0$}
					\put(67,20){\small $\varepsilon_1'$}
					\put(75,20){\small $\varepsilon_1$}
				\end{picture}
			\end{equation*}
			\vspace{0.2cm}
		\end{center}
		
		and $I^\sharp$ is generated by all possible quadratic relations
		\begin{eqnarray*}
			&&\varepsilon_0\varepsilon_0', \,\,\,\,\varepsilon_0'\varepsilon_0, \,\,\,\,\varepsilon_1'\varepsilon_1, \,\,\,\,\varepsilon_1\varepsilon_1',\,\,\,\, \alpha' \varepsilon_0 -\varepsilon_1 \alpha,\,\,\,\, \alpha \varepsilon_0'- \varepsilon_1' \alpha'\,\,\,\, \beta' \varepsilon_0 -\varepsilon_1 \beta,\,\,\,\, \beta \varepsilon_0'- \varepsilon_1' \beta'.
		\end{eqnarray*}
	\end{example}

	\begin{example}
		\label{example diagonal}
		Let $Q=(Q_0,Q_1)$ be a quiver, $Q^{\sharp}$ be its double framed quiver. Let $Q^{\dbl} =Q\sqcup  Q^{\diamond}$,  where $Q^{\diamond}$ is an identical copy of $Q$ with a vertex set $\{i^{\diamond} \mid i\in Q_0\}$ and an arrow set $\{ \alpha^{\diamond} \mid \alpha \in Q_1\}$.
		Let $\Lambda=\bfk Q\otimes_\bfk R_2$, $\Lambda^{\diamond} =k Q^{\diamond} \otimes_\bfk R_2$, and $\Lambda^{\dbl}:=\bfk Q^{\dbl}\otimes R_2$. Then the double framed quiver $(Q^{\dbl})^{\sharp}$ of $Q^{\dbl}$ is $Q^{\sharp}\sqcup (Q^{\sharp})^{\diamond}$, and
		$\Lambda^{\dbl}=\Lambda\times \Lambda^{\diamond} \cong \Lambda\times \Lambda$.
	\end{example}

	\subsection{$\imath$Quiver algebras}
	\label{subsec:i-QA}
	
	The involution $\btau$ of a quiver $Q$ induces an involution ${\btau}^{\sharp}$ of $Q^{\sharp}$ defined by
	\begin{itemize}
		\item ${\btau}^{\sharp}(i)=(\btau i)'$, ${\btau}^{\sharp}(i') =\btau i$ for any $i\in Q_0$;
		\item ${\btau}^{\sharp}(\varepsilon_i)= \varepsilon_{\btau i}'$, ${\btau}^{\sharp}(\varepsilon_i')= \varepsilon_{\btau i}$ for any $i\in Q_0$;
		\item ${\btau}^{\sharp}(\alpha)= (\btau\alpha)'$, ${\btau}^{\sharp}(\alpha')=\btau\alpha$ for any $\alpha\in Q_1$.
	\end{itemize}
	So starting  from the $\imath$quiver $(Q, \btau)$ we have constructed a new $\imath$quiver $(Q^{\sharp}, {\btau}^{\sharp})$.
	
	The action of ${\btau}^{\sharp}$ preserves $I^{\sharp}$. Hence ${\btau}^{\sharp}$ induces an involution, again denoted by ${\btau}^{\sharp}$, of the algebra $\Lambda$.

	\begin{definition}
		The fixed point subalgebra of $\Lambda$ under ${\btau}^{\sharp}$,
		\begin{equation}
			\label{eq:iLa}
			\Lambda^\imath
			= \{x\in \Lambda\mid {\btau}^{\sharp}(x) =x\},
		\end{equation}
		is called the {\rm $\imath$quiver algebra} of $(Q, \btau)$.
	\end{definition}

	According to \cite[Proposition 2.7]{LW22}, we have $\Lambda^\imath \cong \bfk \ov{Q} / \ov{I}$, where
	\begin{itemize}
		\item[(i)] $\ov{Q}$ is constructed from $Q$ by adding a loop $\varepsilon_i$ at the vertex $i\in Q_0$ if $\btau i=i$, and adding an arrow $\varepsilon_i: i\rightarrow \btau i$ for each $i\in Q_0$ if $\btau i\neq i$;
		\item[(ii)] $\ov{I}$ is generated by
		\begin{itemize}
			\item[(1)] (Nilpotent relations) $\varepsilon_{i}\varepsilon_{\btau i}$ for any $i\in\I$;
			\item[(2)] (Commutative relations) $\varepsilon_i\alpha-\btau(\alpha)\varepsilon_j$ for any arrow $\alpha:j\rightarrow i$ in $Q_1$.
		\end{itemize}
	\end{itemize}
	

	By \cite[Corollary 2.12]{LW22}, $\bfk Q$ is naturally a subalgebra and also a quotient algebra of $\Lambda^\imath$.
	Viewing $\bfk Q$ as a subalgebra of $\Lambda^{\imath}$, we have a restriction functor
	\[
	\res: \mod (\Lambda^{\imath})\longrightarrow \mod (\bfk Q).
	\]
	Viewing $\bfk Q$ as a quotient algebra of $\Lambda^{\imath}$, we obtain a pullback functor
	\begin{equation}\label{eqn:rigt adjoint}
		\iota:\mod(\bfk Q)\longrightarrow\mod(\Lambda^{\imath}).
	\end{equation}
	Hence a simple module $S_i (i\in Q_0)$ of $\bfk Q$ is naturally a simple $\iLa$-module.

	For each $i\in Q_0$, define a $\bfk$algebra (which can be viewed as a subalgebra of $\iLa$)
	\begin{align}\label{dfn:Hi}
		\BH _i:=\left\{ \begin{array}{cc}  \bfk[{\textcolor{purple}{\varepsilon_i}}]/({\textcolor{purple}{\varepsilon_i}}^2) & \text{ if }i=\btau i,
			\\
			\bfk(\xymatrix{i \ar@<0.5ex>@[purple][r]^{\textcolor{purple}{\varepsilon_i}} & \btau i \ar@<0.5ex>@[purple][l]^{\textcolor{purple}{\varepsilon_{\btau i}}}})/( \varepsilon_i\varepsilon_{\btau i},\varepsilon_{\btau i}\varepsilon_i)  &\text{ if } \btau i \neq i .\end{array}\right.
	\end{align}
	Note that $\BH _i=\BH _{\btau i}$ for any $i\in Q_0$.
	Choose one representative for each $\btau$-orbit on $\I$, and let
	\begin{align}
		\label{eq:ci}
		\ci = \{ \text{the chosen representatives of $\btau$-orbits in $\I$} \}.
	\end{align}
	
	Define the following subalgebra of $\Lambda^{\imath}$:
	\begin{equation}  \label{eq:H}
		\BH =\bigoplus_{i\in \ci }\BH _i.
	\end{equation}
	Note that $\BH $ is a radical square zero selfinjective algebra. Denote by
	\begin{align}
		\res_\BH :\mod(\iLa)\longrightarrow \mod(\BH )
	\end{align}
	the natural restriction functor.
	On the other hand, as $\BH $ is a quotient algebra of $\iLa$, every $\BH $-module can be viewed as a $\iLa$-module.
	
	Recall the algebra $\BH _i$ for $i \in \ci$ from \eqref{dfn:Hi}. For $i\in Q_0 =\I$, define the indecomposable module over $\BH _i$ (if $i\in \ci$) or over $\BH_{\btau i}$ (if $i\not \in \ci$)
	\begin{align}
		\label{eq:E}
		\E_i =\begin{cases}
			\bfk[\varepsilon_i]/(\varepsilon_i^2), & \text{ if }i=\btau i;
			\\
			\xymatrix{\bfk\ar@<0.5ex>[r]^1 & \bfk\ar@<0.5ex>[l]^0} \text{ on the quiver } \xymatrix{i\ar@<0.5ex>[r]^{\varepsilon_i} & \btau i\ar@<0.5ex>[l]^{\varepsilon_{\btau i}} }, & \text{ if } i\neq \btau i.
		\end{cases}
	\end{align}
	Then $\E_i$, for $i\in Q_0$, can be viewed as a $\iLa$-module and will be called a {\em generalized simple} $\iLa$-module.
	
	\begin{example}\label{example 2}
		(a) We continue Example \ref{exDoubleQuiver}(a), with $\btau=\Id$. Then $\Lambda^\imath$ is isomorphic to the algebra with its quiver $\ov Q$ and relations as follows:
		\begin{center}\setlength{\unitlength}{0.7mm}
			\begin{picture}(50,13)(0,0)
				\put(0,-2){$1$}
				\put(4,0){\vector(1,0){14}}
				\put(10,0){$^{\alpha}$}
				\put(20,-2){$2$}
				\color{purple}
				\put(0,9){\small $\varepsilon_1$}
				\put(20,9){\small $\varepsilon_2$}
				\qbezier(-1,1)(-3,3)(-2,5.5)
				\qbezier(-2,5.5)(1,9)(4,5.5)
				\qbezier(4,5.5)(5,3)(3,1)
				\put(3.1,1.4){\vector(-1,-1){0.3}}
				\qbezier(19,1)(17,3)(18,5.5)
				\qbezier(18,5.5)(21,9)(24,5.5)
				\qbezier(24,5.5)(25,3)(23,1)
				\put(23.1,1.4){\vector(-1,-1){0.3}}
			\end{picture}
			\vspace{0.2cm}
		\end{center}
		\[
		\varepsilon_1^2=0=\varepsilon_2^2, \quad \varepsilon_2 \alpha=\alpha\varepsilon_1.
		\]
		Indeed, $\Lambda^\imath$ is isomorphic to $\bfk Q\otimes_\bfk R_1$. This class of algebras was independently studied by Ringel-Zhang and Geiss-Leclerc-Schr\"{o}er \cite{GLS} from very different perspectives. 
		
		(b) We continue Example \ref{exDoubleQuiver}(b), with the involution $\tau$ of $Q$ swapping vertices $1$ and $3$.
		Then $\Lambda^\imath$ is isomorphic to the algebra with its quiver $\ov Q$ and relations as follows:
		\begin{center}\setlength{\unitlength}{0.7mm}
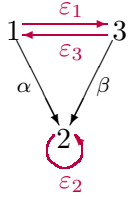
			\begin{picture}(50,20)(0,-10)
				\put(0,-2){$1$}
				\put(20,-2){$3$}
				\put(2,-11){$_{\alpha}$}
				\put(17,-11){$_{\beta}$}
				\put(2,-2){\vector(1,-2){8}}
				\put(20,-2){\vector(-1,-2){8}}
				\put(9.5,-22.5){$2$}
				\color{purple}
				\put(3,1){\vector(1,0){16}}
				\put(19,-1){\vector(-1,0){16}}
				\put(10,3){\small ${\varepsilon_1}$}
				\put(10,-5){\small ${\varepsilon_3}$}
				\put(10,-30){\small ${\varepsilon_2}$}
				\begin{picture}(50,23)(-10,19)
					\color{purple}
					\qbezier(-1,-1)(-3,-3)(-2,-5.5)
					\qbezier(-2,-5.5)(1,-9)(4,-5.5)
					\qbezier(4,-5.5)(5,-3)(3,-1)
					\put(3.1,-1.4){\vector(-1,1){0.3}}
				\end{picture}
			\end{picture}
			\vspace{1.4cm}
		\end{center}
		\[
		\varepsilon_1\varepsilon_3=0=\varepsilon_3\varepsilon_1,
		\quad
		\varepsilon_2^2=0,
		\quad
		\varepsilon_2 \beta=\alpha\varepsilon_3,
		\quad
		\varepsilon_2 \alpha=\beta\varepsilon_1.
		\]

		(c) We continue Example \ref{exDoubleQuiver}(c) with $\btau =\Id$. Then, $\Lambda^\imath\cong R_1\times R_1$.
		
		(d) We continue Example \ref{exDoubleQuiver}(d) of Jordan quiver with $\btau =\Id$.
		Then $\Lambda^\imath$ is isomorphic to the algebra with its quiver $\ov{Q}$ and relations as follows:
		\vspace{-4mm}
		\begin{equation}
			\label{fig:Jordan}
			\begin{picture}(100,30)(0,10)
				%
				\begin{tikzpicture}[scale=1.1]
					\draw[color=purple] (7,0) .. controls (7.2,0.35) and (7.6,0.35) .. (7.6,-0.05);
					\draw[-latex,color=purple] (7.6,-0.05) .. controls (7.6,-0.45) and (7.2,-0.45) .. (7,-0.1);
					\draw[-] (6.7,0) .. controls (6.5,0.35) and (6.1,0.35) .. (6.1,-0.05);
					\draw[-latex] (6.1,-0.05) .. controls (6.1,-0.45) and (6.5,-0.45) .. (6.7,-0.1);
					\node at (6.85,0 ){ \small$1$};
					\node at (5.9,0 ){ \small$\alpha$};
					\node at (7.8,0 ){ \small\color{purple}{$\varepsilon$}};
				\end{tikzpicture}
			\end{picture}
		\end{equation}
		\vspace{-2mm}
		
		\[\alpha\epsilon=\epsilon\alpha,\quad \epsilon^2=0.\]
		
		(e) We continue Example \ref{exDoubleQuiver}(e) of Kronecker quiver with $\btau =\Id$.
		Then $\LaK^\imath$ is isomorphic to the algebra with its quiver $\ov \QK$ and relations as follows:
		\begin{center}\setlength{\unitlength}{0.7mm}
			\vspace{-0.2cm}
			\begin{equation}
				\label{eq:iquiver}
				\begin{picture}(30,13)(0,0)
					\put(0,-2){\small $0$}
					\put(2.5,0.5){\vector(1,0){17}}
					\put(2.5,-2){\vector(1,0){17}}
					\put(10,-0.5){$^{\alpha}$}
					\put(10,-8.5){$^{\beta}$}
					\put(20,-2){\small $1$}
					\color{purple}
					\put(0,9){\small $\varepsilon_0$}
					\put(20,9){\small $\varepsilon_1$}
					
					\qbezier(-1,1)(-3,3)(-2,5.5)
					\qbezier(-2,5.5)(1,9)(4,5.5)
					\qbezier(4,5.5)(5,3)(3,1)
					\put(3.1,1.4){\vector(-1,-1){0.3}}
					
					\qbezier(19,1)(17,3)(18,5.5)
					\qbezier(18,5.5)(21,9)(24,5.5)
					\qbezier(24,5.5)(25,3)(23,1)
					\put(23.1,1.4){\vector(-1,-1){0.3}}
				\end{picture}
			\end{equation}
			\vspace{-0.2cm}
		\end{center}
		\[
		\varepsilon_0^2=0=\varepsilon_1^2, \quad \varepsilon_1 \alpha=\alpha\varepsilon_0, \varepsilon_1 \beta=\beta\varepsilon_0.
		\]
	\end{example}
	
	\begin{example}
		\label{example diagonal 2}
		{\rm ($\imath$quiver of diagonal type)}
		Continuing Example \ref{example diagonal}, we let $\rm{swap}$ be the involution of $Q^{\rm dbl}$ uniquely determined by $\swa(i)=i^\diamond$ for any $i\in Q_0$. 
		Then $(\Lambda^{\dbl})^\imath$ is isomorphic to $\Lambda$. Explicitly, let $(\ov{Q}^{\dbl},\ov{I}^{\dbl})$ be the bound quiver of $(\Lambda^{\dbl})^\imath$. Then $(\ov{Q}^{\dbl},\ov{I}^{\dbl})$ coincides with the double $\imath$quiver $(Q^{\sharp},I^{\sharp})$. So we just use $(Q^{\sharp},I^{\sharp})$ as the bound quiver of $(\Lambda^{\dbl})^\imath$ and identify $(\Lambda^{\dbl})^\imath$ with $\Lambda$.
	\end{example}

	\subsection{Gorenstein properties of $\imath$quiver algebras}
	\label{subsec:iQAGorenstein}
	
	Let $d$ be a positive integer. Recall \cite{Ha3,EJ} a noetherian algebra $A$ is  {\em $d$-Gorenstein} if $\ind{}_AA\leq d$ and $\ind A_A\leq d$.
	
	\begin{proposition}[\text{\cite[Proposition 2.2]{LW20}}]
		\label{proposition of 1-Gorenstein}
		For a general $\imath$quiver $(Q,\btau)$, $\Lambda$ and $\Lambda^\imath$ are $1$-Gorenstein algebras.
	\end{proposition}

	\begin{lemma}[\text{\cite[Lemma 2.4]{LW20}}]
		\label{lem: resolution}
		For any $M\in\mod(\Lambda^\imath)$, there exist short exact sequences
		\begin{align}
			\label{resolution 1}
			0\longrightarrow M\longrightarrow H^M\longrightarrow X^M\longrightarrow0\\
			\label{resolution 2}
			0\longrightarrow X_M\longrightarrow H_M\longrightarrow M\longrightarrow 0,
		\end{align}
		with $H^M,H_M\in \cp^{\leq1}(\Lambda^\imath)$.
	\end{lemma}

	By definition of $(\ov{Q},\ov{I})$, $\Lambda^\imath$ is positively graded by $\deg\varepsilon_i=1$, $\deg\alpha=0$ for any $i\in Q_0$, $\alpha\in Q_1$. Note that $\Lambda^\imath_0=\bfk Q$. Let $\mod^{\Z}(\Lambda^\imath)$ be the finite-dimensional nilpotent $\Z$-graded $\Lambda^\imath$-modules. The singularity category $\cd_{sg}(\mod^{\Z}(\Lambda^\imath))$ is defined to be the Verdier quotient of $\cd^b(\mod^{\Z}(\Lambda^\imath))$ modulo the thick subcategory generated by modules with finite projective dimensions. 
	
	Let $G$ be the composition
	\[
	\cd^b(\mod(\bfk Q))\xrightarrow{\cd^b(\iota)} \cd^b(\mod^\Z (\Lambda^{\imath}))\stackrel{\pi}{\longrightarrow}  \cd_{sg}(\mod^\Z (\Lambda^{\imath})),
	\]
	where $\cd^b(\iota)$ is the derived functor of $\iota$ since $\iota$ is exact. 
	
	Let $\widehat{\btau}$ be the triangulated auto-equivalence of $D^b(\mod(\bfk Q))$ induced by $\btau$. By \cite[Theorem 3.18]{LW22}, we have
	\begin{align}
		\label{eqn: shift}
		(1)\circ G\simeq G\circ \Sigma\circ\widehat{\btau},
	\end{align}
	where $(1)$ is the degree shift functor, and $\Sigma$ is the suspension functor of $\cd^b(\mod(\bfk Q))$. 
	
	By \cite[Lemma 2.10]{LW20}, the restriction of $G$ provides us a triangulated equivalence
	\[\cd^b(\mod(\bfk Q))\stackrel{\simeq}{\longrightarrow} \cd_{sg}(\fdmz(\Lambda^\imath)).\]
	Combing with \eqref{eqn: shift}, we have the following theorem.

	\begin{theorem}[\text{\cite[Theorem 2.11]{LW20}}]
		\label{thm:sigma}
		Let $(Q, \btau)$ be an $\imath$quiver. Then $D^b(\mod(\bfk Q))/\Sigma \circ \widehat{\btau}$ is a triangulated orbit category \`a la Keller \cite{Ke05}, and we have the following triangulated equivalence
		\[
		\cd_{sg}(\mod(\Lambda^{\imath}))\simeq \cd^b(\mod(\bfk Q))/\Sigma \circ \widehat{\btau}.
		\]
	\end{theorem}


	As a corollary,	for any $M\in D_{sg}(\mod(\Lambda^{\imath}))$, there exists a unique (up to isomorphisms) module $N\in \mod(\bfk Q)\subseteq \mod(\Lambda^{\imath})$ such that
	$M\cong N$ in $D_{sg}(\mod(\Lambda^{\imath}))$; see  \cite[Corollary 3.21]{LW22}.

	\begin{corollary}[\text{\cite[Corollary 2.13]{LW20}}]
		\label{cor: res proj}
		For any $M\in\mod(\Lambda^{\imath})$ the following are equivalent.
		\begin{itemize}
			\item[(i)] $\pd M<\infty$;
			\item[(ii)] $\ind M<\infty$;
			\item[(iii)] $\pd M\leq1$;
			\item[(iv)] $\ind M\leq1$;
			\item[(v)] $\res_\BH (M)$ is projective as an $\BH $-module.
		\end{itemize}
	\end{corollary}
	
	Summarizing, we have established the following.
	\begin{proposition}
		The		$\mod(\Lambda^\imath)$ is an abelian category satisfying (Ea)-(Ed). 
	\end{proposition}

	
	\section{$\imath$Hall algebra for the Jordan quiver} 
	\label{sec:iJordan}
	
	In this short section, we formulate the connection between the $\imath$Hall algebra for the Jordan quiver and ring of symmetric functions. 
	
	Let $ \QJ$ be the Jordan quiver, i.e., a quiver with a single vertex $1$ and a single loop $\alpha:1\rightarrow 1$.
	It is well known that $\mod(\bfk\QJ)$ is a uniserial category. Let $S$ be the simple object in $\mod(\bfk\QJ)$.
	Then any indecomposable object of $\mod(\bfk\QJ)$  (up to isomorphisms) is of the form $S^{(n)}$ of length $n\geq1$.
	Thus the set of isomorphism classes of $\mod(\bfk\QJ)$ is canonically isomorphic to the set $\cp$ of all partitions, via the assignment
	\begin{align}
		\la =(\la_1,\la_2,\cdots,\la_r)\mapsto S^{(\la)}= S^{(\la_1)}\oplus \cdots \oplus S^{(\la_r)}.
	\end{align}
	
	Let $\Lambda_{\texttt{J}}^\imath$ be the $\imath$quiver algebra of the $\imath$quiver $\QJ$ equipped with trivial involution \cite{LW22, LW20}, see Example~\ref{example 2}(d) and the $\imath$quiver in \eqref{fig:Jordan}. 
	In particular, $\Lambda_{\texttt{J}}^\imath$ is a commutative $\bfk$-algebra. 
	
	The quotient algebra
	\begin{align}
		\tMHLJ:=\ch\big(\Lambda_{\texttt{J}}^\imath\big)/I
	\end{align}
	is called the {\em$\imath$Hall algebra} of $\Lambda_{\texttt{J}}^\imath$, where $I$ is the two-sided ideal defined as in \eqref{def:I}; comparing semi-derived Hall algebras deinfed in Section~\ref{sec:Semi-derivedHall}, we do not apply any localization here.
	
	\begin{lemma}[\text{\cite[Lemma 4.5]{LRW20}, \cite[Proposition 4.4]{LRW21}}]
		\label{lem:comm Jordan}
		The algebra $\tMHLJ$ is commutative, which is isomorphic to
		\begin{enumerate}
			\item
			the polynomial $\Q$-algebra in the infinitely many generators
			$[\E_1], [S^{(1^r)}],$ for $r\ge 1$;
			\item
			the polynomial $\Q$-algebra in the infinitely many generators
			$[\E_1], [S^{(r)}],$ for $r\ge 1$.
		\end{enumerate}
	\end{lemma}

	Similar to the Hall-Littlewood functions (see \S\ref{subsec:SHconstruction}), we define the $\imath$Hall-Littlewood functions as follows.

	Let $1\leq i<j$. 
	Similar to the raising operator $R_{ij}$, we define the lowering operator $L_{ij}$ acting on integer sequences by
	\[
	L_{ij} \alpha =(\ldots, \alpha_i-1, \ldots, \alpha_j-1, \ldots).
	\]
	Recall that $q_r=Q_r$ for any $r\geq1$; see \S\ref{subsec:SHconstruction}. 
	The actions of  lowering operators on $q_\alpha$ are given by letting
	\[
	L_{ij} q_\alpha =q_{L_{ij}\alpha}.
	\]
	Note that all raising and lowering operators commute with each other.
	
	Let $t, \vth$ be indeterminates.
	We consider the ring $\La_{t,\vth} = \Q(t)[\vth] [q_1, q_2, \ldots]$.
	
	\begin{definition}[$\imath$HL functions; \cite{LRW21}]
		\label{def:V}
		For any integer vector $\alpha$, we define $Q^\imath_\alpha \in \La_{t,\vth}$ by:
		\begin{align}  \label{eq:LRq}
			Q^\imath_\alpha = \prod_{1\le i<j} \frac{1 - \vth L_{ij}}{1 - \vth t L_{ij}} \frac{1 -R_{ij}}{1 -tR_{ij}} q_\alpha.
		\end{align}
		Note that $Q^\imath_{(r)} =Q_{r} =q_r$.
	\end{definition}

	\begin{lemma} [\cite{LRW21}]
		\label{def:iQ}
		For any composition $\la$, $\iQ_\la$ is the coefficient of $u^\la=u_1^{\la_1}u_2^{\la_2}\cdots$ in
		\begin{align}
			\label{HL functions}
			\iQ(u_1,u_2,\dots)&= \prod_{i\geq1} Q(u_i)\prod_{i<j} F(u_i^{-1}u_j)F(\vth  u_iu_j)
			\\\notag
			&= \prod_{i\geq1} Q(u_i) \prod_{i<j} \frac{(1-u_i^{-1}u_j)(1- \vth  u_iu_j)}{(1- tu_i^{-1}u_j)(1-t \vth  u_iu_j)}.
		\end{align}
	\end{lemma}
	
	Denote $\La_{\vth} :=\La_{t,\vth} \otimes_{\Q(t)} \Q$, where $t$ acts on $\Q$ by $q^{-1}$.
	\begin{theorem}[\text{\cite[Theorem 4.9]{LRW21}}]
		\label{thm:iso}
		There exists a $\Q$-algebra isomorphism $\Phi_{(q)}:   \iH(\bfk \QJ)\rightarrow \La_{\vth}$ such that
		\begin{align*}
			\Phi_{(q)}(\E_1)=q\vth,\qquad
			\Phi_{(q)}([S^{(r)}])=q^{r}\iQ_{(r)} =q^r Q_{r} \quad (r\ge 1).
		\end{align*}
		Moreover, for any partition $\la$, we have
		\begin{align}
			\Phi_{(q)}([S^{(\lambda)}])=q^{|\lambda|+n(\lambda)}\iQ_\la .
		\end{align}
	\end{theorem}

	\section{$\imath$Hall algebras for Dynkin $\imath$quivers}
	\label{sec:iHallDynkin}
	
	In this section, we study the $\imath$Hall algebras for $\imath$quivers without loops (i.e., the twisted semi-derived Hall algebras of $\imath$quiver algebras), and establish their Hall bases. We show that $\imath$Hall algebra for a Dynkin $\imath$quiver is isomorphic to the corresponding universal $\imath$quntum group.

	\subsection{$\imath$Hall algebras and $\imath$Hall bases}
	
	Recall we have defined an $\imath$quiver algebra $\Lambda^\imath$ \eqref{eq:iLa} associated to an $\imath$quiver $(Q,\btau)$; we shall assume that $Q$ contains no loops in this section. According  to \S\ref{subsec:iQAGorenstein}, $\mod(\Lambda^\imath)$ is a  $1$-Gorenstein algebra satisfying (Ea)-(Ed). Thus, we can define the semi-derived Ringel-Hall algebra $\utMH$ of $\iLa$. 
	
	For $K,M\in \mod(\iLa)$ with $K\in\cp^{\leq 1}(\Lambda^\imath)=\cp^{<\infty}(\Lambda^\imath)$, by Corollary \ref{cor: res proj}, the following pairings are well defined:
	\begin{align}
		\label{left Euler form}
		\langle K,M\rangle &=\sum_{i=0}^{\infty}(-1)^i \dim_\bfk\Ext^i(K,M)=\dim_\bfk\Hom(K,M)-\dim_\bfk\Ext^1(K,M),
		\\
		\label{right Euler form}
		\langle M,K\rangle &=\sum_{i=0}^{\infty}(-1)^i \dim_\bfk\Ext^i(M,K)=\dim_\bfk\Hom(M,K)-\dim_\bfk\Ext^1(M,K).
	\end{align}
	These pairings induce bilinear pairings (called Euler forms) between the Grothendieck groups $K_0(\cp^{\leq 1}(\Lambda^\imath))$ and $K_0(\mod(\iLa))$:
	\begin{align}
		\label{eq:Euler1}
		\langle\cdot,\cdot\rangle: K_0(\cp^{\leq 1}(\Lambda^\imath))\times K_0(\mod(\Lambda^\imath))\longrightarrow \Z,
		\\
		\label{eq:Euler2}
		\langle\cdot,\cdot\rangle: K_0(\mod(\Lambda^\imath))\times K_0(\cp^{\leq 1}(\Lambda^\imath))\longrightarrow \Z.
	\end{align}
	
	Denote by $\langle\cdot,\cdot\rangle_Q$ the Euler form of $\bfk Q$, and 
	$(\cdot,\cdot)_Q$ the symmetrized Euler form. Denote by $S_i$ the simple $\bfk Q$-module (respectively, $\Lambda^{\imath}$-module) corresponding to the vertex $i\in Q_0$ (respectively, $i\in\ov{Q}_0$).
	
	\begin{lemma}
		\label{lem:Euler}
		For $K,K'\in \cp^{\leq1}(\Lambda^\imath)$, $M\in\mod(\Lambda^\imath)$, $i,j\in\I$, we have
		\begin{align}
			\label{Eform1}
			\langle K,M\rangle =\langle \res_{\BH}(K),M\rangle,& \qquad \langle M,K\rangle =\langle M,\res_{\BH}(K)\rangle,
			\\
			\label{Eform2}
			\langle \bK_i,S_j\rangle=\langle S_i,S_j\rangle_Q,&\qquad \langle S_j,\bK_i\rangle =\langle S_j, S_{\btau i} \rangle_Q,
			\\
			\label{Eform3}
			\langle K,K'\rangle&= \frac{1}{2}\langle \res(K),\res(K')\rangle_Q.
		\end{align}
	\end{lemma}
	
	Via the restriction functor $\res: \mod(\Lambda^{\imath})\rightarrow\mod (\bfk Q)$, we define the twisted semi-derived Ringel-Hall algebra to be the $\Q(\sqq)$-algebra on the same vector space as $\utMH$ with twisted multiplication given by
	\begin{align}
		\label{eqn:twsited multiplication}
		[M]* [N] =\sqq^{\langle \res(M),\res(N)\rangle_Q} [M]\diamond[N].
	\end{align}
	We shall denote this algebra $(\utMH, *)$ by $\tMH$ or preferably,
	$\tMHk$, and call it the {\em Hall algebra associated to the $\imath$quiver $(Q, \btau)$}, (or an {\em $\imath$Hall algebra}, for short). 
	The {\em twisted quantum torus} $\tTL$ is defined to be the subalgebra of $\tMHk$ generated by $[K]$, $K\in\cp^{\leq1}(\Lambda^\imath)$. 	For any $\alpha=\sum\limits_{i\in\I} a_i\widehat{S_i}\in K_0(\mod(\bfk Q))=\Z^\I$, define
	\[
	\bK_\alpha:= [X]* [Y]^{-1} \in \tMHk,
	\]
	where $X= \bigoplus\limits_{i\in\I:a_i\geq0} \bK_i^{\oplus a_i}$ and $Y=\bigoplus\limits_{i\in\I:a_i<0} \bK_i^{\oplus (-a_i)}$.
	Then  $\tTL$ is a Laurent polynomial algebra generated by $[\bK_i]$, for $i\in \I$; and $[\bK_\alpha]*[\bK_\beta]=[\bK_{\alpha+\beta}]$ for any $\alpha,\beta\in\Z^\I$.

	By definition and Lemma \ref{lem:Euler}, we have 	
	\begin{align}
		[M]* [\E_i]&= \sqq^{\langle \res M,S_i\rangle_Q-\langle \res M,S_{\tau i}\rangle_Q} [M\oplus \E_i]
		\label{right module structure},
		\\
		\quad
		[\E_i]*[M]&= \sqq^{\langle S_{\tau i},\res M\rangle_Q-\langle S_{i} ,\res M\rangle_Q} [\E_i\oplus M],
	\end{align}
	and then 
	\begin{align}
		\label{eq:commKM}
		[\E_i]*[M]=\sqq^{(\res M,S_{\tau i})_Q-(\res M, S_i)_Q} [M]*[\E_i].
	\end{align}
	for any $i\in\I$, $M\in\mod(\Lambda^\imath)$.

	\begin{theorem}[$\imath$Hall basis; see \text{\cite[Theorem 3.6]{LW20}}]
		\label{thm:utMHbasis}
		Let $(Q,\btau)$ be an $\imath$quiver. Then
		\begin{equation}
			\label{eq:Hall basis}
			\big\{ [X]* \bK_\alpha ~\big |~ [X]\in\Iso(\mod(\bfk Q))\subseteq \Iso(\mod(\Lambda^{\imath})), \alpha\in \Z^\I \big\}
		\end{equation}
		is a basis of $\tMHk$.
	\end{theorem}

	There is a filtered algebra structure on $\tMHk$ by setting $\deg [M]=\widehat{M}\in K_0(\mod(\bfk Q))=\Z^{\I}$, and $\deg \K_\alpha=0$ for any $M\in\mod(\bfk Q)$ and $\alpha\in\Z^{\I}$. We denote the associated graded algebra
	\[
	\widetilde{\ch}(\bfk Q, \btau)^{\gr} = \bigoplus_{\alpha \in K_0(\mod(\bfk Q))} \widetilde{\ch}(\bfk Q, \btau)_{\alpha}^{\gr}.
	\]
	It is natural to view the quantum torus $\widetilde{\ct}(\Lambda^\imath)$ as a subalgebra of $\widetilde{\ch}(\bfk Q, \btau)^{\gr}$.
	By \cite[Lemma~ 5.4(ii)]{LW22}, the linear map
	\begin{align}
		\label{def:morph}
		\varphi: \widetilde{\ch}(\bfk Q)\longrightarrow \widetilde{\ch}(\bfk Q, \btau)^{\gr},\quad
		\varphi([M])=[M],\;  \forall M\in \mod(\bfk Q),
	\end{align} is an algebra embedding. Then we have 
	$\tMHk\cong \widetilde{\ch}(\bfk Q)*\widetilde{\ct}(\Lambda^\imath)$ as linear spaces.

	Recall that $\langle\cdot,\cdot\rangle_Q$ is the Euler form of $Q$. Define
	\begin{align*}
		(x,y)&= \langle x,y\rangle_Q+\langle y,x\rangle_Q.
	\end{align*}
	In particular, $(S_i,S_j)=c_{ij}$ for any $i,j\in \I$, the entries for the Cartan matrix $C$.

	\subsection{$\imath$Hall algebras and $\imath$quantum groups of Dynkin type}
	
	Let $(Q, \btau)$ be a Dynkin $\imath$quiver. In this case, the presentation of the universal $\imath$quantum group $\tUi$ in Theorem \ref{thm:Serre} can be simplified. That is, $\tUi$ is a $\Q(v)$-algebra with generators $B_i, \tk_i$ $(i\in \I)$, where $\tk_i$ are invertible, subject to the following relations \eqref{Dynrelation1}--\eqref{Dynrelation2}: for $\ell \in \I$, and $i\neq j \in \I$,
	\begin{align}
		\tk_i \tk_\ell =\tk_\ell \tk_i,
		\quad
		\tk_\ell B_i & = v^{c_{\btau \ell,i} -c_{\ell i}} B_i \tk_\ell,
		\label{Dynrelation1}
		\\
		B_iB_{j}-B_jB_i &=0, \quad \text{ if } c_{ij} =0 \text{ and }\btau i\neq j,
		\label{DynrelationBB}
		\\
		\label{DynrelationSerre}
		\sum_{s=0}^{1-c_{ij}} (-1)^s \qbinom{1-c_{ij}}{s}& B_i^{s}B_jB_i^{1-c_{ij}-s} =0, \quad \text{ if } j \neq \btau i\neq i,
		\\
		B_{\btau i}B_i -B_i B_{\btau i}& =   \frac{\tk_i -\tk_{\btau i}}{v-v^{-1}},
		\quad \text{ if } \btau i \neq i \text{ (and hence } c_{i,\tau i}=0),
		\label{Dynrelation5}
		\\
		B_i^2B_{j} - [2] B_iB_{j}B_i +&B_{j}B_i^2 = v \tk_i B_{j},
		%
		\quad \text{ if }  c_{ij} = -1 \text{ and }\btau i=i.
		\label{Dynrelation2}
	\end{align}
	
	Denote $\tUi_{\sqq} =\Q(\sqq) \otimes_{\Q(v)} \tUi$. 
	The following theorem provides an $\imath$Hall algebra realization of $\tUi$ (which also shows why $\tUi$ is more natural than $\Ui$ from a  categorical viewpoint). 
	
	\begin{theorem}
		\label{thm:main}
		Let $(Q, \btau)$ be a Dynkin $\imath$quiver. Then there exists a $\Q(\sqq)$-algebra isomorphism
		\begin{align*}
			\widetilde{\psi}: \tUi_{\sqq} &\longrightarrow \tMHk,
		\end{align*}
		which sends
		\begin{align}
			B_j \mapsto \frac{-1}{q-1}[S_{j}],\text{ if } j\in\ci,
			&\qquad\qquad
			\tk_i \mapsto - q^{-1}[\bK_i], \text{ if }\btau i=i \in \I;
			\label{eq:split}
			\\
			B_{j} \mapsto \frac{{\sqq}}{q-1}[S_{j}],\text{ if }j\notin \ci,
			&\qquad\qquad
			\tk_i \mapsto [\bK_i],\quad \text{ if }\btau i\neq i \in \I.
			\label{eq:extra}
		\end{align}
	\end{theorem}

	\begin{proof}
		It follows from \eqref{eq:commKM} that $\widetilde{\psi}$ preserves the relation \eqref{Dynrelation1}.
		
		For \eqref{DynrelationBB}, we have $\Hom_{\Lambda^\imath}(S_i,S_j)=0=\Ext^1_{\Lambda^\imath}(S_i,S_j)$ for any $i,j$ such that $c_{ij}=0$ and $\tau i\neq j$. So $[S_i]*[S_j]=[S_i\oplus S_j]=[S_j]*[S_i]$. Then $\widetilde{\psi}$ preserves the relation \eqref{DynrelationBB}.

		For \eqref{Dynrelation5}, we have $\Hom_{\Lambda^\imath}(S_i,S_{\tau i})=0$ and $\Ext^1_{\Lambda^\imath}(S_i,S_{\tau i})=\bfk$, and $K\cong \E_i$ for any non-split exact sequence
		$0\rightarrow S_{\tau i}\rightarrow K\rightarrow S_i\rightarrow$. So
		$[S_{\tau i}]*[S_i]-[S_i]*[S_{\tau i}]=(q-1) ([\E_{\tau i}]-[\E_i])$.

		The verification that  $\widetilde{\psi}$ preserves the relations \eqref{DynrelationSerre} and \eqref{Dynrelation2} is reduced to the rank $2$ Dynkin $\imath$quivers, and the details will be given in \S\ref{subsec:computationrank2I}--\S\ref{subsec:computationrank2II} below.
		
		It remains to prove that $\widetilde{\psi}: \tUi_{\sqq} \longrightarrow \tMHk$ is an isomorphism. The homomorphism $\widetilde{\psi}: \tUi_{\sqq} \rightarrow \tMHk$ restricts to an algebra homomorphism
		\begin{align*}
			\widetilde{\psi}: \widetilde{\bU}^{\imath 0}_{\sqq} &\longrightarrow  \tTL,
			\\
			\tk_i\mapsto - q^{-1} [\bK_i], \text{ if }\btau i=i,
			&\qquad
			\tk_i\mapsto [\bK_i], \text{ if }\btau i\neq i.
		\end{align*}
		Since both $\widetilde{\bU}^{\imath 0}_{\sqq}$ and $\tTL$ are Laurent polynomial algebras in the same number of generators, $\widetilde{\psi}:\widetilde{\bU}^{\imath 0}_{\sqq} \rightarrow \tTL$ is an isomorphism.

		We observe that $\widetilde{\psi}$ is a morphism of filtered algebras. Let $\widetilde{\psi}^{\gr}:  \tU^{\imath,\gr}_{\sqq} \longrightarrow \widetilde{\ch}(\bfk Q, \btau)^{\gr}$ be its associated graded morphism, and we obtain the following commutative diagram
		\begin{equation*}
			\xymatrix{ \U^-_{\sqq} \ar[r]^\phi \ar[d]^R &   \tU^{\imath,\gr}_{\sqq} \ar[d]^{ \widetilde{\psi}^{\gr} }  \\
				\widetilde{\ch}(\bfk Q) \ar[r]^\varphi & \widetilde{\ch}(\bfk Q, \btau)^{\gr} }
		\end{equation*}
		It follows that $\widetilde{\psi}^{\gr}\circ \phi$ is injective since $\varphi$ and $R: \U^-_{\sqq}\rightarrow \widetilde{\ch}(\bfk Q)$ are injective by Theorem~ \ref{thm:Ringel-Green} and \eqref{def:morph}.
		
		We claim that $\widetilde{\psi}^{\gr}$ is injective. Indeed, any element in $\tU^{\imath,\gr}$ is of form $\sum_{\alpha\in\Z^\I}\phi(V_\alpha)\cdot \tk_\alpha$, for $V_\alpha\in \U^-$. Here $\tk_\alpha=\prod_{i\in\I} \tk_i^{a_i}$ for $\alpha=\sum_{i\in\I} a_i \alpha_i$. Assume $\widetilde{\psi}^{\gr}(\sum_{\alpha}\phi(V_\alpha)\cdot \tk_\alpha)=0$, i.e.,  $\sum_{\alpha\in\Z^\I}\widetilde{\psi}^{\gr}(\phi(V_\alpha))* \bK_\alpha=0$.  By Theorem \ref{thm:utMHbasis}, we obtain $\widetilde{\psi}^{\gr}(\phi(V_\alpha))* \bK_\alpha=0$ for any $\alpha$, and then $\widetilde{\psi}^{\gr}(\phi(V_\alpha))=0$ since $\K_\alpha$ is invertible. So $V_\alpha=0$ by using the above commutativev diagram. It follows that $\widetilde{\psi}^{\gr}$ is injective.
		
		Now by a standard filtered algebra argument, we obtain that
		$\widetilde{\psi}: \tUi_{\sqq} \longrightarrow \tMHk$ is an algebra monomorphism. 
		
		Since $Q$ is of Dynkin type, we have $\widetilde{\ch}(\bfk Q)$ is generated by $[S_i]$ ($i\in\I$) by Theorem \ref{thm:Ringel-Green}. Then $\tMHk$ is generated by $[S_i]$ ($i\in\I$) and $\E_\alpha$ ($\alpha\in\Z^{\I}$) since $\tMHk\cong \widetilde{\ch}(\bfk Q)*\widetilde{\ct}(\Lambda^\imath)$ as linear spaces. So $\widetilde{\psi}$ is surjective. The theorem is proved.
	\end{proof}

	\subsection{Computations for rank 2 Dynkin $\imath$quivers, I}
	\label{subsec:computationrank2I}

	\begin{proposition}
		\label{prop:A2}
		Let $Q$ be the quiver $1\longrightarrow 2$, with $\btau=\Id$. Then in $\tMH$ we have
		\begin{align}
			[S_2]*[S_1]*[S_1]-({\sqq}+{\sqq}^{-1})[S_1]*[S_2]*[S_1]+[S_1]*[S_1]*[S_2]=-\frac{(q-1)^2}{{\sqq}}[S_2]*[ \E_1],
			\label{eq:S211}
			\\
			[S_1]*[S_2]*[S_2]-({\sqq}+{\sqq}^{-1})[S_2]*[S_1]*[S_2]+[S_2]*[S_2]*[S_1]=-\frac{(q-1)^2}{{\sqq}}[S_1]*[ \E_2].
			\label{eq:S122}
		\end{align}
	\end{proposition}
	
	\begin{proof}
		Recall from Example \ref{example 2}(a) the quiver and relations of $\Lambda^{\imath}$ for this $\imath$quiver.  We shall only prove the first identity \eqref{eq:S211} while skipping a similar proof of the identity \eqref{eq:S122}.
		
		Denote by $U_i$ the indecomposable projective $\Lambda^{\imath}$-module corresponding to $i$. Denote by $X$ the unique indecomposable $\Lambda^{\imath}$-module with $\widehat{S_1}+\widehat{S_2}$ as its class in $K_0(\mod(\Lambda^{\imath}))$. Then we have
		\begin{align*}
			[S_2]*[S_1]*[S_1]&= {\sqq}^{\langle S_2,S_1\rangle_Q} [S_1\oplus S_2]*[S_1]\\
			&= {\sqq}^{\langle S_2,S_1\rangle_Q} {\sqq}^{\langle S_1\oplus S_2,S_1\rangle_Q}
			\Big(\frac{1}{q} [S_2\oplus S_1\oplus S_1] +\frac{q-1}{q}[S_2\oplus \E_1] \Big)\\
			&= \frac{1}{{\sqq}} [S_2\oplus S_1\oplus S_1] +\frac{q-1}{{\sqq}}[S_2\oplus \E_1];
		\end{align*}
		\begin{align*}
			[S_1]*[S_2]*[S_1]&= {\sqq}^{\langle S_1,S_2\rangle_Q} ([S_1\oplus S_2]*[S_1]+(q-1)[X] *[S_1]) \\
			&= {\sqq}^{\langle S_1,S_2\rangle_Q} {\sqq}^{\langle S_1\oplus S_2,S_1\rangle_Q}
			\Big( \frac{1}{q}[S_1\oplus S_1\oplus S_2] +\frac{q-1}{q}[S_2\oplus \E_1] \Big)\\
			& \quad +{\sqq}^{\langle S_1,S_2\rangle_Q} {\sqq}^{\langle S_1\oplus S_2,S_1\rangle_Q}
			\Big(\frac{q-1}{q}[X\oplus S_1]+ \frac{(q-1)^2}{q}[U_1/S_2] \Big)\\
			&= \frac{1}{q}[S_1\oplus S_1\oplus S_2] +\frac{q-1}{q}[S_2\oplus \E_1]
			\\
			&\qquad\qquad\qquad\qquad\, +\frac{q-1}{q}[X\oplus S_1]+ \frac{(q-1)^2}{q}[U_1/S_2];
		\end{align*}
		\begin{align*}
			[S_1]*[S_1]*[S_2]&= [S_1]* \big({\sqq}^{\langle S_1,S_2\rangle_Q} ([S_1\oplus S_2]+(q-1)[X]) \big) \\
			&= {\sqq}^{\langle S_1,S_2\rangle_Q} {\sqq}^{\langle S_1,S_1\oplus S_2\rangle_Q}  \\
			&\qquad
			\times \Big(\frac{1}{q} [S_1\oplus S_1\oplus S_2]+\frac{q-1}{q}[\E_1\oplus S_2]+\frac{q-1}{q}[S_1\oplus X] \Big)\\
			& \quad +{\sqq}^{\langle S_1,S_2\rangle_Q} {\sqq}^{\langle S_1,S_1\oplus S_2\rangle_Q}
			\Big(\frac{(q-1)^2}{q}[U_1/S_2]+(q-1)[S_1\oplus X] \Big)\\
			&= \frac{1}{q{\sqq}} [S_1\oplus S_1\oplus S_2]+\frac{q-1}{q{\sqq}}[\E_1\oplus S_2]
			+\frac{q^2-1}{q{\sqq}}[S_1\oplus X]+\frac{(q-1)^2}{q{\sqq}}[U_1/S_2].
		\end{align*}
		
		By the short exact sequence $0\rightarrow S_2\rightarrow U_1/S_2\rightarrow
		\E_1\rightarrow 0$ with $\E_1 \in P^{\leq 1}(\La^\imath)$, we have
		$[U_1/S_2]=[\E_1\oplus S_2]$ in $\tMH$, and
		then by Lemma \ref{lem:Euler}, we have
		\[
		[S_2]*[ \E_1]={\sqq}^{\langle S_2,S_1\oplus S_1\rangle_Q} q^{-\langle S_2,\E_1\rangle}[U_1/S_2]=[U_1/S_2].
		\]
		Hence, in $\tMH$ we obtain that
		\begin{align*}
			[S_2]*[S_1]*[S_1]  - & ({\sqq}+{\sqq}^{-1})[S_1]* [S_2]*[S_1]+[S_1]*[S_1]*[S_2]
			\\
			&= -\frac{(q-1)^2}{{\sqq}}[U_1/S_2]
			= -\frac{(q-1)^2}{{\sqq}}[S_2]*[ \E_1].
		\end{align*}
		The proposition is proved.
	\end{proof}

	\subsection{Computations for rank 2 Dynkin $\imath$quivers, II}
	\label{subsec:computationrank2II}

	\begin{proposition}
		\label{prop:iA3}
		Let $Q$ be the quiver $1\longrightarrow 2\longleftarrow3$ with $\btau$ being the nontrivial involution. Then in $\tMH$ we have, for $i=1, 3$,
		\begin{align}
			[S_i]*[S_i]*[S_2]-({\sqq}+{\sqq}^{-1})[S_i]*[S_2]*[S_i]+[S_2]*[S_i]*[S_i] &=0,
			\label{eq:Serre112} \\
			[S_2]*[S_2]*[S_i]-({\sqq}+{\sqq}^{-1})[S_2]*[S_i]*[S_2]+[S_i]*[S_2]*[S_2]
			& =-\frac{(q-1)^2}{{\sqq}}[ S_i]*[\E_2].      \label{eq:Serre221}
		\end{align}
	\end{proposition}
	(Using the opposite quiver $1\longleftarrow 2\longrightarrow3$ yields the same formulas above.)
	
	\begin{proof}
		Recall from Example \ref{example 2}(b) the quiver and relations of $\Lambda^{\imath}$.
		Denote by $U_i$ the indecomposable projective $\Lambda^{\imath}$-module corresponding to $i \in \{1,2,3\}$. Denote by $X$ the unique indecomposable $\Lambda^{\imath}$-module with $\widehat{S_1}+\widehat{S_2}$ as its class in $K_0(\Lambda^{\imath})$.
		We shall only prove the formulas for $i=1$, as the remaining case with $i=3$ follows by symmetry.
		
		In $\tMH$, we have
		\begin{align*}
			[S_1]*[S_1]*[S_2]&= {\sqq}^{\langle S_1,S_2\rangle_Q}[S_1]* \big([S_1\oplus S_2] + (q-1)[X] \big) \\
			&= {\sqq}^{\langle S_1,S_2\rangle_Q} {\sqq}^{\langle S_1,S_1\oplus S_2\rangle_Q}   \\
			&\qquad \times \Big(\frac{1}{q}[S_1\oplus S_1\oplus S_2]+ \frac{q-1}{q}[S_1\oplus X]+(q-1)[S_1\oplus X] \Big)\\
			&= \frac{1}{q{\sqq}}[S_1\oplus S_1\oplus S_2]+ \frac{q^2-1}{q{\sqq}}[S_1\oplus X];
		\end{align*}
		\begin{align*}
			[S_1]*[S_2]*[S_1]&= {\sqq}^{\langle S_1,S_2\rangle_Q} \big([S_1\oplus S_2] + (q-1)[X] \big)*[S_1]\\
			&= {\sqq}^{\langle S_1,S_2\rangle_Q} {\sqq}^{\langle S_1\oplus S_2,S_1\rangle_Q}
			\Big(\frac{1}{q}[S_1\oplus S_1\oplus S_2]+ \frac{q-1}{q}[S_1\oplus X] \Big)\\
			&= \frac{1}{q}[S_1\oplus S_1\oplus S_2]+ \frac{q-1}{q}[S_1\oplus X];
		\end{align*}
		and
		\begin{align*}
			[S_2]*[S_1]*[S_1]&= {\sqq}^{\langle S_2,S_1\rangle_Q}[S_1\oplus S_2]*[S_1]\\
			&= {\sqq}^{\langle S_2,S_1\rangle_Q} {\sqq}^{\langle S_1\oplus S_2,S_1\rangle_Q} \frac{1}{q}[S_1\oplus S_1\oplus S_2]\\
			&= \frac{1}{{\sqq}}[S_1\oplus S_1\oplus S_2].
		\end{align*}
		The first identity \eqref{eq:Serre112} follows from combining the above computations.
		
		On the other hand, we have
		\begin{align*}
			[S_2]*[S_2]*[S_1]&= {\sqq}^{\langle S_2,S_1\rangle_Q} [S_2]*[S_2\oplus S_1]\\
			&= {\sqq}^{\langle S_2,S_1\rangle_Q}{\sqq}^{\langle S_2,S_1\oplus S_2\rangle_Q}
			\Big(\frac{1}{q}[S_2\oplus S_2\oplus S_1]+\frac{q-1}{q}[\E_2\oplus S_1] \Big)\\
			&= \frac{1}{{\sqq}}[S_2\oplus S_2\oplus S_1]+\frac{q-1}{{\sqq}}[\E_2\oplus S_1];
		\end{align*}
		\begin{align}
			\label{eqn:S212}
			[S_2]*[S_1]*[S_2]&= {\sqq}^{\langle S_2,S_1\rangle_Q} [S_2\oplus S_1]*[S_2]\\
			&= {\sqq}^{\langle S_2,S_1\rangle_Q}{\sqq}^{\langle S_2\oplus S_1,S_2\rangle_Q} \notag
			\Big(\frac{1}{q}[S_2\oplus S_1\oplus S_2]+\frac{q-1}{q}[\E_2\oplus S_1] \\\notag
			& \qquad +\frac{q-1}{q}[S_2\oplus X]+\frac{(q-1)^2}{q}[\rad (U_3)] \Big)\\\notag
			&= \frac{1}{q}[S_2\oplus S_1\oplus S_2]+\frac{q-1}{q}[\E_2\oplus S_1]
			+\frac{q-1}{q}[S_2\oplus X]+\frac{(q-1)^2}{q}[\rad (U_3)] \notag
			\\
			&= \frac{1}{q}[S_2\oplus S_1\oplus S_2]+(q-1)[\E_2\oplus S_1] +\frac{q-1}{q}[S_2\oplus X],\notag
		\end{align}
		where we have used $[\rad(U_3)]=[\E_2\oplus S_1]$ in $\tMH$ thanks to a short exact sequence $0\rightarrow \E_2\rightarrow \rad(U_3)\rightarrow S_1\rightarrow0$ with $\E_2 \in P^{\leq 1}(\La^\imath)$. In addition, \eqref{eqn:S212} also implies that
		\[
		[S_1\oplus S_2]*[S_2]=\frac{1}{q}[S_2\oplus S_1\oplus S_2]+(q-1)[\E_2\oplus S_1] +\frac{q-1}{q}[S_2\oplus X].
		\]
		Then we have
		\begin{align*}
			[S_1]*[S_2]*[S_2]&= {\sqq}^{\langle S_1,S_2\rangle_Q} ([S_1\oplus S_2]*[S_2]+ (q-1)[X]*[S_2])\\
			&= {\sqq}^{\langle S_1,S_2\rangle_Q}{\sqq}^{\langle S_2\oplus S_1,S_2\rangle_Q}
			\Big(\frac{1}{q}[S_2\oplus S_1\oplus S_2]+(q-1)[\E_2\oplus S_1]  \\
			& \qquad\quad +\frac{q-1}{q}[S_2\oplus X] +(q-1)[X\oplus S_2] \Big)\\
			&= \frac{1}{q{\sqq}}[S_2\oplus S_1\oplus S_2]+\frac{q-1}{{\sqq}}[\E_2\oplus S_1] +\frac{q^2-1}{q{\sqq}}[S_2\oplus X].
		\end{align*}
		Summarizing, we have obtained
		\begin{align*}
			[S_2]*[S_2]*[S_1]& -({\sqq}+{\sqq}^{-1})[S_2]*[S_1]*[S_2]+[S_1]*[S_2]*[S_2] \\
			&=  -\frac{(q-1)^2}{{\sqq}}[\E_2\oplus S_1]
			= -\frac{(q-1)^2}{{\sqq}}[ S_1]*[\E_2],
		\end{align*}
		whence the identity~\eqref{eq:Serre221}.
	\end{proof}

	\section{$\imath$Hall algebras for acyclic $\imath$quivers}
	\label{sec:iHallKM}
	
	In this section, we establish the $\imath$Hall algebra realization of the universal $\imath$quantum group of Kac-Moody type. The relations \eqref{relation6} and \eqref{relation5} in $\tUi$ are highly nontrivial in contrast to the Dynkin case, and their verification in the $\imath$Hall algebra setting is highly nontrivial and requires different techniques. We outline the proof of the relation \eqref{relation6} for $\imath$Hall algebra, by first formulating the $\imath$Hall algebra counterpart of $\imath$divided powers. 

	\subsection{Rank 1 and $\imath$divided powers}
	
	Consider the $\imath$quiver consisting of a single vertex with a trivial involution. The associated $\imath$quiver algebra is given by $\Lambda^\imath = \bfk[x]/(x^2)$, and the corresponding $\imath$Hall algebra $\tMHk$ is commutative and isomorphic to $\Q(\sqq)[[S],\bk^{\pm1}]$. The split $\imath$quantum group $\tUi$ of rank one is the algebra $\Q(v)[\ff,\tk^{\pm1}]$.
	The $\Q(\sqq)$-algebra isomorphism
	$\widetilde{\psi}:\tUi_{\sqq}\longrightarrow \tMHk$ (of rank one) (see Theorem~\ref{thm:main}) sends
	\begin{align}
		\label{rank1map}
		\ff\mapsto \frac{-1}{q-1}[S], \qquad \tk\mapsto -\frac{\bK}{q}.
	\end{align}

	Inspiring by \eqref{eq:iDPodd}--\eqref{eq:iDPev} and \eqref{rank1map}, we define the $\imath$-divided power of $[S]$ in $\tMHk$ to be
	\begin{align*}
		&[S]_{\odd}^{(m)}:=\frac{1}{[m]_{\sqq}!}\left\{ \begin{array}{ll} [S]\prod_{j=1}^k ([S]^2+\sqq^{-1}(\sqq^2-1)^2[2j-1]_{\sqq}^2 [\bK] ) & \text{if }m=2k+1,\\
			\prod_{j=1}^k ([S]^2+\sqq^{-1}(\sqq^2-1)^2[2j-1]_{\sqq}^2[\bK]) &\text{if }m=2k; \end{array}\right.
		\\
		&[S]_{\ev}^{(m)}:= \frac{1}{[m]_{\sqq}!}\left\{ \begin{array}{ll} [S]\prod_{j=1}^k ([S]^2+\sqq^{-1}(\sqq^2-1)^2[2j]_{\sqq}^2[\bK] ) &\text{if }m=2k+1,\\
			\prod_{j=1}^{k} ([S]^2+\sqq^{-1}(\sqq^2-1)^2[2j-2]_{\sqq}^2[\bK]) &\text{if }m=2k. \end{array}\right.
	\end{align*}
	
	These $\imath$-divided powers satisfy the following recursive relations:
	\begin{align}
		[S]*[S]^{(2m)}_{\odd} &=[2m+1] [S]^{(2m+1)}_{\odd},
		\label{rec:2modd}
		\\
		[S]*[S]^{(2m+1)}_{\odd} &=[2m+2] [S]^{(2m+2)}_{\odd} -\sqq (\sqq-\sqq^{-1})^2[2m+1] [S]^{(2m)}_{\odd}*[\bK],
		\label{rec:2m+1odd}
		\\
		[S]*[S]^{(2m-1)}_{\ev} &=[2m][S]^{(2m)}_{\ev},
		\label{rec:2m-1ev}
		\\
		[S]*[S]^{(2m)}_{\ev} &=[2m+1] [S]^{(2m+1)}_{\ev} -\sqq (\sqq-\sqq^{-1})^2 [2m] [S]^{(2m-1)}_{\ev}*[\bK].
		\label{rec:2mev}
	\end{align}
	
	The following lemma follows by definition.
	
	\begin{lemma}
		The isomorphism $\widetilde{\psi}$ satisfies that, for $m\in \N$,
		\begin{align}
			\widetilde{\psi}(\ff_{\odd}^{(m)})=\frac{[S]_{\odd}^{(m)}}{(1-\sqq^2)^m},
			\qquad
			\widetilde{\psi}(\ff_{\ev}^{(m)})=\frac{[S]_{\ev}^{(m)}}{(1-\sqq^2)^m}.
		\end{align}
	\end{lemma}

	We denote by $[0]_\sqq ^{^{!!}}=1$, and for any $k\in \Z_{\ge 1}$,
	\[
	[2k]_\sqq^{^{!!}}=[2k]_\sqq[2k-2]_\sqq \cdots [4]_\sqq[2]_\sqq.
	\]
	We denote by $\lfloor x \rfloor$ the largest integer not exceeding $x$, for $x\in \mathbb R$.
	
	\begin{proposition}
		\label{prop:iDPev}
		The following identity holds in $\tMHk$, for $n\in \N$:
		\begin{align}
			\label{eq:nev}
			&[S]^{(n)}_{\ev}=\sum_{k=0}^{\lfloor \frac{n}2\rfloor} \frac{\sqq^{k(k - (-1)^n)-\binom{n-2k}{2}}\cdot(\sqq-\sqq^{-1})^k}{[n-2k]_{\sqq}^{!}[2k]_\sqq^{!!}}[(n-2k)S]*[\bK]^k,
			\\
			\label{eq:nodd}
			&[S]^{(n)}_{\odd}=\sum_{k=0}^{\lfloor \frac{n}2\rfloor}  \frac{\sqq^{k(k+(-1)^n)-\binom{n-2k}{2}}\cdot(\sqq-\sqq^{-1})^k}{[n-2k]_{\sqq}^{!}[2k]_\sqq^{!!}}[(n-2k)S]*[\bK]^k.
		\end{align}
	\end{proposition}
	
	\begin{proof}
		The required Euler form is given by $\langle S,S^{\oplus m}\rangle_Q=m=\dim_\bfk\Hom_{\Lambda^\imath}(S,S^{\oplus m})$. For any non-split short exact sequence
		$0\rightarrow S^{\oplus m} \rightarrow M\rightarrow S\rightarrow0$
		in $\mod(\Lambda^\imath)$, we have $M\cong S^{\oplus(m-1)}\oplus \bK$. Note that $\Ext^1_{\Lambda^\imath}(S,S^{\oplus m})=m$. Then we have
		\begin{align}
			\label{SmS}
			[S]*[mS]&= \sqq^{-m}[(m+1)S]+\sqq^{-m}(q^m-1) [(m-1)S\oplus \bK]
			\\
			&= \sqq^{-m}[(m+1)S]+(\sqq^m-\sqq^{-m})[(m-1)S]*[\bK].
			\notag
		\end{align}
		
		The identities  \eqref{eq:nev}--\eqref{eq:nodd} can now be proved by induction on $n$ using the recursive identities \eqref{rec:2modd}--\eqref{rec:2mev} and \eqref{SmS}.
	\end{proof}

	\subsection{Rank 2 computations}
	\label{subsec:rank2KM}
	
	The following multiplication formula in the $\imath$Hall algebra is useful.
	
	\begin{proposition}[\text{\cite[Proposition 3.10]{LW20}}]
		\label{prop:iHallmult}
		Let $(Q,\tau)$ be an $\imath$quiver with $\tau=\Id$.
		For any $A,B\in\mod(\bfk Q)\subset \mod(\Lambda^\imath)$, the following identity holds in $\tMHk$. 
		\begin{align}
			\label{Hallmult1}
			[A]*[B]&= 
			\sum_{[L],[M],[N]\in\Iso({\mod(\bfk Q)})} \sqq^{\langle A,B\rangle_Q}  q^{\langle N,L\rangle_Q -\langle A,B\rangle_Q}\frac{|\Ext^1(N, L)_{M}|}{|\Hom(N,L)| }
			\\
			\notag
			&\qquad \times |\{s\in\Hom(A,B)\mid \Ker s\cong N, \coker s\cong L
			\}|\cdot [M]*[\bK_{\widehat{A}-\widehat{N}}].
		\end{align}
	\end{proposition}

	Consider the $\imath$quiver
	\begin{align}
		\label{diag:split2}
		Q=(\xymatrix{1\ar@<1ex>[r]|-{a}& 2\ar@<1ex>[l]|-{b}}),
		\quad
		\btau=\Id,
		\qquad \text{where } a+b=-c_{12}.
	\end{align}
	The arrows $1\rightarrow 2$ are denoted by $\alpha_i$, $1\leq i\leq a$, and $2\rightarrow 1$ are denoted by $\beta_j$, $1\leq j\leq b$.
	Then the corresponding $\imath$quiver algebra $\Lambda^\imath$ has its quiver $\ov{Q}$ as
	\begin{center}\setlength{\unitlength}{0.8mm}
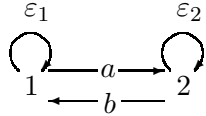
		\begin{equation*}
			\begin{picture}(50,14)(0,-8)
				\put(0,-3){$1$}
				\put(4,1){\line(1,0){8}}
				\put(12.5,0){$a$}
				\put(15.5,1){\vector(1,0){7.5}}
				\put(23,-4){\line(-1,0){7}}
				\put(12,-4){\vector(-1,0){8}}
				\put(13,-6){$b$}
				\put(25,-3){$2$}
				\qbezier(-1,1)(-3,3)(-2,5.5)
				\qbezier(-2,5.5)(1,9)(4,5.5)
				\qbezier(4,5.5)(5,3)(3,1)
				\put(3.1,1.4){\vector(-1,-1){0.3}}
				\qbezier(24,1)(22,3)(23,5.5)
				\qbezier(23,5.5)(26,9)(29,5.5)
				\qbezier(29,5.5)(30,3)(28,1)
				\put(28.1,1.4){\vector(-1,-1){0.3}}
				\put(0,10){$\varepsilon_1$}
				\put(25,10){$\varepsilon_2$}
			\end{picture}
		\end{equation*}
	\end{center}

	For any $\Lambda^\imath$-module $M=(M_i, M(\alpha_j), M(\beta_k), M(\varepsilon_i))_{i=1,2;1\leq j\leq a,1\leq k\leq b}$ such that $\dim_\bfk M_2=1$, we define
	\begin{align}
		\label{eq:UW}
		U_M:= \bigcap_{1\leq i\leq a} \Ker M(\alpha_i),\qquad W_M:= \sum_{j=1}^b \Im M(\beta_j),
	\end{align}
	and let
	\begin{align}  \label{eq:uw}
		u_M:=\dim U_M,\qquad w_M:=\dim W_M.
	\end{align}
	For any $r\geq0$, define
	\begin{align}
		\cI_{r}:&= \{[M]\in\Iso(\mod(\bfk Q))\mid \widehat{M}=r\widehat{S_1}+\widehat{S_2}, W_M\subseteq U_M\},
		\\
		\cJ_r:&= \{[M]\in\Iso(\mod(\bfk Q))\mid \widehat{M}=r\widehat{S_1}+\widehat{S_2}, S_2\in\add \Top(M)\}.
	\end{align}
	
	The following proposition  generalizes \cite[Proposition 7.3]{LW20}, which treated the special cases when either $a$ or $b$ is $0$. 
	
	\begin{proposition}[cf. \text{\cite[Proposition 7.3]{LW20}}]
		\label{prop:SSSM}
		The following identity holds in $\tMHk$, for $s, t \ge 0$:
		\begin{align}
			\label{eq:SSS}
			[s S_1]*[S_2]*[t S_1] &= \sqq^{-tb}\sum_{r=0}^{\min\{s,t\}}\sum\limits_{[M]\in\cI_{s+t-2r}}\sqq^{\tilde{p}(a,b,r,s,t)}(\sqq-\sqq^{-1})^{s-r+t+1}
			\frac{[s]_\sqq^![t]_\sqq^!}{[r]_\sqq^!} 
			\\
			& \qquad\qquad\qquad  \times \qbinom{u_M-w_M}{(t-r)-w_M}_\sqq\frac{1}{|\Aut(M)|}[M]*[\E_1]^r,
			\notag
		\end{align}
		where
		\begin{align*}
			\tilde{p}(a,b,r,s,t)&= 2(s-r)(t-r-a)-s(t-a)+2br-2r(t-r)-r^2+rs+t^2+\tbinom{t}{2}+
			\\&(s-r)^2+\tbinom{s-r}{2}+\big(u_M-(t-r)\big)\big((t-r)-w_M\big)+1.
		\end{align*}
	\end{proposition}
	
	\begin{proof}
		The computation in cases either $a=0$ or $b=0$ can be found in the proof of \cite[Proposition 7.3]{LW20}.
		The computation required in the general case is similar, again by using crucially Proposition \ref{prop:iHallmult}; we skip the details. 
	\end{proof}
	
	We can formulate the $\imath$Serre relations in the $\imath$Hall algebra, and we postpone its sketched proof in the next subsection. 
	
	\begin{theorem}
		\label{thm:iSerreHall}
		The following identity holds in $\tMHk$, for any $\overline{p}\in\Z_2$:
		\begin{align}
			\sum_{n=0}^{1+a+b} (-1)^n  [S_1]_{\overline{p}}^{(n)}*[S_2] *[S_1]_{\overline{a}+\ov{b}+\overline{p}}^{(1+a+b-n)} =0.
			\label{eqn:iserre1}
		\end{align}
	\end{theorem}

	\subsection{Proof of Theorem~\ref{thm:iSerreHall}}
	\label{subsec:rank2KMproof}
	
	In this subsection, we sketch the proof of Theorem~\ref{thm:iSerreHall}. 
	
	Note that \eqref{eqn:iserre1} can be written as
	\begin{align}
		\sum_{n=0}^{1+a+b} (-1)^n  [S_1]_{\overline{0}}^{(n)}*[S_2] *[S_1]_{\overline{a}+\ov{b}}^{(1+a+b-n)} =0,
		\label{eqn:iserre10}
		\\
		\sum_{n=0}^{1+a+b} (-1)^n  [S_1]_{\overline{1}}^{(n)}*[S_2] *[S_1]_{\overline{a}+\ov{b}+\ov{1}}^{(1+a+b-n)} =0.
		\label{eqn:iserre11}
	\end{align}
	We only sketch the proof of \eqref{eqn:iserre10} here, since the other one is similar.
	
	Using \eqref{eq:nev}--\eqref{eq:nodd} and Proposition \ref{prop:SSSM}, we obtain
	\begin{align}
		&\sum_{n=0}^{a+b+1}(-1)^n[S_1]_{\overline{0}}^{(n)}*[S_2] *[S_1]_{\overline{a}+\ov{b}}^{(1+a+b-n)}\label{eqn:iserre111}
		\\ \notag &= \sum_{n=0,2\mid n}^{a+b+1}\sum_{k=0}^{\frac{n}{2}} \sum_{m=0}^{\lfloor\frac{a+b+1-n}{2}\rfloor} \sqq^{-(1+a+b-n-2m)b} \sum_{r=0}^{\min\{n-2k,1+a+b-n-2m\}}\sum\limits_{[M]\in\cI_{1+a+b-2k-2m-2r}}
		\\ \notag &\quad\quad\frac{\sqq^{\tilde z}(\sqq-\sqq^{-1})^{2+a+b-k-m-r}}{[2k]_\sqq^{!!}[2m]_\sqq^{!!}[r]_\sqq^{!}}\begin{bmatrix} u_M-w_M \\ 1+a+b-n-2m-r-w_M \end{bmatrix}_\sqq
		\frac{[M]*[\E_1]^{k+m+r}}{|\Aut(M)|}
		\\ \notag &-\sum_{n=0,2\nmid n}^{a+b+1}\sum_{k=0}^{\frac{n-1}{2}} \sum_{m=0}^{\lfloor\frac{a+b+1-n}{2}\rfloor} \sqq^{-(1+a+b-n-2m)b} \sum_{r=0}^{\min\{n-2k,1+a+b-n-2m\}}\sum\limits_{[M]\in\cI_{1+a+b-2k-2m-2r}}
		\\ \notag &\quad\frac{\sqq^{\tilde z+2k-2m}(\sqq-\sqq^{-1})^{2+a+b-k-m-r}}{[2k]_\sqq^{!!}[2m]_\sqq^{!!}[r]_\sqq^{!}}\begin{bmatrix} u_M-w_M \\ 1+a+b-n-2m-r-w_M \end{bmatrix}_\sqq \frac{[M]*[\E_1]^{k+m+r}}{|\Aut(M)|}.
	\end{align}

	Set
	\begin{align*}
		d=k+m+r.
	\end{align*}
	Now we only need to prove that the coefficient of $\frac{[M]*[\E_1]^d}{|\Aut(M)|}$ in the \text{RHS} of \eqref{eqn:iserre111} is zero, for any given $[M]\in\cI_{1+a+b-2d}$ and any $d\in \mathbb{N}$. For simplicity, we denote by $u=u_M$, $w=w_M$ in the following.

	Note the powers of $(\sqq-\sqq^{-1})$ in all terms are the same (equals to $2+a+b-d$). Denote
	\begin{align}
		\widetilde{T}(a,b,d,u,w)=\sum_{n=0,2\mid n}^{a+b+1}\sum_{k=0}^{\frac{n}{2}} \sum_{m=0}^{\lfloor\frac{a+b+1-n}{2}\rfloor}& \sqq^{-(1+a+b-n-2m)b}\delta\{0\leq r\leq n-2k\}
		\label{eqn:identityT}
		\\ \notag &\times \frac{\sqq^{\tilde z}}{[2k]_\sqq^{!!}[2m]_\sqq^{!!}[r]_\sqq^{!!}}\begin{bmatrix} u-w \\ 1+a+b-n-2m-r-w \end{bmatrix}_\sqq
		\\ \notag -\sum_{n=0,2\nmid n}^{a+b+1}\sum_{k=0}^{\frac{n-1}{2}} \sum_{m=0}^{\lfloor\frac{a+b+1-n}{2}\rfloor}& \sqq^{-(1+a+b-n-2m)b}\delta\{0\leq r\leq n-2k\}
		\\ \notag &\times \frac{\sqq^{\tilde z+2k-2m}}{[2k]_\sqq^{!!}[2m]_\sqq^{!!}[r]_\sqq^{!!}}\begin{bmatrix} u-w \\ 1+a+b-n-2m-r-w \end{bmatrix}_\sqq,
	\end{align}
	where we set $\delta\{X\}=1$ if the statement $X$ holds and $\delta\{X\}=0$ if $X$ is false. Then the coefficient of $\frac{[M]*[\E_1]^d}{|\Aut(M)|}$ in the \text{RHS} of \eqref{eqn:iserre111} is equal to $(\sqq-\sqq^{-1})^{2+a+b-d}\widetilde{T}(a,b,d,u,w)$. So the identity \eqref{eqn:iserre10} is equivalent to the identity $\widetilde{T}(a,b,d,u,w)=0$, for any nonnegative integers $a,b,d,u,w$ subject to the constraints
	\begin{align}
		0\leq d\leq (a+b+1)/2,\quad0\leq w\leq b,\quad b+1-2d\leq u\leq 1+a+b-2d,
		\label{eqn:identityTT}
	\end{align}
	where $u,d$ are not both zeros.

	\begin{lemma}[\text{\cite[Proposition 8.1]{LW20}}]
		\label{lem:identity=0}
		For any nonnegative integers $a,b,d,u,w$ satisfying the constraint \eqref{eqn:identityTT}, the following identity holds:
		\begin{align}
			\widetilde{T}(a,b,d,u,w)=0.
		\end{align}
	\end{lemma}
	
	The proof of Lemma \ref{lem:identity=0} can be reduced to the following identities: for $p, d\ge 1$,
	\begin{align*}
		\sum_{k=0}^{p} v^{-k(p-k+1)} \qbinom{p}{k}  & =  \prod_{j=1}^p (1+v^{-j}),
		\qquad
		\sum_{\stackrel{k,m,r \in \N}{k+m+r =d}} (-1)^r \frac{v^{{r+1 \choose 2} -2(k-1)m}}{[r]^! [2k]^{!!} [2m]^{!!}} =0.
	\end{align*}
	
	Thanks to Lemma \ref{lem:identity=0}, the proof of the identity~ \eqref{eqn:iserre10} is complete. This proves Theorem~\ref{thm:iSerreHall}.

	
	\subsection{$\imath$Hall algebra realization of $\tUi$ of Kac-Moody type}
	
	\begin{theorem}[\text{\cite[Theorem 9.6]{LW20}}]
		\label{thm:mainKM}
		Let $(Q, \btau)$ be an $\imath$quiver without loops. Then there exists a $\Q(\sqq)$-algebra monomorphism
		\begin{align*}
			\widetilde{\psi}: \tUi_{\sqq} &\longrightarrow \tMHk,
		\end{align*}
		which sends
		\begin{align}
			B_j \mapsto \frac{-1}{q-1}[S_{j}],\text{ if } j\in\ci,
			&\qquad\qquad
			\tk_i \mapsto - q^{-1}[\bK_i], \text{ if }\btau i=i \in \I;
			\label{eq:split}
			\\
			B_{j} \mapsto \frac{{\sqq}}{q-1}[S_{j}],\text{ if }j\notin \ci,
			&\qquad\qquad
			\tk_i \mapsto \sqq^{\frac{-c_{i,\btau i}}{2}}[\bK_i],\quad \text{ if }\btau i\neq i \in \I.
			\label{eq:extra}
		\end{align}
	\end{theorem}
	Due to the constraint arising from definition of $\imath$quivers, the Cartan matrix $(c_{ij})_{i,j\in \I}$ for Theorem~\ref{thm:mainKM} must satisfy the condition $c_{i,\tau i} \in 2 \Z$, for all $i\in \I$. 
	
	\begin{proof}
		Let us outline the proof. 
		First, we prove that $\widetilde{\psi}$ is an algebra homomorphism by verifying that it preserves the relations \eqref{relation1}--\eqref{relation5}. 
		
		The verification of relations \eqref{relation1}--\eqref{relation2} is the same as for Theorem \ref{thm:main}.
		The verification of relation \eqref{relation3} can be reduced to $\widetilde{\ch}(\bfk Q)$, which holds by Theorem \ref{thm:Ringel-Green}.
		
		The verification of relation \eqref{relation5} requires long homological computations (compare to Theorem \ref{thm:Ringel-Green}), and the detail can be found in \cite[\S 5]{LW20}.
		
		The verification of relation \eqref{relation6} can be reduced to the rank $2$ $\imath$quivers, while the highly nontrivial rank $2$ computation was sketched in \S\ref{subsec:rank2KM}--\ref{subsec:rank2KMproof}.
		
		The injectivity of $\widetilde{\psi}$ is proved in the same way as for Theorem \ref{thm:main}. This completes the proof of the theorem. 
	\end{proof}

	
	\section{Reflection functors in $\imath$Hall algebras and relative braid group symmetries}
	\label{sec:braid}
	
	In this section, we formulate and compute the reflection functors in the framework of $\imath$Hall algebras. The reflection functors give rise to (relative braid group) automorphisms of $\imath$quantum groups, and we provide explicit formulas for these symmetries on generators. 
	
	\subsection{Modulated graphs for $\imath$quivers}
	\label{subsec:modulated}
	
	In this subsection, we study representations of modulated graphs associated to $\imath$quivers, cf.  \cite{LW20,LW22}. 
	
	Let $(Q, \btau)$ be an $\imath$quiver, and $\Lambda^\imath=\bfk\ov{Q}/\ov{I}$ with $(\ov{Q}, \ov{I})$ as defined in \S\ref{subsec:i-quivers}.
	For each $i\in Q_0$, define a $\bfk$-algebra
	\begin{align}\label{dfn:bHi}
		{\ov{\BH}} _i:=\left\{ \begin{array}{ll}  \bfk[{\color{purple}{\varepsilon_i}}]/(\textcolor{purple}{\varepsilon_i}^2) & \text{ if }i=\btau i,
			\\
			\bfk(\xymatrix{ i  \ar@/^0.5pc/@<0.8ex>@[purple][r]^{\color{purple}{\varepsilon_i}} \ar@<0.75ex>[r]|-{r_i} & \btau i\ar@<0.75ex>[l]|-{r_i}  \ar@/^0.5pc/@<0.8ex>@[purple][l]^{\color{purple}{\varepsilon_{\btau i}}} })/( \varepsilon_i\varepsilon_{\btau i},\varepsilon_{\btau i}\varepsilon_i, \alpha_j\varepsilon_{\btau i}-\varepsilon_i\beta_{j},\beta_j\varepsilon_i-\varepsilon_{\btau i}\alpha_j\mid 1\leq j\leq r)  &\text{ if } \btau i \neq i ,  n_{i,\btau i}=2r_i.\end{array}\right.
	\end{align}
	Define the following subalgebra of $\Lambda^{\imath}$:
	\begin{equation}  \label{eq:H}
		\ov{\BH} =\bigoplus_{i\in \ci }\ov{\BH} _i.
	\end{equation}
	
	Define
	\begin{align}
		\label{eqn:orientation}
		\Omega:=\Omega(Q) =\{(i,j) \in Q_0\times Q_0\mid  \exists (\alpha:i\rightarrow j)\in Q_1, \btau i\neq j\}.
	\end{align}
	Then $\Omega$ represents the orientation of $Q$. 
	We also use $\Omega(i,-)$ to denote the subset $\{j\in Q_0 \mid \exists (\alpha:i\rightarrow j)\in Q_1\}$, and $\Omega(-,i)$ is defined similarly.
	
	For any $(i,j)\in \Omega$, we define
	\begin{equation}  \label{eq:jHi}
		{}{}_j{\ov{\BH}}_i := {\ov{\BH}} _j\Span_\bfk\{\alpha,\btau\alpha\mid(\alpha:i\rightarrow j)\in Q_1\text{ or } (\alpha:i\rightarrow \btau j)\in Q_1\}{\ov{\BH}} _i.
	\end{equation}
	Note that ${}{}_j{\ov{\BH}}_i ={}_{\btau j} {\ov{\BH}} _{\btau i}={}_{j} {\ov{\BH}} _{\btau i}={}_{\btau j} {\ov{\BH}} _{i}$ for any $(i,j)\in \Omega$.

	Hence ${}_j\oH_i $ is an $\oH _j\mbox{-}\oH _i$-bimodule, which is free as a left $\oH _j$-module (and respectively, right $\oH _i$-module), with a basis $_j\LL_i$ (and respectively, $_j\RR_i$) defined in the following.
	\begin{eqnarray}
		\label{basis of Hij left}
		_j\LL_i&=&\left\{ \begin{array}{cc}
			\{\alpha \mid (\alpha:i\rightarrow j)\in Q_1\} & \text{ if }i=\btau i, \btau j=j,\\
			\{\alpha+\btau\alpha \mid(\alpha:i\rightarrow j)\in Q_1\} & \text{ if }i=\btau i, \btau j\neq j,\\
			\{\alpha,\btau \alpha \mid (\alpha:i\rightarrow j)\in Q_1\} & \text{ if }i\neq \btau i,\btau j=j,\\
			\{\alpha+\btau \alpha \mid (\alpha:i\rightarrow j) \text{ or }(\alpha:i\rightarrow \btau j)\in Q_1\} & \text{ if }i\neq \btau i,\btau j\neq j;\label{eqn:basis of L}
		\end{array}\right.\\
		_j\RR_i&=& \left\{ \begin{array}{cc}\label{basis of Hij right}
			\{\alpha \mid (\alpha:i\rightarrow j)\in Q_1\} & \text{ if }i=\btau i,  \btau j=j,\\
			\{\alpha,\btau \alpha \mid (\alpha:i\rightarrow j)\in Q_1\} & \text{ if }i=\btau i, \btau j\neq j,\\
			\{\alpha+\btau \alpha \mid (\alpha:i\rightarrow j)\in Q_1\} & \text{ if }i\neq \btau i,\btau j=j,\\
			\{\alpha+\btau \alpha \mid (\alpha:i\rightarrow j) \text{ or }(\alpha:i\rightarrow \btau j)\in Q_1\} & \text{ if }i\neq \btau i,\btau j\neq j.\end{array}\right. \label{eqn:basis of R}
	\end{eqnarray}
	
	%
	%
	
	Denote
	\begin{equation}  \label{eq:ovOmega}
		\overline{\Omega}:=\{(i,j)\in \ci \times \ci \mid (i,j)\in\Omega\text{ or }(i,\btau j)\in\Omega\}.
	\end{equation}
	Recall that $\ci $ is a (fixed) subset of $Q_0 =\I$ consisting of the representatives of $\btau$-orbits. The tuple $(\oH_i,\,_j\oH_i):=(\oH_i,\,_j\oH_i)_{i\in\ci ,(i,j)\in\ov{\Omega}}$  is called a \emph{modulation} of $(Q, \btau)$ and is denoted by $\cm(Q, \btau)$.
	
	A representation $(N_i,N_{ji}):=(N_i,N_{ji})_{i\in\ci ,(i,j)\in\ov{\Omega}}$ of $\cm(Q, \btau)$ is defined by assigning to each $i\in \ci$ a finite-dimensional ${\ov{\BH}} _i$-module $N_i$ and to each $(i,j)\in \overline{\Omega}$ an ${\ov{\BH}} _j$-morphism
	$N_{ji}:{}_j{\ov{\BH}}_i \otimes_{{\ov{\BH}} _i} N_i\rightarrow N_j$. A morphism $f:L\rightarrow N$ between representations $L=(L_i,L_{ji})$ and $N=(N_i,N_{ji})$ of $\cm(Q, \btau)$ is a tuple $f=(f_i)_{i\in \ci}$ of ${\ov{\BH}} _i$-morphisms $f_i:L_i\rightarrow N_i$ such that the following diagram is commutative for each $(i,j)\in\overline{\Omega}$:
	\[\xymatrix{_j{\ov{\BH}}_i\otimes_{{\ov{\BH}} _i} L_i \ar[r]^{1\otimes f_i} \ar[d]^{L_{ij}}&  _j{\ov{\BH}}_i \otimes_{{\ov{\BH}} _i} N_i\ar[d]^{N_{ij}}\\
		L_j\ar[r]^{f_j} & N_j}\]
	
	\begin{proposition}[\text{\cite[Proposition 2.16]{LW22}}]
		\label{prop:modulated representation}
		The categories $\rep(\cm(Q, \btau))$ and $\rep(\ov{Q},\ov{I})$ are isomorphic.
	\end{proposition}
	
	\subsection{Reflection functors}
	\label{subsec:APR}
	
	In this subsection, we shall introduce the reflection functors in the setting of $\imath$quivers; see \cite{LW21a,LW22b}.
	
	
	Let $Q^*$ be the quiver constructed from $Q$ by reversing all the arrows $\alpha:i\rightarrow j$ such that $\btau i\neq j$. For any $i,j\in\I$ such that $\btau i\neq j$, we have $(i,j)\in\Omega$ if and only if $(j,i)\in \Omega^*:=\Omega(Q^*)$. For any $\alpha:i\rightarrow j$ in $Q$ such that $\btau i\neq j$, denote by $\widetilde{\alpha}:j\rightarrow i$ the corresponding arrow in $Q^*$.   Then $\btau$ induces an involution $\btau^*$ of $Q^*$. Clearly, $\btau^* i=\btau i$ for any vertex $i\in Q_0$. Then similarly we can define $\Lambda^*=\bfk Q^*\otimes_\bfk R_2$, and an involution $\btau^{*\sharp}$ for $\Lambda^*$, and its $\btau^{*\sharp}$-fixed point subalgebra $(\Lambda^*)^\imath$. Note that $\oH$ 
	is also a subalgebra of $(\Lambda^*)^\imath$.
	
	It is worth noting that $\ci$ is also a set of representatives of $\btau^*$-orbits. In this way, one can define $\ov{\Omega}^*$ (cf. \eqref{eq:ovOmega} for $\ov{\Omega}$).
	
	For any $(j,i)\in\Omega^*$, we can define
	$_i\oH_j$ as follows:
	\begin{align*}
		_i\oH_j:= \oH_i\Span_\bfk\{\widetilde{\alpha},\btau^* \widetilde{\alpha} \mid(\widetilde{\alpha}:j\rightarrow i)\in Q^*_1\text{ or }(\widetilde{\alpha}:j\rightarrow \btau^* i)\in Q^*_1\}\oH_j;
	\end{align*}
	
	Recall from \eqref{basis of Hij left}--\eqref{basis of Hij right} the basis $_i\LL_j$ (and respectively, $_i\RR_j$) for ${}_i\oH_j $ as a left $\oH _i$-module (and respectively, right $\oH _j$-module).
	Let $_j\LL_i^*$ and $_j\RR_i^*$ be the dual bases of $\Hom_{\oH_j}(_j \oH_i, \oH_j)$ and $\Hom_{\oH_i}(_j\oH_i,\oH_i)$, respectively. Denote by $b^*$ the corresponding dual basis vector for any $b\in{_j\LL_i}$ or $b\in{_j\RR_i}$.
	
	Since ${}_i\oH_j$ and $\Hom_{\oH_j}(_j\oH_i,\oH_j)$ are right free $\oH_j$-modules with bases given by ${}_i\RR_j$ and ${}_j\LL_i^*$ respectively, there is a right $\oH_j$-module isomorphism
	\[
	\rho:\, _i\oH_j \longrightarrow \Hom_{\oH_j}(_j\oH_i,\oH_j)
	\]
	such that $\rho(\widetilde{b})=b^*$ for any $b\in{}_j\LL_i$.
	It is then routine to check that $\rho$ is actually an $\oH_i$-$\oH_j$-bimodule isomorphism.
	Similarly, there is an $\oH_i$-$\oH_j$-bimodule isomorphism
	$$\lambda: \, _i\oH_j\longrightarrow \Hom_{\oH_i}(_j\oH_i,\oH_i).$$
	These two isomorphisms satisfy that $\rho(_i\RR_j)={}_j\LL_i^*$ and $\lambda(_i\LL_j)={}_j\RR_i^*$. We sometimes identify the spaces $\Hom_{\oH_j}(_j\oH_i,\oH_j)$, $_i\oH_j$ and $\Hom_{\oH_i}(_j\oH_i,\oH_i)$ via $\rho$ and $\lambda$.
	
	If $N_j$ is an $\oH_j$-module, then we have a natural isomorphism of $\oH_i$-modules
	$$\Hom_{\oH_j}(_j\oH_i,N_j)\longrightarrow\, _i\oH_j\otimes_{\oH_j}N_j$$
	defined by
	$$f\mapsto \sum_{b\in_j\LL_i} b^*\otimes f(b).$$
	Furthermore, for any $\oH_i$-module $L_i$, there is a natural isomorphism of $\bfk$-vector spaces:
	$$\Hom_{\oH_j}(_j\oH_i\otimes_{\oH_i} L_i, N_j)\longrightarrow \Hom_{\oH_i}(L_i,\Hom_{\oH_j}(_j\oH_i,N_j)).$$
	Composing the two maps above, we obtain the following.
	
	\begin{lemma}[cf. \text{\cite[Lemma 3.2]{LW21a}}]
		\label{lem:ad}
		There exists a canonical $\bfk$-linear isomorphism
		\begin{align*}
			{\rm ad}_{ji}  ={\rm ad}_{ji}(L_i,N_j)  : & \Hom_{\oH_j}(_j\oH_i\otimes_{\oH_i} L_i,N_j)\longrightarrow \Hom_{\oH_i}(L_i,\,_i\oH_j\otimes_{\oH_j}N_j)
			\\
			{\rm ad}_{ji}: & f\mapsto \big(f^\vee:l\mapsto \sum_{b\in_j\LL_i} b^*\otimes f(b\otimes l) \big).
		\end{align*}
		The inverse ${\rm ad}_{ji}^{-1}$ is given by
		${\rm ad}_{ji}^{-1} (g) = \big(g^\vee:h\otimes l\mapsto \sum_{b\in_j\LL_i} b^*(h)(g(l))_b \big),$ where the elements $(g(l))_b\in N_j$ are uniquely determined by
		$g(l)=\sum_{b\in_j\LL_i} b^*\otimes(g(l))_b.$
	\end{lemma}

	Let $(Q, \btau)$ be an $\imath$quiver. Without loss of generality, we assume $Q$ to be connected and of rank $\ge 2$. Recall $\Omega=\Omega(Q)$ is the orientation of $Q$. For any sink $\ell \in Q_0$, define the quiver $s_\ell (Q)$ by reversing all the arrows ending to $\ell $. Note that in this case, we have $n_{\ell,\tau \ell}=0$.
	By definition, $\ell $ is a sink of $Q$ if and only if $\btau \ell $ is a sink of $Q$.  Define the quiver
	\begin{align}  \label{eq:QQ}
		Q' := \bs_\ell Q :=\left\{ \begin{array}{cc} s_\ell (Q) & \text{ if } \btau \ell =\ell ,\\ s_{\ell }s_{\btau\ell  }(Q) &\text{ if }\btau \ell \neq \ell .  \end{array}\right.
	\end{align}
	Note that $s_\ell s_{\btau \ell }(Q)=s_{\btau \ell }s_\ell (Q)$.
	Then $\btau$ also induces an involution $\btau$ on the quiver $Q'$. In this way, we can define $\Lambda'=\bfk Q'\otimes_\bfk R_2$ with an involution $\btau^\sharp$, and denote the $\btau^\sharp$-fixed point subalgebra by $\Lambda'^\imath =\bs_{\ell }\Lambda^{\imath}$. Note that
	$\bs_{\ell }\Lambda^{\imath}=\bs_{\btau \ell }\Lambda^{\imath}$
	for any sink $\ell \in Q_0$.
	The quiver $\ov{Q'}$ of $\bs_\ell \Lambda^{\imath}$ can be constructed from $\ov{Q}$ by reversing all the arrows ending to $\ell $ and $\btau \ell $. Denote by $\Omega':=\Omega(Q')$ the orientation of $Q'$.
	
	We shall define a {\rm reflection functor} associated to a sink $\ell \in Q_0$ 
	\begin{align}
		F_\ell ^+:\mod(\Lambda^{\imath}):=\rep(\ov{Q},\ov{I})\longrightarrow\mod(\bs_\ell \iLa):=\rep(\ov{Q'},\ov{I'}),
	\end{align}
	in \eqref{eq:F+} below.
	Using Proposition \ref{prop:modulated representation}, we shall identify the category $\rep(\ov{Q},\ov{I})$ with $\rep(\cm(Q, \btau))$, and respectively, $\rep(\ov{Q'},\ov{I'})$ with $\rep(\cm(Q',\btau))$.
	
	Without loss of generality, we assume that the sink $\ell \in \ov{\I}_\btau$. Let $L=(L_i,L_{ji})\in \rep(\cm(Q, \btau))$. Then $L_{ji}:\,_j\oH_i\otimes_{\oH_i} L_i\rightarrow L_j$ is an $\oH_i$-morphism for any $(i,j)\in\ov{\Omega}$. Denote by
	\[
	L_{\ell ,{\rm in}}:= (L_{\ell i})_i:\bigoplus_{i\in\ov{\Omega}(-,\ell )}  \,_\ell \oH_i\otimes_{\oH_i} L_i\longrightarrow L_\ell.
	\]
	
	Let $N_\ell :=\ker(L_{\ell ,{\rm in}})$. Recall $\BH_\ell$ from \eqref{dfn:Hi}, and note that $\oH_\ell=\BH_\ell$ is finite-dimensional by our assumption that $\ell$ is a sink. We have $\dim N_\ell<\infty$. By definition, there exists an exact sequence
	\begin{align}
		\label{def:reflection}
		0\longrightarrow N_\ell \longrightarrow \bigoplus_{i\in\ov{\Omega}(-,\ell )}  \,_\ell \oH_i\otimes_{\oH_i} L_i\xrightarrow{L_{\ell ,{\rm in}}} L_\ell .
	\end{align}
	Denote by
	$(N_{i\ell }^\vee)_i$ the inclusion map $N_\ell \rightarrow \bigoplus_{i\in\ov{\Omega}(-,\ell )}\,_\ell \oH_i\otimes_{\oH_i} L_i $.

	For any $L\in \rep(\cm(Q, \btau))$, define
	\begin{align}
		\label{eq:F+}
		F_\ell ^+(L)=(N_r,N_{rs})\in \rep(\cm(Q',\btau)),
	\end{align}
	where
	\[
	N_r:=\left\{ \begin{array}{cc} L_r & \text{ if }r\neq \ell ,
		\\
		N_{\ell } &\text{ if }r=\ell,
	\end{array} \right.
	\qquad
	N_{rs}:=\left\{ \begin{array}{cc} L_{rs} &\text{ if } (s,r)\in \ov{\Omega} \text{ with }r\neq \ell ,
		\\(N_{r\ell }^\vee)^\vee & \text{ if } (s,r)\in\ov{\Omega}^* \text{ and } s=\ell.
	\end{array}  \right.
	\]
	Here $(N_{r\ell }^\vee)^\vee={\rm ad}_{r\ell }^{-1}(N_{r\ell }^\vee)$; see Lemma~\ref{lem:ad}.
	
	Dually, associated to any source $\ell \in Q_0$, we have a reflection functor
	\begin{align}
		F_\ell ^-:\mod(\Lambda^{\imath}) \longrightarrow\mod(\bs_\ell \Lambda^{\imath}).
	\end{align}
	
	As proved in \cite[Proposition 2.4]{LW22b},  the pair $(F_\ell^-,F_\ell^+)$ is a pair of adjoint functors, i.e., there is a functorial isomorphism
	\begin{align*}
		\Hom_{\Lambda^\imath}(F_\ell^-(M),N)\cong \Hom_{\bs_\ell\Lambda^\imath}(M,F_\ell^+(N)). 
	\end{align*}

	Let $\ell\in Q$ be a sink. For any $A\in\{\Lambda^\imath,\bs_\ell\Lambda^\imath\}$ and $j\in Q_0$, let
	\begin{align*}
		\ct_j^A:=\{X\in\mod(A)\mid \Hom_A(X,S_j\oplus S_{\btau j})=0\},\\
		\cs_j^A:=\{X\in\mod(A)\mid \Hom_A(S_j\oplus S_{\btau j},X)=0\}.
	\end{align*}
	There is 
	an equivalence of subcategories:
	\begin{align}
		\label{eqn:YY}
		F^+_\ell: \ct_\ell^{\Lambda^\imath}\stackrel{\simeq}{\longrightarrow} \cs_\ell^{\bs_\ell\Lambda^\imath},
	\end{align}
	with its inverse given by $F^-_\ell$; see \cite[Corollary 2.5]{LW22b}.

	Let $\ct:=\ct_\ell^{\Lambda^\imath}$, and let $\cf$ be the extension closed subcategory of $\mod(\Lambda^\imath)$ generated by $S_\ell$ and $S_{\btau \ell}$.
	
	\begin{lemma}[\text{\cite[Lemma 2.7]{LW22b}}]
		\label{lem: torsion pair}
		{\quad}
		\begin{itemize}
			\item[(a)] $(\ct,\cf)$ is a torsion pair in $\mod(\Lambda^\imath)$;
			\item[(b)] For any $M\in \mod(\Lambda^\imath)$, there exists a short exact sequence
			\[
			0 \longrightarrow M \longrightarrow T_M^0 \longrightarrow T_M^1 \longrightarrow 0
			\]
			with $T_M^0,T_M^1\in \ct$ and $T_M^1\in\cp^{\leq1}(\Lambda^\imath)$.
		\end{itemize}
	\end{lemma}
	
	We have a $\Z$-linear isomorphism $\dimv\colon K_0(\mod(\bfk Q)) \rightarrow \Z^\I$, which sends an isoclass to its dimension vector. By identifying $i \in \I$ with a simple root $\alpha_i$ and thus $\Z^\I$ with the root lattice (of a Kac-Moody algebra $\fg$), we have a simple reflection $s_i$ acting on $\Z^\I$. 
	Then  we denote
	\begin{align*}
		\bs_i= \left\{
		\begin{array}{ll}
			s_{i}, & \text{ if } i=\btau i
			\\
			s_is_{\btau i}, & \text{ if } i\neq \btau i.
		\end{array}
		\right.
	\end{align*}
	
	\begin{lemma}[\text{\cite[Lemma 2.9]{LW22b}}]
		\label{lem:reflecting dimen}
		Let $(Q, \btau)$ be an $\imath$quiver with a sink $\ell$. Let $L\in\mod(\bfk Q)\subseteq \mod(\iLa)$ be an indecomposable module. Then either $F_\ell^+(L)=0$ (equivalently $L\cong S_\ell$ or $S_{\btau \ell}$) or $F_\ell^+(L)$ is indecomposable with
		$\dimv F_\ell^+(L)=\bs_\ell(\dimv L)$.
	\end{lemma}

	\begin{theorem}[\text{\cite[Theorem 4.3]{LW21a}}]
		\label{thm:Gamma}
		Let $\ell$ be a sink of $Q$. Then we have an isomorphism of algebras:
		\begin{eqnarray}
			\Gamma_{\ell}:\tMH & \stackrel{\cong}{\longrightarrow} & \tMHl,
			\label{eqn:reflection functor 2} \\
			\, [M]&\mapsto& \sqq^{\langle \res (T_M),\res(M)\rangle_Q}q^{-\langle T_M,M\rangle} [F_{\ell}^+(T_M)]^{-1}* [F_{\ell}^+(X_M)].
			\notag
		\end{eqnarray}
		where $M\in\mod(\Lambda^\imath)$ and $X_M,T_M\in \ct$ ($T_M\in\cp^{\leq1}(\Lambda^\imath)$) fit into a short exact sequence
		$0\rightarrow M\rightarrow X_M\rightarrow T_M\rightarrow0$.
	\end{theorem}
	
	Dually, the functor $F_\ell^-$ for any source $\ell$ of $Q$ induces an isomorphism of algebras:
	$$\Gamma_{\ell}^-:\tMH  \stackrel{\cong}{\longrightarrow}  \tMHl.$$
	\subsection{Actions of reflection functors on simple modules}
	
	We first have the following formulas for the isomorphism $\Gamma_{\ell}:\tMH\rightarrow \tMHl$, which basically follow from the definitions:
	\begin{align}
		\Gamma_{\ell}([M])&= [F_{\ell}^+(M)], \quad \forall M\in\ct,
		\label{eqn:reflection 1}
		\\
		\Gamma_{\ell}([S_\ell])&= \left\{
		\begin{array}{cc}\sqq[\E'_\ell]^{-1}* [S_{\btau \ell}'], &\text{ if }\btau \ell\neq \ell,
			\\ \,[\E'_\ell]^{-1}* [S_{\btau \ell}'], &\text{ if }\btau \ell= \ell, \end{array}\right.
		\label{eqn:reflection 2}
		\\
		\Gamma_{\ell}([S_{\btau \ell}])&=
		\sqq [\E'_{\btau \ell}]^{-1} * [S_{\ell}'], \quad\quad \text{ if }\btau \ell\neq \ell,
		\label{eqn:reflection 3}
		\\
		\Gamma_{\ell}([\E_\alpha])&= [\E'_{\bs_{\ell}\alpha}],
		\quad \forall\alpha\in K_0(\mod(\bfk Q)).\label{eqn:reflection 4}
	\end{align}

	Next, we compute the formulas for $\Gamma_\ell$ on simple modules $S_j$, for $j\neq i, \tau i$. 
	
	\begin{theorem}[\text{\cite[Theorem 3.2]{LW22b}}]
		\label{thm:Braidsplit}
		Let $(Q,\btau)$ be an $\imath$quiver.
		For any sink $i\in Q_0$ such that $i=\btau i \neq j$, we have
		\begin{align}
			\label{eqn:Braidsplit}
			\Gamma_i([S_j])&= \sum_{r+s=-c_{ij}}(-1)^r\sqq^{r}(1-\sqq^2)^{c_{ij}} [S'_i]_{\ov{p}}^{(r)}*[S'_j]*[S'_i]_{\ov{c_{ij}}+\ov{p}}^{(s)}\\
			\notag&+ (-1)^{p}  \sum_{t\geq1}\sum_{\stackrel{r+s+2t=-c_{ij} }{\ov{r}=\ov{p} }}   \sqq^{r}(1-\sqq^2)^{c_{ij}+2t}[S'_i]_{\ov{p}}^{(r)}*[S'_j]*[S'_i]_{\ov{c_{ij}}+\ov{p}}^{(s)}*[\E'_i]^t.
		\end{align}
	\end{theorem}
	
	The proof of Theorem \ref{thm:Braidsplit} is reduced to the consideration of the following two types of rank 2 $\imath$quivers (nevertheless, the computations are the same for both types):

	\begin{center}\setlength{\unitlength}{0.7mm}
		\begin{equation}
			\label{diag:split KM}
			\begin{picture}(50,13)(0,0)
				\put(-40,0){\begin{picture}(50,13)
						\put(-23,0){$\ov{Q}=$}
						\put(0,-3){$i$}
						\put(12,-1){\vector(-1,0){9}}
						\put(12.5,-2.2){$a$}
						\put(16,-1){\line(1,0){8}}
						
						\put(25,-3){$j$}
						\color{purple}
						\qbezier(-1,1)(-3,3)(-2,5.5)
						\qbezier(-2,5.5)(1,9)(4,5.5)
						\qbezier(4,5.5)(5,3)(3,1)
						\put(3.1,1.4){\vector(-1,-1){0.3}}
						\qbezier(24,1)(22,3)(23,5.5)
						\qbezier(23,5.5)(26,9)(29,5.5)
						\qbezier(29,5.5)(30,3)(28,1)
						\put(28.1,1.4){\vector(-1,-1){0.3}}
						\put(-1,10){$\varepsilon_i$}
						\put(24,10){$\varepsilon_j$}
				\end{picture}}
				\put(60,0){\begin{picture}(50,13)
						\put(-23,0){$\ov{\bs_iQ}=$}
						\put(0,-3){$i$}
						\put(4,-1){\line(1,0){8}}
						\put(12.5,-2.2){$a$}
						\put(16,-1){\vector(1,0){8}}
						
						\put(25,-3){$j$}
						\color{purple}
						\qbezier(-1,1)(-3,3)(-2,5.5)
						\qbezier(-2,5.5)(1,9)(4,5.5)
						\qbezier(4,5.5)(5,3)(3,1)
						\put(3.1,1.4){\vector(-1,-1){0.3}}
						\qbezier(24,1)(22,3)(23,5.5)
						\qbezier(23,5.5)(26,9)(29,5.5)
						\qbezier(29,5.5)(30,3)(28,1)
						\put(28.1,1.4){\vector(-1,-1){0.3}}
						\put(-1,10){$\varepsilon_i$}
						\put(24,10){$\varepsilon_j$}
				\end{picture}}
			\end{picture}
			\vspace{0.2cm}
		\end{equation}
	\end{center}

	\begin{center}\setlength{\unitlength}{1mm}
		\vspace{.1cm}
		\begin{equation}
			\label{diag:qsplit KM}
			\begin{picture}(50,10)(0,-10)
				\put(-20,0){
					\begin{picture}(50,10)
						\put(0,-2){$j$}
						\put(20,-2){$\btau j$}
						
						\put(-15,-10){$\ov{Q}=$}
						
						\put(7,-12){\vector(1,-2){3.5}}
						\put(5.5,-12){$a$}
						\put(6,-10){\line(-1,2){3.5}}

						\put(3,1){\line(1,0){7.5}}
						\put(10.5,0){$r$}
						\put(12,1){\vector(1,0){7.5}}

						\put(10.5,-1.5){\vector(-1,0){7.5}}
						\put(10.5,-2.5){$r$}
						\put(19.5,-1.5){\line(-1,0){7.5}}

						\put(15.5,-12){\vector(-1,-2){3.5}}
						\put(15,-12){$a$}
						\put(16.5,-10){\line(1,2){3.5}}
						
						\put(10,-22){$i$}
						\color{purple}
						\qbezier(3,2)(11,5)(19,2)
						\put(19,2){\vector(2,-1){0.2}}
						
						\qbezier(3,-2)(11,-5)(19,-2)

						\put(3,-2){\vector(-2,1){0.2}}
						\put(10,4){$^{\varepsilon_j}$}
						\put(10,-5){$_{\varepsilon_{\btau j}}$}
						\put(10,-28){$_{\varepsilon_i}$}
						\begin{picture}(50,23)(-10,19)
							\color{purple}
							\qbezier(-1,-1)(-3,-3)(-2,-5.5)
							\qbezier(-2,-5.5)(1,-9)(4,-5.5)
							\qbezier(4,-5.5)(5,-3)(3,-1)
							\put(3.1,-1.4){\vector(-1,1){0.3}}
						\end{picture}
				\end{picture}}
				\put(50,0){
					\begin{picture}(50,10)
						\put(0,-2){$j$}
						\put(20,-2){$\btau j$}
						
						\put(-15,-10){$\ov{\bs_i Q}=$}
						
						\put(10,-18){\line(-1,2){3}}
						\put(5.5,-12){$a$}
						\put(6,-10){\vector(-1,2){3.5}}

						\put(3,1){\line(1,0){7.5}}
						\put(10.5,0){$r$}
						\put(12,1){\vector(1,0){7.5}}

						\put(10.5,-1.5){\vector(-1,0){7.5}}
						\put(10.5,-2.5){$r$}
						\put(19.5,-1.5){\line(-1,0){7.5}}

						\put(12.5,-18){\line(1,2){3}}
						\put(15,-12){$a$}
						\put(16.5,-10){\vector(1,2){3.5}}
						
						\put(10,-22){$i$}
						\color{purple}
						\qbezier(3,2)(11,5)(19,2)
						\put(19,2){\vector(2,-1){0.2}}
						
						\qbezier(3,-2)(11,-5)(19,-2)

						\put(3,-2){\vector(-2,1){0.2}}
						\put(10,4){$^{\varepsilon_j}$}
						\put(10,-5){$_{\varepsilon_{\btau j}}$}
						\put(10,-28){$_{\varepsilon_i}$}
						\begin{picture}(50,23)(-10,19)
							\color{purple}
							\qbezier(-1,-1)(-3,-3)(-2,-5.5)
							\qbezier(-2,-5.5)(1,-9)(4,-5.5)
							\qbezier(4,-5.5)(5,-3)(3,-1)
							\put(3.1,-1.4){\vector(-1,1){0.3}}
						\end{picture}
				\end{picture}}
			\end{picture}
			\vspace{1.8cm}
		\end{equation}
	\end{center}
	
	\begin{theorem}[\text{\cite[Theorem 5.1]{LW22b}}]
		\label{thm:braidqsplit}
		Let $(Q,\btau)$ be an $\imath$quiver.
		For any sink $i\in Q_0$ such that $c_{i,\btau i}=0$ and $j \neq i,\tau i$, we have
		\begin{align}
			\label{eq:braidqsplit}
			\Gamma_i([S_j])
			&= \sum^{-\max(c_{ij},c_{\tau i,j})}_{u=0} \sum^{-c_{ij}-u}_{r=0} \sum_{s=0}^{-c_{\tau i,j}-u}
			(-1)^{r+s} \sqq^{-r-s+(r-s)u} (\sqq-\sqq^{-1})^{c_{ij}+c_{\tau i,j}+2u}
			\\
			&\qquad \times [S'_{i}]^{(-c_{ij}-r-u)} [S'_{\tau i}]^{(-c_{\tau i,j}-s-u)} *[S'_j]* [S'_{\tau i}]^{(s)} *[S'_{ i}]^{(r)} *[\E'_{\tau i}]^u.
			\notag
		\end{align}
	\end{theorem}
	
	The proof of Theorem \ref{thm:Braidsplit} is reduced to the consideration of the following two types of rank 2 $\imath$quivers (nevertheless, the computations are the same for both types):

	\begin{equation}
		\label{eq:figKM2}
		\xymatrix{&i   \ar@[purple][rr]<0.5ex>^{\textcolor{purple}{\varepsilon_i}}  && {\btau i} \ar@[purple][ll]<0.5ex>^{\textcolor{purple}{\varepsilon_{\btau i}}}   \\
			\ov{Q}=& \\
			&j\ar[uu]|-{a} \ar[uurr]|-{\textcolor{white}{b}}\ar@/^0.5pc/@[purple][rr]<0.5ex>^{\textcolor{purple}{\varepsilon_j}}  \ar@<0.5ex>[rr]|-{r}& &{\btau j} \ar[uull]|-{b}  \ar[uu]|-{a} \ar@<0.5ex>[ll]|-{r} \ar@/^0.5pc/@[purple][ll]<0.5ex>^{\textcolor{purple}{\varepsilon_{\btau j}}}}  \qquad\qquad\qquad \xymatrix{&i\ar[dd]|-{a} \ar@[purple][rr]<0.5ex>^{\textcolor{purple}{\varepsilon_i}} \ar[ddrr]|-{\textcolor{white}{b}}&& {\btau i} \ar@[purple][ll]<0.5ex>^{\textcolor{purple}{\varepsilon_{\btau i}}}  \ar[dd]|-{a} \ar[ddll]|-{b} \\
			\ov{\bs_i Q}=& \\
			&j \ar@/^0.5pc/@[purple][rr]<0.5ex>^{\textcolor{purple}{\varepsilon_j}}  \ar@<0.5ex>[rr]|-{r}& &{\btau j}\ar@<0.5ex>[ll]|-{r} \ar@/^0.5pc/@[purple][ll]<0.5ex>^{\textcolor{purple}{\varepsilon_{\btau j}}}}
	\end{equation}
	
	\begin{center}\setlength{\unitlength}{1mm}
		\vspace{.1cm}
		\begin{equation}
			\label{diag:qsplit KM2}
			\begin{picture}(40,10)(0,-10)
				
				\put(-20,0){
					\begin{picture}(50,10)
						\put(0,-2){$i$}
						\put(20,-2){$\btau i$}
						
						\put(-15,-10){$\ov{Q}=$}
						
						\put(10,-18){\line(-1,2){3}}
						\put(5.5,-12){$a$}
						\put(6,-10){\vector(-1,2){3.5}}

						\put(12.5,-18){\line(1,2){3}}
						\put(15,-12){$a$}
						\put(16.5,-10){\vector(1,2){3.5}}
						
						\put(10,-22){$j$}
						\color{purple}
						
						\put(2.5,1){\vector(1,0){17}}
						\put(19.5,-1){\vector(-1,0){17}}
						\put(10,1){$^{\varepsilon_i}$}
						\put(10,-3){$_{\varepsilon_{\btau i}}$}
						\put(10,-28){$_{\varepsilon_j}$}
						\begin{picture}(50,23)(-10,19)
							\color{purple}
							\qbezier(-1,-1)(-3,-3)(-2,-5.5)
							\qbezier(-2,-5.5)(1,-9)(4,-5.5)
							\qbezier(4,-5.5)(5,-3)(3,-1)
							\put(3.1,-1.4){\vector(-1,1){0.3}}
						\end{picture}
				\end{picture}}
				
				\put(50,0){
					\begin{picture}(50,10)
						\put(0,-2){$i$}
						\put(20,-2){$\btau i$}
						
						\put(-15,-10){$\ov{\bs_i Q}=$}
						\put(7,-12){\vector(1,-2){3.5}}
						\put(5.5,-12){$a$}
						\put(6,-10){\line(-1,2){3.5}}

						\put(15.5,-12){\vector(-1,-2){3.5}}
						\put(15,-12){$a$}
						\put(16.5,-10){\line(1,2){3.5}}
						\put(10,-22){$j$}

						\color{purple}
						\put(2.5,1){\vector(1,0){17}}
						\put(19.5,-1){\vector(-1,0){17}}
						\put(10,1){$^{\varepsilon_i}$}
						\put(10,-3){$_{\varepsilon_{\btau i}}$}
						\put(10,-28){$_{\varepsilon_j}$}
						\begin{picture}(50,23)(-10,19)
							\color{purple}
							\qbezier(-1,-1)(-3,-3)(-2,-5.5)
							\qbezier(-2,-5.5)(1,-9)(4,-5.5)
							\qbezier(4,-5.5)(5,-3)(3,-1)
							\put(3.1,-1.4){\vector(-1,1){0.3}}
						\end{picture}
				\end{picture}}
			\end{picture}
			\vspace{1.8cm}
		\end{equation}
	\end{center}
	
	\subsection{Relative braid group symmetries of $\imath$quantum groups}
	
	Let $\btau$ be an involution of $Q$, which induces an involution on $\fg$ again denoted by $\btau$. We shall define the {\em relative (or restricted) Weyl group} associated to the quasi-split symmetric pair $(\fg, \fg^{\omega\tau})$ to be the following subgroup $W^{\btau}$ of $W$:
	\begin{align}
		\label{eq:Wtau}
		W^{\btau} =\{w\in W\mid \btau w =w \btau\}
	\end{align}
	where $\btau$ is regarded as an automorphism of $\Aut(C)$. In finite type, it is well known that the relative Weyl group defined in this way coincides with the one arising from real groups (cf., e.g., \cite{Lus03, KP11}).
	
	Recall the subset $\I_\tau$ of $\I$ from \eqref{eq:ci}, and define
	\begin{align}
		\label{eq:Itau2}
		\ov{\I}_\btau:=\{i\in\I_\btau\mid c_{i,\btau i}=0 \text{ or }2 \}.
	\end{align}
	In our setting, $\ov{\I}_\btau$ consists of exactly those $i \in \I_\btau$ such that the $\btau$-orbit of $i$ is of finite type. Note that $\ov{\I}_\btau=\I_\btau$ if $(Q,\btau)$ is acyclic. We denote by $\bs_{i}$, for $i\in\ov{\I}_\btau$, the following element of order 2 in the Weyl group $W$
	\begin{align}
		\label{def:si}
		\bs_i= \left\{
		\begin{array}{ll}
			s_{i}, & \text{ if } i=\btau i
			\\
			s_is_{\btau i}, & \text{ if } i\neq \btau i.
		\end{array}
		\right.
	\end{align}

	\begin{lemma}
		[\text{\cite[Appendix]{Lus03}}]
		\label{lem:iWeyl}
		The relative Weyl group $W^{\btau}$ can be identified with a Coxeter group with $\bs_i$ ($i\in \ov{\I}_\btau$) as its generators.
	\end{lemma}

	Let $i$ be a sink of an $\imath$quiver $(Q,\btau)$. By Theorem \ref{thm:main}, there are injective homomorphisms: $\widetilde{\Psi}_{Q}: \tUi_{\sqq}\rightarrow \tMH$,  $\widetilde{\Psi}_{Q'}: \tUi_{\sqq}\rightarrow \tMHi$.
	
	Recall $\ov{\I}_\btau$ from \eqref{eq:Itau2}. We define algebra automorphisms
	\[
	\TT''_{i,1} \in \Aut (\tUi), \qquad
	\text{ for } i \in \ov{\I}_\btau,
	\]
	such that
	specializing at $v=\sqq$ we have the following commutative diagram:
	\begin{equation}
		\label{eq:Up=T}
		\xymatrix{
			\tUi _{\sqq}\ar[r]^{\TT''_{i,1}}  \ar[d]^{\widetilde{\Psi}_Q} & \tUi_{\sqq}\ar[d]^{\widetilde{\Psi}_{Q'}}
			\\
			\tMH \ar[r]^{\Gamma_i} &  \tMHi
		}
	\end{equation}
	
	Using Theorems \ref{thm:Braidsplit}--\ref{thm:braidqsplit}, we have the following theorem. 
	
	\begin{theorem}[\text{\cite[Theorem 6.8]{LW22b}}]
		\label{thm:BG}
		We have a $\Q(v)$-algebra automorphism $\TT''_{i,1}$ of $\tUi$, for $i\in\bar{\I}_\btau$, such that
		\begin{enumerate}
			\item
			$\underline{(i=\btau i)}:$ \;
			$\TT''_{i,1}(\tk_j)= (-v^2 \tk_i)^{-c_{ij}}\tk_j$, and
			\begin{align*}
				\TT''_{i,1}(B_i) &=(-v^{2} \tk_{i} )^{-1}B_i,
				\\
				\TT''_{i,1}(B_j) &= \sum_{r+s=-c_{ij}}(-1)^{r}v^{r}B_{i,\ov{p}}^{(r)}B_jB_{i,\ov{c_{ij}}+\ov{p}}^{(s)}\\
				& +\sum_{u\geq1}\sum_{\stackrel{r+s+2u=-c_{ij}}{
						\ov{r}=\ov{p}}}(-1)^{r}v^{r}B_{i,\ov{p}}^{(r)}B_jB_{i,\ov{c_{ij}}+\ov{p}}^{(s)}(-v^2 \tk_i )^u, \qquad \text{ for }j\neq i;
			\end{align*}
			\item
			$\underline{(i\neq \btau i)}:$ \;
			$\TT''_{i,1}(\tk_j)= \tk_i^{-c_{ij}} \tk_{\btau i}^{-c_{\tau i,j}} \tk_j$,
			\begin{align*}
				\TT''_{i,1}(B_j) &=
				\begin{cases}  -\tk_{i}^{-1}B_{\btau i},  & \text{ if }j=i \\
					-B_i\tk_{\btau i}^{-1}  ,  &\text{ if }j=\btau i,\end{cases}
			\end{align*}
			and for $j\neq i,\btau i$,
			\begin{align*}
				\TT''_{i,1}(B_j)
				&= \sum^{-\max(c_{ij},c_{\tau i,j})}_{u=0} \; \sum^{-c_{ i,j}-u}_{r=0} \; \sum_{s=0}^{-c_{\tau i,j}-u} (-1)^{r+s} v^{ r-s+(-c_{ij}-r-s-u)u } \\
				&\qquad\qquad\qquad\qquad\qquad
				\times B_i^{(r)} B_{\tau i}^{(-c_{\tau i,j}-u-s)} B_j B_{\tau i}^{(s)} B_i^{(-c_{ij}-r-u)}\tk_{\tau i}^u.
			\end{align*}
		\end{enumerate}
	\end{theorem}
	
	The following lemma follows by inspection of the defining relations for $\tUi$ in Theorem~\ref{thm:Serre};  cf. \cite{CLW21c}.
	\begin{lemma}
		\label{lem:bar}
		(a) 
		There exists a $\Q$-algebra involution $ \psi_\imath: \tUi\rightarrow \tUi$ (called a bar involution) such that
		\[
		\psi_\imath(v)=v^{-1}, \quad
		\psi_\imath(\tk_i)=v^{c_{i,\btau i}}\tk_{\btau i}, \quad
		\psi_\imath (B_i)=B_i,  \quad
		\forall i\in \I.
		\]
		
		(b)
		There exists a $\Q(v)$-algebra anti-involution $\sigma_\imath: \tUi\rightarrow \tUi$ such that
		\begin{align}
			\label{eq:sigma}
			\sigma_\imath(B_i)=B_{i}, \quad \sigma_\imath(\tk_i)= \tk_{\btau i},
			\quad \forall i\in \I.
		\end{align}
	\end{lemma}
	
	\begin{theorem}[\text{\cite[Theorems 6.11--6.12]{LW22b}}]
		\label{thm:BG1}
		For $i\in \bar{\I}_\btau$ and $e \in \{\pm 1\}$,
		there are automorphisms $\TT'_{i,e},\TT''_{i,e}$ on $\tUi$  such that
		\begin{align}
			\label{eq:TTsig}
			\TT_{i,e}'  &= \sigma_\imath \TT_{i,-e}'' \sigma_\imath,
			\qquad
			\psi_\imath \TT_{i,e}' \psi_\imath =\TT_{i,-e}'.
		\end{align}
		Moreover, we have $\TT_{i,e}' =(\TT_{i,-e}'')^{-1}$, for any $i\in\bar{\I}_\btau$ and $e=\{\pm1\}$.
	\end{theorem}
	
	\begin{proof}
		We only need to prove the last statement. 
		
		It suffices to prove for $e=-1$ by \eqref{eq:TTsig}.
		Given an $\imath$quiver $(Q,\btau)$ such that  $i$ is a sink, we have two isomorphisms, $\Gamma_i:\tMH  \stackrel{\cong}{\rightarrow}  \tMHi$ and $\Gamma_i^-:\tMHi\stackrel{\cong}{\rightarrow} \tMH$, which are inverses to each other.
		In particular, we have the following commutative diagram:
		\begin{align}
			\label{eq:defT2}
			\xymatrix{
				\tUi_{\sqq}  \ar[r]^{\TT'_{i,-1}}  \ar[d]^{\widetilde{\Psi}_{Q'}} & \tUi_{\sqq} \ar[d]^{\widetilde{\Psi}_{Q}} \\
				\ar[r]^{\Gamma_i^-}  \tMHi & \tMH
			}
		\end{align}
		We conclude that $\TT_{i,-1}'$ and $\TT_{i,1}''$ are inverses to each other.
	\end{proof}
	
	\begin{remark}
		The results in Theorem~\ref{thm:BG} and Theorem \ref{thm:BG1} verify part of \cite[Conjecture~ 6.5]{CLW21b} in case $i=\tau i$ and  \cite[Conjecture~3.7]{CLW21c} in case $i\not =\tau i$, for quasi-split $\imath$quantum groups $\tUi$ associated with symmetric generalized Cartan matrices such that all $c_{i,\tau i}$ even. (Note that in Theorem \ref{thm:BG}(1), our notation here $\TT_{i,e}', \TT_{i,e}''$ follow \cite[Theorem {\bf A}(1)]{LW22b}, which was swapped in \cite[Conjecture~ 6.5]{CLW21b}.) 
		
		A complete proof of \cite[Conjecture~ 6.5]{CLW21b} and  \cite[Conjecture~3.7 \& \S3.5]{CLW21c} is obtained by Zhang \cite{Z22}, who has further developed the approach initiated in \cite{WZ22}. In addition, Zhang shows that $\TT'_{i,e}$ (or $\TT''_{i,e}$) satisfy the relative braid group relations. 
	\end{remark}

	\section{$\imath$Hall algebras of weighted projective lines and  $\imath$quantum loop algebras}
	\label{sec:iWPL}
	
	In this section, we outline the connection between $\imath$Hall algebras of (weighted) projective lines and affine $\imath$quantum groups in Drinfeld type presentations.

	\subsection{Drinfeld type presentation of $q$-Onsager algebra}
	
	The (universal) $\imath$quantum group $\tUi(\hA)$ associated to the Satake diagram of split affine type $A_1^{(1)}$ is known as {\em (universal) $q$-Onsager algebra}. By definition, $\tUi(\hA)$ is generated by $B_i, \tk_i$, for $i=0, 1$, subject to relations \eqref{relation1}--\eqref{relation6}; in $\imath$Hall algebra setting, the renormalized generators $\K_i:= -v^2 \tk_i$ are more natural as they correspond to (generalized) simple modules. 
	
	\begin{definition} [\cite{LW21c}]
		Let $\tUiD(\hA)$ be the $\Q(v)$-algebra  generated by $\K_1^{\pm1}$, $C^{\pm1}$, $H_{m}$ and $B_{1,r}$, where $m\geq1$, $r\in\Z$, subject to the following relations, for $r,s\in \Z$ and $m,n\ge 1$:
		\begin{align}
			\K_1\K_1^{-1}=1, C C^{-1}& =1, \quad \K_1, C \text{ are central, }
			\\
			[H_m,H_n] &=0,  \label{iDR1}
			\\
			[H_m, B_{1,r}] &=\frac{[2m]}{m} B_{1,r+m}-\frac{[2m]}{m} B_{1,r-m}C^m,
			\label{iDR2}
			\\
			\label{iDR3}
			[B_{1,r}, B_{1,s+1}]_{v^{-2}}  -v^{-2} [B_{1,r+1}, B_{1,s}]_{v^{2}}
			&= v^{-2}\Theta_{s-r+1} C^r \K_1-v^{-4} \Theta_{s-r-1} C^{r+1} \K_1 \\\notag
			&\quad +v^{-2}\Theta_{r-s+1} C^s \K_1-v^{-4} \Theta_{r-s-1} C^{s+1} \K_1.\notag
		\end{align}
		Here
		\begin{align}
			\label{eq:exp1}
			1+ \sum_{m\geq 1} (v-v^{-1})\Theta_{m} z^m  = \exp\Big( (v-v^{-1}) \sum_{m\ge 1}  H_m z^m \Big).
		\end{align}
	\end{definition}
	
	\begin{theorem}   [\text{\cite[Theorem 2.16]{LW21c}}] 
		\label{thm:Dr1}
		There is an isomorphism of $\Q(v)$-algebras 
		\begin{align*}
			\Phi: \tUiD (\hA) \longrightarrow\tUi (\hA). 
		\end{align*} 
	\end{theorem}
	We shall call $\tUiD(\hA)$ the Drinfeld type  
	presentation of the $q$-Onsager algebra. To establish Theorem~\ref{thm:Dr1}, one needs to construct the real root vectors $B_{1,r}$ and imaginary root vectors $\Theta_{m}$ in $\tUi(\hA)$ via braid group symmetries (cf. \cite{BK20, LW21c}), and the map $\Phi$ matches the generators in the same notations and identifes $C=\K_0\K_1$. To show $\Phi$ is an algebra homomorphism amounts to showing that all the relations in $\tUiD(\hA)$ are satisfied by these root vectors in $\tUi$. The surjectivity of $\Phi$ is easy, while the injectivity follows by some filtration argument.

	\subsection{$\imath$Hall algebra of the projective line} 
	
	Recall $\bfk$ denotes a finite field of $q$ elements. 
	The coordinate ring of the projective line $\PL$ over $\bfk$ is the $\bbZ$-graded ring $\bS=\bfk[X_0,X_1]$ with $\deg(X_0)=\deg(X_1)=1$. A closed point $x$ of $\PL$ is given by a prime homogeneous ideal of $\bS$ generated by an irreducible polynomial in $\bS$. The degree of $x$, denoted by $\deg(x)$ or $d_x$, is defined to be the degree of the defining irreducible polynomial associated to $x$.
	We denote
	
	$\triangleright$
	$\mod^{\bbZ}(\bS)$ -- category of finitely generated $\bbZ$-graded $\bS$-modules
	
	$\triangleright$
	$\mod_0^{\bbZ}(\bS)$-- the full subcategory of $\mod^{\bbZ}(\bS)$ of finite-dimensional graded $\bS$-modules
	
	$\triangleright$
	$\coh(\PL)$ -- category of cohorent sheaves on $\PL$

	We can associate a coherent sheaf $\widetilde{M}$ on $\PL$ to any $M \in \mod^{\bbZ}(\bS)$. This gives rise to a category equivalence    (which goes back to Serre):
	$\mod^{\mathbb{Z}}(\bS) \big / \mbox{mod}_{0}^{\mathbb{Z}}(\bS)
	\cong \coh(\PL).$
	The category $\coh(\PL)$ is a finitary hereditary abelian Krull-Schmidt category.
	
	Any indecomposable vector bundle on $\PL$ is a line bundle; more precisely \cite{Gro57}, it is of the form $\co(n)=\widetilde{\bS[n]}$, for $n\in\bbZ$, where $\bS[n]$ is the $n$-th shift of the trivial module $\bS$, i.e., $\bS[n]_i=\bS_{n+i}$. The homomorphism between two line bundles are given by
	\begin{align}   \label{eq:OmOn}
		\Hom(\co(m),\co(n))\cong \bS_{n-m},
	\end{align}
	which then has dimension $n-m+1$ if $n\geq m$ and $0$ otherwise.
	
	For a hereditary abelian category $\ca$, let $\cc_{\Z_1}(\ca)$ be the category of $1$-periodic complexes of $\ca$. Recall that $\cs\cd\widetilde{\ch}(\cc_{\Z_1}(\ca))$ is the twisted semi-derived Ringel-Hall algebra of $\cc_{\Z_1}(\ca)$. 
	
	For any $X\in\ca$, denote the stalk complex in $\cc_{\Z_1}(\ca)$ by
	\[
	C_X =(X,0)
	\]
	(or just by $X$ when there is no confusion), and denote by $K_X$ the following acylic complex:
	\[
	K_X:=(X\oplus X, d),
	\qquad \text{ where }
	d=\left(\begin{array}{cc} 0&\Id \\ 0&0\end{array}\right).
	\]
	For any $\alpha\in K_0(\ca)$,  there exist $X,Y\in\ca$ such that $\alpha=\widehat{X}-\widehat{Y}$. Define $[K_\alpha]:=[K_X]* [K_Y]^{-1}\in\cs\cd\widetilde{\ch}(\cc_{\Z_1}(\ca))$.

	Return to $\ca=\coh(\PL)$. We shall use a shorthand notation $\tMHX$ to denote the twisted semi-derived Ringel-Hall algebra of  $\cc_{\Z_1}(\coh(\PL))$.

	For $m\ge 1$, we define
	\begin{align}  
		\label{eq:hTm}
		\haT_{m}= \frac{1}{(q-1)^2\sqq^{m-1}}\sum_{0\neq f:\co\rightarrow \co(m) } [\coker f].
	\end{align}
	
	The classes $\widehat{\co}$ and $\widehat{\co(1)}$ form a basis of $K_0(\PL):= K_0(\coh(\PL))$. Denote by
	\begin{align}
		\delta:=\widehat{\co(1)}-\widehat{\co}.
		\label{eq:delta}
	\end{align}
	Then $\{\widehat{\co}, \delta\}$ is also a basis.
	
	\begin{theorem}  [\text{\cite[Theorem 4.2]{LRW20}}]
		\label{thm:morphi}
		There exists a $\Q(\sqq)$-algebra homomorphism
		\begin{align}
			\label{eq:phi}
			\Omega: {}^{\text{Dr}}\tUi_{\sqq}(\widehat{\mathfrak{sl}}_2) \longrightarrow \tMHX
		\end{align}
		which sends, for all $r\in \Z$ and $m \ge 1$,
		\begin{align*}
			\K_1\mapsto [K_\co],  \quad
			C\mapsto [K_\delta],  \quad
			B_{1,r} \mapsto -\frac{1}{q-1}[\co(r)],  \quad
			\Theta_{m}\mapsto  \haT_m.
		\end{align*}
	\end{theorem}
	Moreover, $\widehat{H}_m:=\Omega(H_m)$ is also described in \cite{LRW20}.

	\subsection{$\imath$Hall algebra interpretation of the isomorphism $\Phi$}

	Consider the Kronecker quiver $\QK: \xymatrix{0\ar@<0.5ex>[r]^\alpha \ar@<-0.5ex>[r]_\beta& 1}$. Let $\LaK^\imath$ be the $\imath$quiver algebra of the split Kronecker $\imath$quiver $(\QK,\Id)$ described as in Example \ref{example 2}(e). 
	We can make sense of the $\imath$Hall algebra $\tMHL=\widetilde{\ch}(
	\bfk \QK,\Id)$. 
	Let $T=\co\oplus\co(1)$ and $B=\End_{\PL}(T)$. It is known  that $T$ is a tilting object, and $B^{op}\cong \bfk \QK$. It follows that
	\begin{align}  \label{Bei}
		\RHom_{\PL}(T,-):\cd^b(\coh(\PL))\stackrel{\simeq}{\longrightarrow} \cd^b(\rep_\bfk( \QK))
	\end{align}
	is a derived equivalence. 
	Let $\cv$ be the subcategory of $\coh(\PL)$ consisting of $M$ such that $\Hom_{\PL}(T,M)=0$. Denote $\cu=\Fac T$, the full subcategory of $\coh(\PL)$ consisting of homomorphic images of objects in $\add T$. Then $(\cu,\cv)$ is a {\em torsion pair} of $\coh(\PL)$. It is proved in \cite[Lemma 5.8]{LRW20} that  $(\cc_1(\cu),\cc_1(\cv))$ is a torsion pair of $\cc_1(\coh(\P_\bfk^1))$. In particular, any $M\in\cc_1(\coh(\P_\bfk^1))$ admits a short exact sequence of the form
	\begin{align}
		\label{T-resol}
		0 \longrightarrow M \longrightarrow X_M \longrightarrow T_M \longrightarrow 0
	\end{align}
	where $X_M\in \cc_1(\cu)$ and $T_M\in\add K_T$. 
	Then applying \cite[Theorem A.22]{Lu22} to our setting gives us an algebra isomorphism:
	\begin{align*}
		\BF:\tMHX&\stackrel{\sim}{\longrightarrow} \tMHL\\\notag
		[M]&\mapsto [F(T_M)]^{-1}* [F(X_M)],
	\end{align*}
	where $X_M\in \cc_1(\cu)$ and $T_M\in\add K_T$, are defined in the short exact sequence \eqref{T-resol}.
	
	\begin{theorem}  [\text{\cite[Theorem 5.11]{LRW20}}]
		\label{main thm2}
		We have the following commutative diagram of algebra homomorphisms
		\[
		\xymatrix{\tUi_{\sqq}(\widehat{\mathfrak{sl}}_2) \ar[d]^{\widetilde{\psi}} \ar[r]^{\Phi^{-1}}   & {}^{\text{Dr}}\tUi_{\sqq}(\widehat{\mathfrak{sl}}_2) \ar[d]^{\Omega}
			\\
			\tMHL \ar[r]^{\BF^{-1}} & \tMHX  }\]
		where $\Phi, \BF$ are isomorphisms.
		In particular, the homomorphism $\Omega$ is injective.
	\end{theorem}

	\subsection{$\imath$Quantum loop algebras}
	
	Associated to a generalized Cartan matrix (GCM) $C=(c_{ij})_{i,j\in \II}$, the $\imath$quantum loop algebra $\tUiD$ of split type is the $\Q(v)$-algebra  generated by $\K_{i}^{\pm1}$, $C^{\pm1}$, $H_{i,m}$, $\Theta_{i,l}$ and $\y_{i,l}$, where  $i\in \II$, $m \in \Z_{+}$, $l\in\Z$, subject to some relations. 
	
	To describe the relations in $\tUiD$,  we introduce some shorthand notations below. 
	Let $k_1, k_2, l\in \Z$ and $i,j \in \II$. Set
	\begin{align}
		\begin{split}
			S(k_1,k_2|l;i,j)
			&=  B_{i,k_1} B_{i,k_2} B_{j,l} -[2] B_{i,k_1} B_{j,l} B_{i,k_2} + B_{j,l} B_{i,k_1} B_{i,k_2},
			\\
			\SS(k_1,k_2|l;i,j)
			&= S(k_1,k_2|l;i,j)  + \{k_1 \leftrightarrow k_2 \}.
			\label{eq:Skk}
		\end{split}
	\end{align}
	Here and below, $\{k_1 \leftrightarrow k_2 \}$ stands for repeating the previous summand with $k_1, k_2$ switched if $k_1\neq k_2$, so the sums over $k_1, k_2$ are symmetric.
	We also denote
	\begin{align}
		\begin{split}
			R(k_1,k_2|l; i,j)
			&=   \K_i  C^{k_1}
			\Big(-\sum_{p\geq0} v^{2p}  [2] [\Theta _{i,k_2-k_1-2p-1},\y_{j,l-1}]_{v^{-2}}C^{p+1}
			\label{eq:Rkk} \\
			&\qquad\qquad -\sum_{p\geq 1} v^{2p-1}  [2] [\y_{j,l},\Theta _{i,k_2-k_1-2p}]_{v^{-2}} C^{p}
			- [\y_{j,l}, \Theta _{i,k_2-k_1}]_{v^{-2}} \Big),
			\\
			\R(k_1,k_2|l; i,j) &= R(k_1,k_2|l;i,j) + \{k_1 \leftrightarrow k_2\}.
		\end{split}
	\end{align}

	
	Now we can give the definition of $\imath$quantum loop algebras $\tUiD$ for simply-laced generalized Cartan matrices $C$, following \cite[Remark 3.17]{LW21c}. If $C$ is a Cartan matrix of simple Lie algebra $\fg$ of type ADE, then $\tUiD$ provides a Drinfeld type presentation of the affine $\imath$quantum group of split type (associated to the affinization of $\fg$). That is, there is a $\Q(v)$-algebra isomorphism ${\Phi}: \tUiD \rightarrow\tUi$; see \cite{LW21c}.
	
	\begin{definition}[$\imath$quantum loop algebras]
		\label{def:iDR}
		Let $C=(c_{ij})_{i,j\in \II}$ be a simply-laced generalized Cartan matrix (GCM). 
		The $\imath$quantum loop algebra $\tUiD$ of split type is the $\Q(v)$-algebra  generated by $\K_{i}^{\pm1}$, $C^{\pm1}$, $H_{i,m}$ and $\y_{i,l}$, where  $i\in \II$, $m \in \Z_{+}$, $l\in\Z$, subject to the following relations, for $m,n \in \Z_{+}$ and $k,l\in \Z$:
		\begin{align}
			& \K_i, C \text{ are central, } \quad \K_i\K_i^{-1}=1, \;\; C C^{-1}=1,
			\label{iDR1a}
			\\
			&[H_{i,m},H_{j,n}]=0,\label{iDR1b}\\
			&[H_{i,m},\y_{j,l}]=\frac{[mc_{ij}]}{m} \y_{j,l+m}-\frac{[mc_{ij}]}{m} \y_{j,l-m}C^m,
			\label{iDR2}
			\\
			&[\y_{i,k} ,\y_{j,l}]=0,   \text{ if }c_{ij}=0,  \label{iDR4}
			\\
			&[\y_{i,k}, \y_{j,l+1}]_{v^{-c_{ij}}}  -v^{-c_{ij}} [\y_{i,k+1}, \y_{j,l}]_{v^{c_{ij}}}=0, \text{ if }i\neq j;
			\label{iDR3a}
			\\ 
			&[\y_{i,k}, \y_{i,l+1}]_{v^{-2}}  -v^{-2} [\y_{i,k+1}, \y_{i,l}]_{v^{2}}
			=v^{-2}\Theta_{i,l-k+1} C^k \K_i-v^{-4}\Theta_{i,l-k-1} C^{k+1} \K_i
			\label{iDR3b} \\
			&\qquad\qquad\qquad\qquad\qquad\qquad\quad\quad\quad
			+v^{-2}\Theta_{i,k-l+1} C^l \K_i-v^{-4}\Theta_{i,k-l-1} C^{l+1} \K_i, \notag
			\\
			\label{iDR5}
			&    \SS(k_1,k_2|l; i,j) =  \R(k_1,k_2|l; i,j), \text{ if }c_{ij}=-1.
		\end{align}
		Here we set
		\begin{align}  \label{Hm0}
			{\Theta}_{i,0} =(v-v^{-1})^{-1}, \qquad {\Theta}_{i,m} =0, \; \text{ for }m<0;
		\end{align}
		and $\Theta_{i,m}$ ($m\geq1$)  are related to  $H_{i,m}$ by the following equation:
		\begin{align}
			\label{exp h}
			1+ \sum_{m\geq 1} (v-v^{-1})\Theta_{i,m} u^m  = \exp\Big( (v-v^{-1}) \sum_{m\geq 1} H_{i,m} u^m \Big).
		\end{align}
	\end{definition}

	\subsection{Weighted projective lines and $\imath$quantum loop algebras}

	Fix a positive integer $\bt$ such that $2\leq \bt\leq q$.
	Let $\bp=(p_1,p_2,\dots,p_\bt)\in\Z_+^{\bt}$, and let $L(\bp)$
	denote the rank one abelian group on generators $\vec{x}_1$, $\vec{x}_2$, $\ldots$, $\vec{x}_{\bt}$ with relations
	$p_1\vec{x}_1=p_2\vec{x}_2=\cdots=p_{\bt}\vec{x}_{\bt}$. We call $\vec{c}=p_i\vec{x}_i$ the \emph{canonical element} of $L(\bp)$.
	Obviously, the polynomial ring $\bfk[X_1,\dots,X_{\bt}]$ is a $L(\bp)$-graded algebra by setting $\deg X_i=\vec{x}_i$, which is denoted by $\bS(\bp)$.
	
	Let $\ul{\bla}=\{\bla_1,\dots,\bla_{\mathbf{t}}\}$ be a collection of distinguished closed points (of degree one) on the projective line $\PL$, normalized such that $\bla_1=\infty$, $\bla_2=0$, $\bla_3=1$.
	Let $I(\bp,\ul{\bla})$ be the $L(\bp)$-graded ideal of $\bS(\bp)$ generated by
	$$\{X_i^{p_i}-(X_2^{p_2}-\bla_iX_1^{p_1})\mid 3\leq i\leq \bt \}.$$
	Then $\bS(\bp,\ul{\bla}):=\bS(\bp)/I(\bp,\ul{\bla})$ is an $L(\bp)$-graded algebra.
	For $1\leq i\leq \bt$, denote by $x_i$ the image of $X_i$ in $\bS(\bp,\ul{\bla})$.
	
	The \emph{weighted projective line} $\X:=\X_{\bp,\ul{\bla}}$ is the set of all non-maximal prime homogeneous ideals of $\bS:=\bS(\bp,\ul{\bla})$, which is also denoted by $\X_\bfk$ to emphasis the base field $\bfk$. 	The classification of the closed points in $\X$ is provided in  \cite[Proposition 1.3]{GL87}. First, each $\bla_i$ corresponds to the prime ideal generated by $x_i$, called the \emph{exceptional point}. Second, any other ideals are of the form $(f(x_1^{p_1},x_2^{p_2}))$, where $f\in \bfk [y_1,y_2]$ is an irreducible homogeneous polynomial in $y_1,y_2$, which is different from $y_1$ and $y_2$, called the \emph{ordinary point}. 
	
	Following \cite{GL87},	the category of coherent sheaves	$\coh(\X)$ over $\X$ is defined to be the Serre quotient $\mod^{L(\bp)}(\bS)/ \mod^{L(\bp)}_0(\bS)$ of the category of finitely generated $L(\bp)$-graded
	$\bS$-modules modulo the category of finite length graded $\bS$-modules.
	Let $\co=\co_\X$ be the structure sheaf, i.e., the image of $\bS$. A sheaf is called a \emph{torsion sheaf} if it is a finite-length object in $\coh(\X)$.
	Let $\scrt$ be the full subcategory consisting of all torsion sheaves. 
	
	In order to describe the category $\scrt$, we shall introduce the representation theory of cyclic quivers.
	We consider the oriented cyclic quiver $\text{Cy}_n$ for $n\geq2$ with its vertex set $\Z_n=\{0,1,2,\dots,n-2,n-1\}$:
	\begin{center}\setlength{\unitlength}{0.5mm}
		\begin{equation}
			\label{fig:Cn}
			\begin{picture}(100,30)
				\put(2,0){\circle*{2}}
				\put(22,0){\circle*{2}}
				
				\put(82,0){\circle*{2}}
				\put(102,0){\circle*{2}}
				\put(52,25){\circle*{2}}
				
				\put(20,0){\vector(-1,0){16}}
				\put(40,0){\vector(-1,0){16}}
				\put(47.5,-2){$\cdots$}
				\put(80,0){\vector(-1,0){16}}
				\put(100,0){\vector(-1,0){16}}
				
				\put(54,24.5){\vector(2,-1){47}}
				\put(3,1){\vector(2,1){47}}
				
				\put(50.5,27){\tiny $0$}
				\put(1,-6){\tiny $1$}
				\put(21,-6){\tiny $2$}
				\put(75,-6){\tiny $n-2$}
				\put(95,-6){\tiny $n-1$}
			\end{picture}
		\end{equation}
		\vspace{-0.2cm}
	\end{center}
	In particular,
	$\text{Cy}_1$ is just the Jordan quiver, i.e., the quiver with only one vertex and one loop arrow.
	Denote by $\rep_\bfk(\text{Cy}_n)$ the category of finite-dimensional nilpotent representations of $\text{Cy}_n$ over the field $\bfk$. Then the structure of $\scrt$ is described in the following.
	\begin{lemma}[\cite{GL87}]
		\label{lem:isoclasses Tor}
		(1) The category $\scrt$ decomposes as a coproduct $\scrt=\coprod_{x\in\X} \scrt_{x}$, where $\scrt_{x}$ is the subcategory of torsion sheaves with support at $x$.
		
		(2) For any ordinary point $x$ of degree $d$, let $\bfk_{x}$ denote the residue field at $x$, i.e., $[\bfk_x:\bfk]=d$. Then $\scrt_{x}$ is equivalent to the category $\rep_{\bfk_{x}}(C_1)$.
		
		(3) For any exceptional point $\bla_i$ ($1\leq i\leq \bt$), the category $\scrt_{\bla_i}$ is equivalent to $\rep_{\bfk}(\text{Cy}_{p_i})$.
	\end{lemma}
	
	For any ordinary point $\blx$ of degree $d$, let $\pi_{\blx}$ be the prime homogeneous polynomial corresponding to $\blx$. The multiplication by $\pi_{\blx}$ gives the exact sequence
	$$0\longrightarrow \co\stackrel{\pi_{\blx}}{\longrightarrow} \co(d\vec{c})\longrightarrow S_{\blx}\rightarrow0,$$
	where $S_{\blx}$ is the unique (up to isomorphism) simple sheaf in the category $\scrt_{\blx}$. 
	For any exceptional point $\bla_i$, multiplication by $x_i$ yields the short exact sequence
	$$0\longrightarrow \co((j-1)\vec{x_i})\stackrel{x_i}{\longrightarrow} \co(j\vec{x_i})\longrightarrow S_{ij}\rightarrow0,\text{ for } 1\le j\le p_i;$$
	where $\{S_{ij}\mid  j\in\Z_{p_i}\}$ is a complete set of pairwise non-isomorphic simple sheaves in the category $\scrt_{\bla_i}$ for any $1\leq i\leq \bt$.
	
	Let $\cc_{\Z_1}(\coh(\X_\bfk))$ be the category of $1$-periodic complexes of $\coh(\X_\bfk)$. We shall use a shorthand notation $\tMHW$ to denote the twisted semi-derived Ringel-Hall algebra of  $\cc_{\Z_1}(\coh(\X_\bfk))$.
	
	For $\bp=(p_1,\dots,p_\bt)\in\Z_{+}^\bt$, let us consider the following star-shaped graph $\Gamma=\T_{p_1,\dots,p_\bt}$:
	
	\begin{center}\setlength{\unitlength}{0.8mm}
		\begin{equation}
			\label{star-shaped}
			\begin{picture}(110,30)(0,35)
				\put(0,40){\circle*{1.4}}
				\put(2,42){\line(1,1){16}}
				\put(20,60){\circle*{1.4}}
				\put(23,60){\line(1,0){13}}
				\put(40,60){\circle*{1.4}}
				\put(43,60){\line(1,0){13}}
				\put(60,58.5){\large$\cdots$}
				\put(70,60){\line(1,0){13}}
				\put(88,60){\circle*{1.4}}

				\put(3,41){\line(4,1){13}}
				\put(20,45){\circle*{1.4}}
				\put(23,45){\line(1,0){13}}
				\put(40,45){\circle*{1.4}}
				\put(43,45){\line(1,0){13}}
				\put(60,43.5){\large$\cdots$}
				\put(70,45){\line(1,0){13}}
				\put(88,45){\circle*{1.4}}
				
				\put(19,30){\Large$\vdots$}
				
				\put(39,30){\Large$\vdots$}
				
				\put(87,30){\Large$\vdots$}
				
				\put(2,38){\line(1,-1){16}}
				\put(20,20){\circle*{1.4}}
				\put(23,20){\line(1,0){13}}
				\put(40,20){\circle*{1.4}}
				\put(43,20){\line(1,0){13}}
				\put(60,18.5){\large$\cdots$}
				\put(70,20){\line(1,0){13}}
				\put(88,20){\circle*{1.4}}
				
				\put(-4.5,39){$\star$}
				
				\put(16,62){\tiny$[1,1]$}
				
				\put(36,62){\tiny$[1,2]$}
				\put(81,62){\tiny$[1,p_1-1]$}

				\put(16.5,47){\tiny$[2,1]$}
				\put(36.5,47){\tiny$[2,2]$}
				\put(81,47){\tiny$[2,p_2-1]$}

				\put(16.5,16){\tiny$[\bt,1]$}
				\put(36.5,16){\tiny$[\bt,2]$}
				\put(81,16){\tiny$[\bt,p_\bt-1]$}

				
			\end{picture}
		\end{equation}
		\vspace{1cm}
	\end{center}
	The set of vertices is denoted by $\II$. Let $C=(c_{ij})_{i,j\in \II}$ be the generalized Cartan matrix (GCM) of $\Gamma$. Then by Definition~\ref{def:iDR}, we have an $\imath$quantum loop algebra $\tUiD$ of split type associated to $C$.

	Recall the cyclic quiver $\text{Cy}_n$ from \eqref{fig:Cn}. We can make sense of the $\imath$Hall algebra $\iH(\bfk\text{Cy}_n)=\widetilde{\ch}(
	\bfk \text{Cy}_n,\Id)$.
	As a special case of \cite[Theorem 9.6]{LW20}, we have an algebra monomorphism:
	\begin{align}
		\label{eq:isoCnsln}
		\widetilde{\psi}_{\text{Cy}_n}: \tUi_{\sqq}(\widehat{\mathfrak{sl}}_n)&\longrightarrow \iH(\bfk \text{Cy}_n)
		\\
		B_j &\mapsto \frac{-1}{q-1}[S_{j}],
		\quad  \K_j \mapsto [K_{S_j}], \text{ for }0\leq j\leq n-1. 
	\end{align}
	
	Define 
	\begin{align}\label{the map Psi A}
		\Omega_{\text{Cy}_n}:=\widetilde{\psi}_{\text{Cy}_n}\circ \Phi:{}^{\text{Dr}}\tUi_{\sqq}(\widehat{\mathfrak{sl}}_n)\longrightarrow \iH(\bfk \text{Cy}_n).
	\end{align}
	Then $\Omega_{\text{Cy}_n}$ sends
	\begin{align}
		\label{eq:HaDrA1}\K_j\mapsto [K_{S_j}],\qquad C\mapsto [K_\de],\qquad
		B_{j,0}\mapsto \frac{-1}{q-1}[S_j], \qquad \forall \,\,1\leq j\leq n-1.
	\end{align}
	For any $1\leq j\leq n-1$, $l\in\Z$ and $r\geq 1$, we define
	\begin{align}
		\label{def:haB}
		\haB_{j,l}:=(1-q)\Omega_{\text{Cy}_n}(B_{j,l}),\qquad
		\widehat{\Theta}_{j,r}:= \Omega_{\text{Cy}_n}(\Theta_{j,r}),\qquad
		\widehat{H}_{j,r}:= \Omega_{\text{Cy}_n}(H_{j,r}).
	\end{align}
	In particular, $\haB_{j,0}=[S_j]$ for any $1\leq j\leq n-1 $.

	Let $\X$ be a weighted projective line of weight type $(\bp,\ul{\bla})$.
	Recall that $\bla_i$ is the exceptional closed point of $\X$ of weight $p_i$ for any $1\leq i\leq \bt$. Recall that $\scrt_{\bla_i}$ is the Serre subcategory of $\coh(\X)$ consisting of torsion sheaves supported at $\bla_i$, and we have an equivalence $\scrt_{\bla_i}\cong\rep_\bfk(\text{Cy}_{p_i})$, which induces  an embedding of $\imath$Hall algebras:
	\begin{equation}
		\label{eq:embeddingx}
		\iota_i: \iH(\bfk \text{Cy}_{p_i})\longrightarrow\iH(\X_\bfk).
	\end{equation}
	
	Inspired by \eqref{def:haB}, we define
	\begin{align}
		\label{def:haBThH}
		\haB_{[i,j],l}:= \iota_i (\haB_{j,l}),\quad \widehat{\Theta}_{[i,j],r}:=\iota_i(\widehat{\Theta}_{j,r}),\quad \widehat{H}_{[i,j],r}:=\iota_{i}(\widehat{H}_{j,r}),
	\end{align}
	for any $1\leq j\leq p_i-1$, $l\in\Z$ and $r>0$.
	
	Let $\mathcal{C}$ be the Serre subcategory of $\coh(\X)$ generated by those simple sheaves $S$ satisfying $\Hom(\co, S)=0$. Then the Serre quotient $\coh(\X)/\mathcal{C}$ is equivalent to the category $\coh(\P^1)$ and the canonical functor $\coh(\X)\rightarrow\coh(\X)/\mathcal{C}$ has an exact fully faithful right adjoint functor
	\begin{equation}
		\label{the embedding functor F}
		\mathbb{F}_{\X,\P^1}: \coh(\P^1)\rightarrow\coh(\X),
	\end{equation}
	which sends $$\co_{\P^1}(l)\mapsto\co(l\vec{c}),\quad S_{\bla_i}^{(r)}\mapsto S_{i,0}^{(rp_i)},\quad S_{x}^{(r)}\mapsto S_{x}^{(r)}$$ 
	for any $l\in\Z, r\geq 1$, $1\leq i\leq \bt$ and $x\in\PL\setminus\{\bla_1,\cdots, \bla_\bt\}$. Then $\mathbb{F}_{\X,\P^1}$ induces an exact fully faithful functor $\cc_1(\coh(\P^1))\rightarrow \cc_1(\coh(\X))$, which is also denoted by $\mathbb{F}_{\X,\P^1}$. This functor $\mathbb{F}_{\X,\P^1}$ induces a canonical embedding of Hall algebras $\ch\big(\cc_1(\coh(\P^1))\big)\rightarrow \ch\big(\cc_1(\coh(\X))\big)$, and then an embedding
	\begin{equation}
		\label{the embedding functor F on algebra}
		F_{\X,\P^1}:\iH(\P^1_\bfk) \longrightarrow \iH(\X_\bfk).
	\end{equation}

	Inspired by \eqref{eq:hTm}, we define
	\begin{align}
		\label{def:Theta star}
		\widehat{\Theta}_{\star,m}:&= F_{\X,\P^1}(\haT_m)=\frac{1}{(q-1)^2\sqq^{m-1}}\sum_{0\neq f:\co(\vec{c})\rightarrow \co(m\vec{c}) } [\coker f],
		\\
		\widehat{H}_{\star,m}:&= F_{\X,\P^1}(\widehat{H}_m).
	\end{align}
	
	Recall the star-shaped graph $\Gamma=T_{p_1,\dots,p_\bt}$ in \eqref{star-shaped}. We define \begin{align}
		\delta:=\widehat{\co(\vec{c})}-\widehat{\co}\in K_0(\coh(\X));
	\end{align}
	compare with \eqref{eq:delta}. 
	
	\begin{theorem}[\text{\cite{LR21}}]
		\label{thm:morphi}
		For any star-shaped graph $\Gamma$, let $C$ be the GCM and $\X$ be the weighted projective line associated to $\Gamma$. Then there exists a $\Q(\sqq)$-algebra homomorphism
		\begin{align}
			\Omega: {}^{\text{Dr}}\tUi_{ \sqq}\longrightarrow \iH(\X_\bfk),
		\end{align}
		which sends
		\begin{align}
			\label{eq:mor1}
			&\K_{\star}\mapsto [K_{\co}], \qquad \K_{[i,j]}\mapsto [K_{S_{ij}}], \qquad C\mapsto [K_\de];&
			\\
			\label{eq:mor2}
			&{B_{\star,l}\mapsto \frac{-1}{q-1}[\co(l\vec{c})]},\qquad\Theta_{\star,r} \mapsto {\widehat{\Theta}_{\star,r}}, \qquad H_{\star,r} \mapsto {\widehat{H}_{\star,r}};\\
			\label{eq:mor3}
			&\y_{[i,j],l}\mapsto {\frac{-1}{q-1}}\haB_{[i,j],l}, \quad \Theta_{[i,j],r}\mapsto\widehat{\Theta}_{[i,j],r}, \quad H_{[i,j],r}\mapsto \widehat{H}_{[i,j],r}
		\end{align}
		for any $[i,j]\in\II-\{\star\}$, $l\in\Z$, $r>0$.
		Moreover, if $\Gamma$ is of finite type or affine type, then $\Omega$ is an embedding.
	\end{theorem}
	
	\begin{remark}
		Assume that $\X$ is of weight type $\bp=(p_1,p_2)$. Let
		$T=\bigoplus_{0\leq \vec{x}\leq \vec{c}} \co_\X(\vec{x})$ be the canonical tilting sheaf in $\coh(\X)$; see \cite{GL87}.
		Then $\End_\X(T)=\bfk Q$ is the canonical algebra with $Q=C_{p_1,p_2}$ as shown below. 
		\vspace{-6mm}
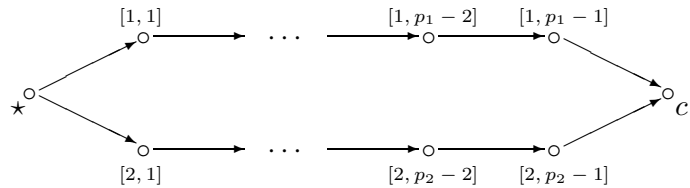
		\begin{center}\setlength{\unitlength}{0.75mm}
			\begin{equation*}
				\begin{picture}(80,20)(0,20)
					\put(-20,20){$\circ$}
					\put(-17,22){\vector(2,1){17}}
					\put(-17,21){\vector(2,-1){17}}
					\put(-22,17.5){$\star$}
					\put(0,10){$\circ$}
					\put(0,30){$\circ$}
					\put(50,10){$\circ$}
					\put(50,30){$\circ$}
					\put(72,10){$\circ$}
					\put(72,30){$\circ$}
					\put(92,20){$\circ$}
					\put(-3,6){\tiny$[2,1]$}
					\put(-3,34){\tiny${[1,1]}$}
					\put(44,6){\tiny $[2,p_2-2]$}
					\put(44,34){\tiny $[1,p_1-2]$}
					\put(67,6){\tiny $[2,p_2-1]$}
					\put(67,34){\tiny $[1,p_1-1]$}
					\put(94.5,17.5){$c$}
					
					\put(3,11.5){\vector(1,0){16}}
					\put(3,31.5){\vector(1,0){16}}
					\put(23,10){$\cdots$}
					\put(23,30){$\cdots$}
					\put(33.5,11.5){\vector(1,0){16}}
					\put(33.5,31.5){\vector(1,0){16}}
					\put(53,11.5){\vector(1,0){18.5}}
					\put(53,31.5){\vector(1,0){18.5}}
					
					\put(75,12){\vector(2,1){17}}
					\put(75,31){\vector(2,-1){17}}

				\end{picture}
				%
				%
			\end{equation*}
		\end{center}
		\vspace{1.1cm}
		Then the corresponding result in Theorem \ref{main thm2} also holds for $\X$ in this case.
	\end{remark}

	\section{Discussions and open problems}
	\label{sec:open}

	\subsection{Categorical symmetric pairs}
	
	Quantum group $\U$ and Drinfeld double $\tU$ are Hopf algebras. A categorical intepretation of the comultiplication was provided by Green \cite{Gr95} (which is extended to Bridgeland's setup \cite{Br}). 
	
	We have a quantum symmetric pair $(\tU, \tUi)$, where $\tUi$ is a coideal subalgebra of $\tU$, i.e., $\Delta: \tUi \rightarrow \tUi \otimes \tU$. 
	We have seen Hall algebra constructions of $\tU$ and $\tUi$, separately. 
	
	\begin{problem}
		Find a categorical framework for the coideal subalgebra structure of $\tUi$ in $\tU$ and for the quantum symmetric pair $(\tU, \tUi)$.  
	\end{problem}
	
	See \cite{FL21, FL+20} for a geometric interpretation of the quantum symmetric pair of (affine) type AIII.

	\subsection{$\imath$Hall basis}
	
	For Dynkin quivers, a suitable Hall basis (up to the normalization of some specific $v$-powers) for Ringel-Hall algebra is a PBW basis; cf. \cite[Chapter 11]{DDPW}.

	\begin{problem}
		Let $(Q, \btau)$ be a Dynkin $\imath$quiver. Identify an $\imath$Hall basis in the $\imath$Hall algebra $\tMHk$ explicitly in terms of $\imath$quantum group $\tUi$, under the isomorphism $ \tUi_{\sqq} \stackrel{\cong}{\rightarrow} \tMHk$ in Theorem~\ref{thm:main}.
	\end{problem}

	\subsection{Beyond quasi-split type}
	
	A Satake diagram consists of a pair $(\I =\I_\bullet \cup \I_\circ, \tau)$, where $\I =\I_\bullet \cup \I_\circ$ is a bi-colored partition of the Dynkin diagram $\I$, and $\tau$ a diagram involution (we allow $\tau=\Id$), subject to suitable conditions. In finite type, the Satake diagrams  are classified by Araki; cf. \cite[Table 4]{BW18b}; they are in bijection with real forms of complex simple Lie algebras.
	
	Associated to a general Satake diagram of Kac-Moody type, the $\imath$quantum groups $\Ui_\bvs$ were defined in \cite{Ko14} and the universal $\imath$quantum groups $\tUi$ in \cite{LW22, WZ22}, respectively.
	
	\begin{problem}
		Provide an $\imath$Hall algebra realization of $\tUi$ associated to a general Satake diagram.
	\end{problem}

	\subsection{Generalizing $\imath$quiver algebras}
	
	So far the $\imath$Hall algebras can only be used to realize quasi-split (i.e., $\I_\bullet =\emptyset$) $\imath$quantum groups under the additional assumption on the evenness of Cartan integers $c_{i,\tau i}$. (This assumption arises from the $\imath$quiver algebra construction in \S\ref{sec:iQA}.)
	
	This simplest example of quasi-split Satake diagrams which do not satisfy the evenness condition on $c_{i,\tau i}$ is the $A_2$ Dynkin diagram with nontrivial involution $\tau$. 
	Note that as $\tau$ does not preserve the orientation of the $A_2$ quiver, this is not an example of $\imath$quiver in the sense of \S\ref{subsec:i-quivers}.

	\begin{problem}
		Extend the formulation of $\imath$quiver algebra and the construction of $\imath$Hall algebra to realize $\tUi$ associated to $A_2$ with nontrivial involution $\tau$. 
	\end{problem}

	\subsection{Valued $\imath$quivers}
	
	In the realization of quantum groups via Hall algebras,	Ringel and Green actually used the Hall algebras of valued quivers to realise quantum groups associated to symmetrizable (which may not be symmetric) GCM. The reflection functors of valued quivers were also introduced, and they were used to realise Lusztig's symmetries of quantum groups; we refer to \cite{DDPW} for a good survey on Hall algebras of valued quivers. We expect a good outcome for the following problem. 
	
	
	\begin{problem}
		Formulate valued $\imath$quiver algebras and their homological properties. Use the $\imath$Hall algebras of valued $\imath$quiver algebras to realize the $\imath$quantum groups associated to symmetrizable GCM.		
	\end{problem}

	\subsection{Geometric $\imath$Hall algebras}
	
	In \cite{LW21b}, we use (singular) Nakajima-Keller-Scherozke varieties \cite{KS16,S19} to realize	the $\imath$quiver algebras, and then construct convolution algebras on the dual graded Grothendieck rings of perverse sheaves to realize the quasi-split $\imath$quantum group of ADE type. As a result, the dual basis of perverse sheaves gives a basis of $\imath$quantum groups, called dual $\imath$canonical basis, which has positive integral structure constants. This is an $\imath$generalization of Qin's geometric realization \cite{Qin16} of Drinfeld double quantum group $\tU$ and its dual canonical basis. 
	
	\begin{problem}
		Provide an algebraic construction of the dual canonical basis for $\tU$ and $\tUi$, respectively. (The rank one cases including $\tU (\mathfrak{sl}_2)$ are already open.)
	\end{problem}
	
	A precise relation (if any) between the $\imath$canonical basis on the modified $\imath$quantum group $\dot \Ui$ \cite{BW18b, BW21} and the dual canonical basis on the universal $\imath$quantum group $\tUi$ \cite{LW21b} remains to be clarified. 

	\subsection{Elliptic $\imath$Hall algebra}
	
	The $\imath$constructions are often related to various constructions of classical types, sometimes unexpectedly; see, e.g., \cite{BW18a} for connections to Kazhdan-Lusztig theory, \cite{BKLW, FL+20} for connections to convolution algebras arising from (affine) flag varieties, and \cite{LRW21} for connections to type C universal characters. 
	
	The definition of $\imath$Hall algebra of the projective line can be extended to higher genus curves. The $\imath$Hall algebra of an elliptic curve will be called the {\em elliptic $\imath$Hall algebra}. Stimulated by work of Schiffmann and Vasserot (see \cite{SV11}), we post the following loosely phrased question.
	
	\begin{problem}
		Is the elliptic $\imath$Hall algebra related to any of the following subjects? (i) spherical DAHA of type C; (ii)  (inhomogeneous) $\imath$Macdonald polynomials. 
	\end{problem}


\end{document}